\documentclass[a4paper,12pt]{article}
\usepackage{hyperref}
\usepackage[english]{babel}
\usepackage[T1]{fontenc}
\usepackage[utf8]{inputenc}
\usepackage{amsfonts}
\usepackage{amssymb}
\usepackage{amsmath}
\usepackage[separate-uncertainty=true]{siunitx}
\usepackage{graphicx}
\usepackage{cite}
\usepackage{epsfig}
\usepackage{psfrag}
\usepackage{tikz}
\usetikzlibrary{matrix,arrows,decorations.pathmorphing,positioning} 
\usepackage{url}
\usepackage{amsthm}
\usepackage{enumerate}
\usepackage{tikz-cd}

\usepackage{float}
\usepackage{appendix}
\usepackage{pgfplots}
\pgfplotsset{compat=1.15}
\usepackage{mathrsfs}
\usetikzlibrary{arrows}
\usepackage[margin=0.5in]{geometry}

\date{}

\title{\bf A Cosine Rule-Based Discrete Sectional Curvature for Graphs} 

\author{J.F. Du Plessis$^{1, }$\footnote{\url{23787295@sun.ac.za}}  {}  and  Xerxes D. Arsiwalla$^{2, 3,  }$\footnote{\url{x.d.arsiwalla@gmail.com}}\\  
{}  \\
{\it \small $^{1}$Stellenbosch University, South Africa}\\ 
{\it \small $^{2}$Pompeu Fabra University, Barcelona, Spain}\\
{\it \small $^{3}$Wolfram Research, USA}    
}

\begin{document}
\maketitle

\begin{abstract}
How does one generalize differential geometric constructs such as curvature of a manifold to the discrete world of graphs and other combinatorial structures? This problem carries significant importance for analyzing models of discrete spacetime in quantum gravity; inferring network geometry in network science; and manifold learning in data science. The key contribution of this paper is to introduce and validate a new estimator of discrete  sectional curvature for random graphs with low metric-distortion. The latter are constructed via a specific graph sprinkling method on different manifolds with constant sectional curvature. We define a notion of metric distortion, which quantifies how well the graph metric approximates the metric of the underlying manifold. We show how graph sprinkling algorithms can be refined to produce hard annulus random geometric graphs with minimal metric distortion. We construct random geometric graphs for spheres,  hyperbolic and euclidean planes; upon which we validate our curvature estimator. Numerical analysis reveals that the error of the estimated curvature diminishes as the mean metric distortion goes to zero, thus demonstrating convergence of the estimate. We also perform comparisons to other existing discrete curvature measures. Finally, we demonstrate two practical applications: (i) estimation of the earth's radius using geographical data; and (ii) sectional curvature distributions of self-similar fractals. 
\end{abstract}

\clearpage

\tableofcontents


\section{Introduction}

As one of the oldest branches of mathematics, pre-17$^{th}$ century geometry largely developed as the study of spatial properties of mostly continuous structures. The advent of real analysis and infinitesimal calculus in the 17$^{th}$ century by Newton and Leibniz (and, Descartes and Fermat before them) led to analytic approaches to geometry, culminating in contemporary differential geometry. It was only in the 19$^{th}$ century that non-Euclidean geometric structures were studied following the works of Riemann and Poincare. But once again, these fields have heavily relied on analysis and differential calculus (though more recently, these frameworks have been expressed in abstract algebraic and category-theoretic foundations). On the other hand, and in contrast to "continuous mathematics", discrete structures, including integers, graphs, countable sets, symbolic systems, etc have been studied using discrete analysis and other combinatorial methods, under the realm of "discrete mathematics". Needless to say, methods from discrete analysis are helping shape contemporary foundations in areas as diverse as theoretical physics, formal logic, computer science, network science, data science, machine learning, operations research and other areas of applied mathematics \cite{rosen2012discrete}. From that perspective, it is natural to ask whether can one carry over geometric concepts from the continuous domain to discrete combinatorial structures. In particular, how can one generalize classical differential geometric notions of curvature, torsion, tangent spaces, fiber bundles, etc to work on graphs, hypergraphs and other discrete systems? Addressing these questions will inevitably have far-reaching consequences across disciplines from quantum gravity to network and data science. In this paper we propose and validate a potential discrete candidate for sectional curvature of graphs.  

As of now, there exist several nonequivalent definitions of discrete curvature (and their respective variations) that have been developed for graphs and have been applied to a variety of problems involving network models  \cite{lin_ricci_2011, jost_olliviers_2014, samal_comparative_2018}. Prominent examples of these include:  `Ollivier-Ricci curvature'   \cite{ollivier_ricci_2007,saucan_discrete_2019,lin_ricci_2011},   `Forman-Ricci curvature'    \cite{forman_bochners_2003}, `combinatorial curvature' \cite{kamtue2018combinatorial},  `resistance curvature'  \cite{devriendt_discrete_2022}, `Wolfram-Ricci curvature' \cite{Wolfram2020, gorard2020some} and a modified Ollivier curvature for mesoscopic graph neighbourhoods \cite{PhysRevResearch.3.013211}, among others. Some of these definitions  have also been extended to the case of hypergraphs  \cite{eidi_ollivier_2020}. Several recent applications investigating geometric properties of real-world complex networks have been reported and have been shown to correlate well with standard network theoretic measures in a wide range of classes of graphs and hypergraphs  \cite{sreejith_forman_2016, samal_comparative_2018}.  These notions of discrete curvature and discrete differential geometry are claimed to be particularly relevant to the foundations of physics and have been employed in approaches to quantum gravity \cite{trugenberger_combinatorial_2017, tee2021enhanced}. The aforementioned proposals seek to specify a curvature for each vertex and/or edge in the network. However, many of these definitions (if not all) do not sense mesoscopic or global features of the network geometry \cite{PhysRevResearch.3.013211}.  This may consequently lead to inconsistent curvature estimates on networks such as random geometric graphs \cite{penrose_random_2003}, constructed from well behaved manifolds. On such well-behaved manifolds one would hope any curvature estimator of a graph with metric structure close to that of the manifold, would give a curvature estimate close to that of the manifold. This is part of the challenge we seek to address with our sectional curvature estimator. In the context of random geometric graphs, a noteworthy recent development is the above-mentioned modified Ollivier-Ricci curvature, which was shown to converge to the Ricci curvature on random geometric graphs in the limit of very large vertex count \cite{PhysRevResearch.3.013211}. This is the only known case where discrete curvature has been tested in such a near-continuum limit. More generally, many of the existing consistency tests of discrete curvature do not fully take into account metric properties of the graph in their construction, particularly, metric embedding and distortion compared to metrics of some underlying manifold, upon which these graphs are obtained from sprinkling, and should therefore also have corresponding curvature estimates. Hence, an important step we undertake in this direction is the systematic construction of hard annulus random geometric graphs \cite{dettmann_random_2016} with low metric-distortion. This is a crucial development for testing and comparing notions of discrete curvature, since it provides a notion of how 'close' a graph is to some manifold with known curvature, which we can then use to judge whether graphs that are close to a manifold with known curvature also gives estimates close to that known curvature. Although the estimator developed here is applicable to any path metric space, we use this new and more rigorous test to show the accuracy of the new estimator in specific situations. We can then also more easily compare the estimate to existing notions. This newly developed test only gives a notion of how well the estimator captures the geometric information of the graph, and is independent of the purely network-theoretic context many existing notions of discrete curvature were developed for. It then gives an idea of how much these network theoretic curvature constructions carry over to geometric problems, and how they compare in such problems with curvature constructions designed exactly for such geometric applications such as the discrete sectional curvature estimate defined in this paper.


The main objective of this paper is to introduce and test a new estimator of sectional curvature, which we validate on random geometric graphs with low metric-distortion. The latter are graphs that mimic the metric structure of continuous manifolds well, and hence enable a rigorous assessment  of the accuracy of our curvature estimates. More specifically, we  construct the so-called "hard annulus random geometric graphs", which are slightly generalized in the sense that they refer to manifolds of given topology and curvature (instead of the usual euclidean setting where random geometric graphs are often  considered), which uses parameters that can be tuned to find graphs with metric structure close to that of the continuous manifold it is supposed to represent. This serves as a proof-of-concept for this new curvature estimator as these graphs are obtained from manifolds with known geometric properties. 

In Riemannian geometry, sectional curvature can be seen as the Gaussian curvature of geodesic planes in some manifold with dimension two or higher. Knowing the sectional curvature completely, also determines the Riemann curvature tensor completely (and vice-versa) meaning it gives an equivalent notion of curvature. From it the Ricci tensor as well as the Ricci scalar curvature can be calculated by an averaging process. If we want to define some form of curvature for structures more general than manifolds, then classical definitions of sectional curvature may not apply. If we consider a space with given minimum length, taking the limit of this length scale to zero is in general not practical.   Taking the limit as the length goes to the minimum length will also fail, as there is in general not enough geometric information left. In these discrete spaces the curvature information is in a sense "smeared" over the structure between many vertices, instead of being well-defined at each point. This can be seen from the fact that manifolds are locally euclidean, so in a fine enough discrete approximation of a manifold, a single vertex or edge is unable to distinguish if it is in a curved space or not, let alone how curved it is. This can be seen in difficulties in showing that edge or vertex based curvatures act as expected in large random graphs with non-zero curvature. If we assume that our discrete space is a path metric space (the distance between two points is the length of the shortest path between them, where a path is some function that parametrizes the path as a function of distance) we can exploit the spherical cosine rule to get an estimate of curvature. If this discrete structure approximates the metric of some manifold, we want this estimate to correspond to the curvature of the manifold, up to some error which is proportional to the "error" between the discrete space and the manifold. This error will properly defined in what follows.

The outline of this paper is as follows: In section \ref{sec:curv} we present a new estimator for discrete sectional curvature, originating from a cosine rule. We discuss the applicability of this estimator on discrete constructions, in particular, random geometric graphs. Taking inspiration from metric embedding theory, a notion of metric distortion is defined, in order to quantify how far a graph's metric deviates from the metric of the manifold it is supposed to represent. We also detail a sprinkling algorithm for generating random geometric graphs with given geometric properties and low metric-distortion. In section \ref{sec:results} we extensively demonstrate how our estimator works on random geometric graphs generated from manifolds of positive, negative and zero curvature; including a real-world example, where our method can be used to estimate the radius of the earth from geographical data. In section \ref{sec:convergence} we discuss how the error of the estimated curvature diminishes as the mean metric distortion goes to zero, thus demonstrating convergence of the estimate in the continuum limit. In section  \ref{sec:vertex} we use our estimator to provide a vertex-based sectional curvature estimate on graphs. In section \ref{sec:compare} we compare this discrete sectional curvature to an existing notion of discrete curvature that uses large scale structure. In section \ref{sec:apps} we demonstrate other applications of our methods to data that may not have underlying constant sectional curvature.  In section \ref{sec:conclusion} we conclude with key results and next steps.

\section{A New Estimator for Discrete Sectional Curvature}\label{sec:curv}

\subsection{Sectional Curvature }
\subsubsection{On (pseudo)Riemannian Manifolds}\label{sec:Riem}
On (pseudo)Riemannian manifolds the Riemann curvature tensor is generally the most general curvature-related object of interest, since both the Ricci curvature tensor and the Ricci scalar curvature can be calculated from it. However, there is a simpler notion of curvature, namely sectional curvature, that turns out to give equivalent curvature information as the Riemann curvature tensor. The sectional curvature can geometrically be thought of the Gaussian curvature of a given geodesic plane. The sectional curvature $K(u,v)$ of a geodesic plane defined by two, non-parallel, vectors $u,v$ in a Riemannian manifold can be defined in terms of the Riemann tensor as
\begin{equation}
    K(u,v)=\frac{\langle \text{Riem}(u,v)v,u\rangle}{u^2 v^2-\langle u, v\rangle^2}
\end{equation}
where $\langle u , v \rangle$ denotes the inner product of the Riemannaian manifold.  It can also be shown that knowing the sectional curvature at a point completely, determines the Riemann tensor at that point completely. So knowing the sectional curvature completely is theoretically equivalent to knowing the Riemann tensor completely, and therefore it should be possible to calculate the Ricci tensor and Ricci scalar from sectional curvatures. We can denote the average sectional curvature with one vector $v$ fixed 
\begin{equation}
    \kappa(v)=\langle K(u,v)\rangle_u
\end{equation} 
where $\langle\cdot \rangle_u$ is the average over all vectors $u$ of some length (their exact length is irrelevant here) in the tangent space, using the solid angle over the hypersphere. Then it turns out a projection of the Ricci tensor is given by
\begin{equation}
    \text{Ric}(v,v)=(n-1)\kappa(v) v^2
\end{equation}
which determines the Ricci tensor completely, since it is symmetric and thus diagonalizable. We can then go one step further and calculate the Ricci scalar by averaging out over all (unit vector) projections, giving us
\begin{equation}
    R= n(n-1)\kappa
\end{equation}
where\begin{equation}
    \kappa=\langle\kappa(v)\rangle_v
\end{equation}
This then allows us to compute an object of great significance in General Relativity; the Einstein tensor
\begin{equation}
    G(u,v)=\text{Ric}(u,v)-\frac12 R g(u,v)
\end{equation}
which we can see is trivially also symmetric, and therefore diagonalizable. We therefor only have to consider projections
\begin{equation}
    G(v,v)=(n-1)\kappa(v) v^2-\frac12 n(n-1)\kappa v^2
\end{equation}
in order to determine it completely. It can be shown that such a projection of the Einstein tensor is proportional to the average sectional curvature of planes orthogonal to the projection vector.
 
Hence, the sectional curvature is a very useful quantity for determining other types of curvatures. Furthermore, what is interesting about the sectional curvature is that it can be fully determined from combinatorial data. This feature will allow us to construct an algorithm to compute a notion of a discrete sectional curvature associated to random graphs.

Consider a geodesic triangle (which therefore lies in a geodesic plane) with side lengths $a,b,c$ and the angle opposite side $c$ is $\gamma$. If the geodesic plane has constant sectional curvature over the area of the triangle, we can use some generalisations of the cosine rule. We know in the Euclidean plane (with sectional curvature $K=0$) we trivially have the law of cosines
 \begin{equation}
     c^2=a^2+b^2-2a b\cos(\gamma).
 \end{equation}
 On a sphere with radius $R$, and therefore sectional curvature $K=\frac{1}{R^2}$, we have the spherical law of cosines
 \begin{equation}
     \cos\left(\frac{c}{R}\right)=\cos\left(\frac{a}{R}\right)\cos(\frac{b}{R})+\sin\left(\frac{a}{R}\right)\sin\left(\frac{b}{R}\right)\cos(\gamma).
 \end{equation}
Similarly on a hyperbolic plane with constant sectional curvature $K=-\frac{1}{k^2}$, we have the hyperbolic law of cosines
\begin{equation}
     \cosh\left(\frac{c}{k}\right)=\cosh\left(\frac{a}{k}\right)\cosh(\frac{b}{k})+\sinh\left(\frac{a}{k}\right)\sinh\left(\frac{b}{k}\right)\cos(\gamma).
\end{equation}
 it is straight-forward to then show that all three of these cases is neatly summarised in the single equation
 \begin{equation}
     \cos(c\sqrt{K})=\cos(a\sqrt{K})\cos(b\sqrt{K})+\sin(a\sqrt{K})\sin(b\sqrt{K})\cos(\gamma)
 \end{equation}
where for $K=0$ a limit needs to be carefully taken to recover the usual law of cosines. This means that the three measurable quantities namely distance, angle and sectional curvature go together, in the sense that knowing two of the three (under the assumption of constant sectional curvature) theoretically allows one to calculate the other. This equation can of course be further simplified if we assume our triangle is a right-angled triangle with $c$ as the hypotenuse, since then $\cos(\gamma)=0$, giving
  \begin{equation}\label{eq:master}
     \cos(c\sqrt{K})=\cos(a\sqrt{K})\cos(b\sqrt{K}).
 \end{equation}
We then need to consider if a given combination of lengths $a,b,c$ (satisfying the triangle inequality) only has a single value of $K$ satisfying Equation \ref{eq:master}. It can be easily seen that $K=0$ is always a trivial solution to the equation, which means we will in general ignore the $K=0$ root\footnote{Since the $K=0$ root is a simple root when considering triangles in curved geometries, dividing the above equation by $K$ gives an equation which only has roots corresponding to the curvature. This expression is a bit more unwieldy, although depending on the root-finding algorithm it might be better to consider for certain applications}. Let us then consider the positive $K$-axis. There one can find infinitely many roots. Notice however, that a triangle with lengths $a,b,c$ can not exist on a sphere with radius smaller than $\frac{\max(a,b,c)}{\pi}$, and therefore the only possible physical roots to consider are those for  $K\leq\frac{\pi^2}{\max(a,b,c)^2}$. With these two constraints it can be argued (for details we refer the reader to Appendix \ref{apn:proof}) that there is only one other root, and since we know that the underlying curvature will always be a root, we then know that the single root satisfying the above conditions must correspond to the sectional curvature.  Therefore, if we have a right-angled triangle on a geodesic plane of constant sectional curvature, we can then calculate that sectional curvature from the edge lengths of that triangle. On the other hand, if we have a geodesic plane of varying sectional curvature, we can then consider triangles very small compared to the magnitude of the gradient of curvature, in order to get an estimate of the curvature in that region.

\subsubsection{On Path Metric Spaces}
Let us consider some path metric space, which we will define as a metric space where the distance between two points coincides with the length of the shortest path between them. There will in general be more than one shortest path, but we only concern ourselves with the fact that at least one exists. We can then take an arbitrary path triangle (three shortest paths between three points, forming a triangle), and form an (possibly only approximate) isosceles triangle from it (by 'cutting' the longest leg to be equal to another side). By dividing the edge not equal to the other two edges in half, we can declare the two triangles formed (as can be seen in Figure \ref{fig:tikzTriang}) to be right-angled triangles. On Riemannian surfaces with constant curvature this will be true, as illustrated in Figure \ref{fig:embedded}. We can use this triangle to calculate a sectional curvature as detailed above, and we know that in a Riemannian manifold with constant sectional curvature over the triangle this will give the correct result. However we can do this calculation in any path metric space and still get a result. We can interpret this result as an estimate of the sectional curvature of a Riemannian manifold whose metric structure is well-approximated by the path metric space. Errors of this estimate then comes from the rate of change of curvature over the triangles considered, how faithfully the described triangles can be constructed and how accurate the metric of the path metric space is to that of a Riemannian manifold (this idea of 'metric distortion' is made more concrete in Section \ref{sec:distortion}).

\definecolor{blue}{rgb}{0,0,1}
\definecolor{red}{rgb}{1,0,0}
\definecolor{orange}{rgb}{0.67,0.33,0}

\begin{figure}[H]
\centering
\begin{tikzpicture}[line cap=round,line join=round,>=triangle 45,x=1cm,y=1cm,scale=1.3]

\node at (3,0.15) {a};
\node at (5.15,1) {b};
\node at (5.15,-1) {b};
\node at (2.75,1.3) {c};
\node at (2.75,-1.3) {c};

\draw [line width=2pt] (0,0)-- (5,-2);
\draw [line width=2pt] (0,0)-- (5,2);
\draw [line width=2pt] (5,2)-- (5,0);
\draw [line width=2pt] (5,0)-- (5,-2);

\draw [line width=1pt] (2.525,-0.96) -- (2.48,-1.04);
\draw [line width=1pt] (2.48,1.04) -- (2.51,0.95);
\draw [line width=1pt] (5.04,1.02) -- (4.95,1.02);
\draw [line width=1pt] (5.04,0.97) -- (4.95,0.97);
\draw [line width=1pt] (5.04,-0.98) -- (4.95,-0.98);
\draw [line width=1pt] (5.05,-1.02) -- (4.95,-1.02);

\draw [line width=1pt,color=orange] (4.83,0)-- (4.83,0.17);
\draw [line width=1pt,color=orange] (4.83,0.17)-- (5,0.17);
\draw [line width=1pt,color=orange] (4.83,0)-- (4.83,-0.17);
\draw [line width=1pt,color=orange] (4.83,-0.17)-- (5,-0.17);

\draw [line width=2pt,color=blue] (0,0)-- (5,0);
\draw [fill=red] (0,0) circle (3pt);
\draw [fill=red] (5,-2) circle (3pt);
\draw [fill=red] (5,2) circle (3pt);
\draw [fill=blue] (5,0) circle (3pt);

\end{tikzpicture}
\caption{Illustration of constructing a right-angle in a region of constant sectional curvature}
\label{fig:tikzTriang}
\end{figure}
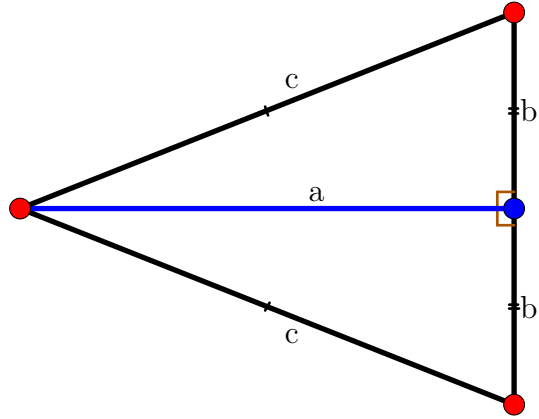

\begin{figure}[H]
    \centering
    \includegraphics[scale=0.095]{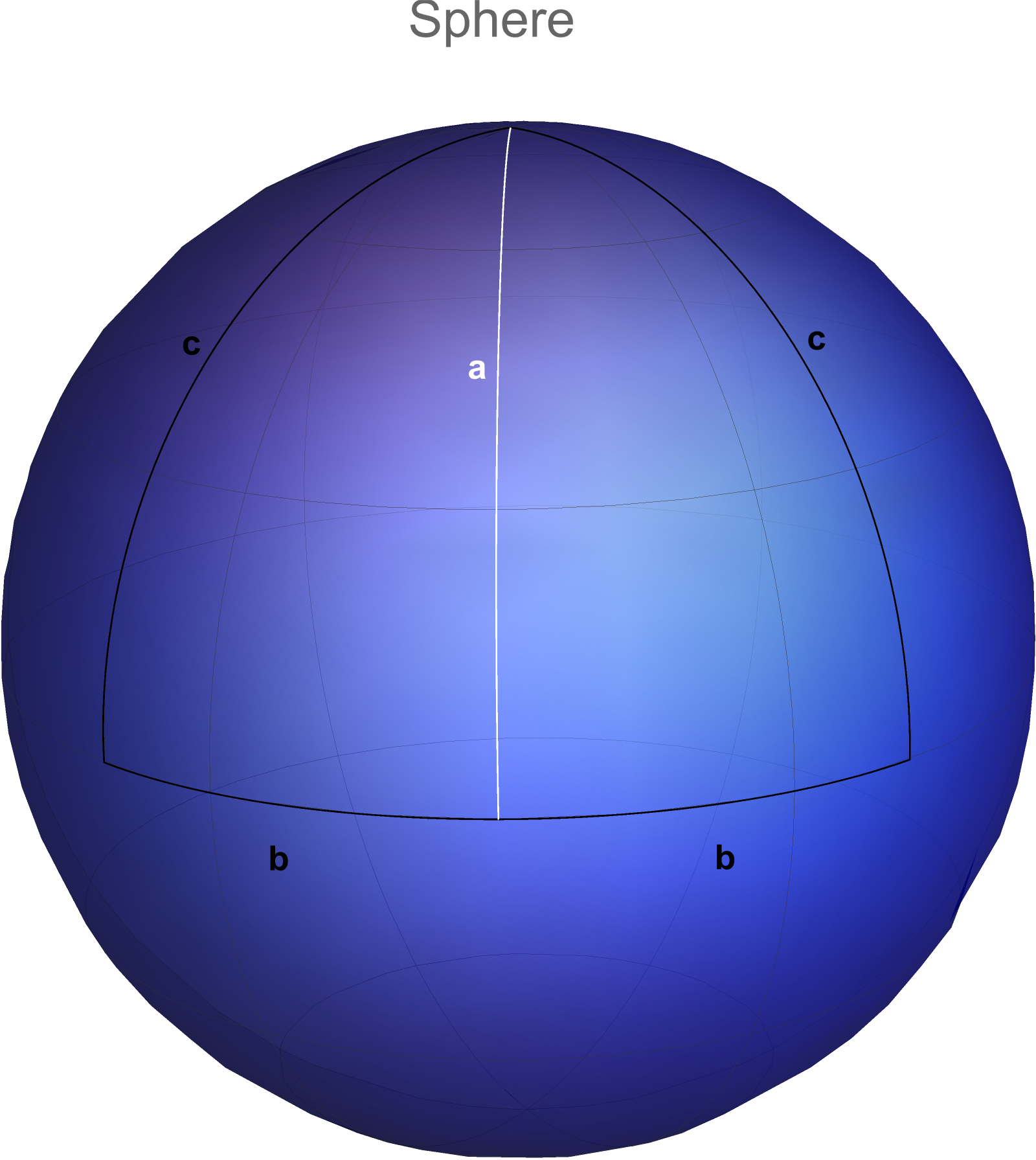}
    \includegraphics[scale=0.095]{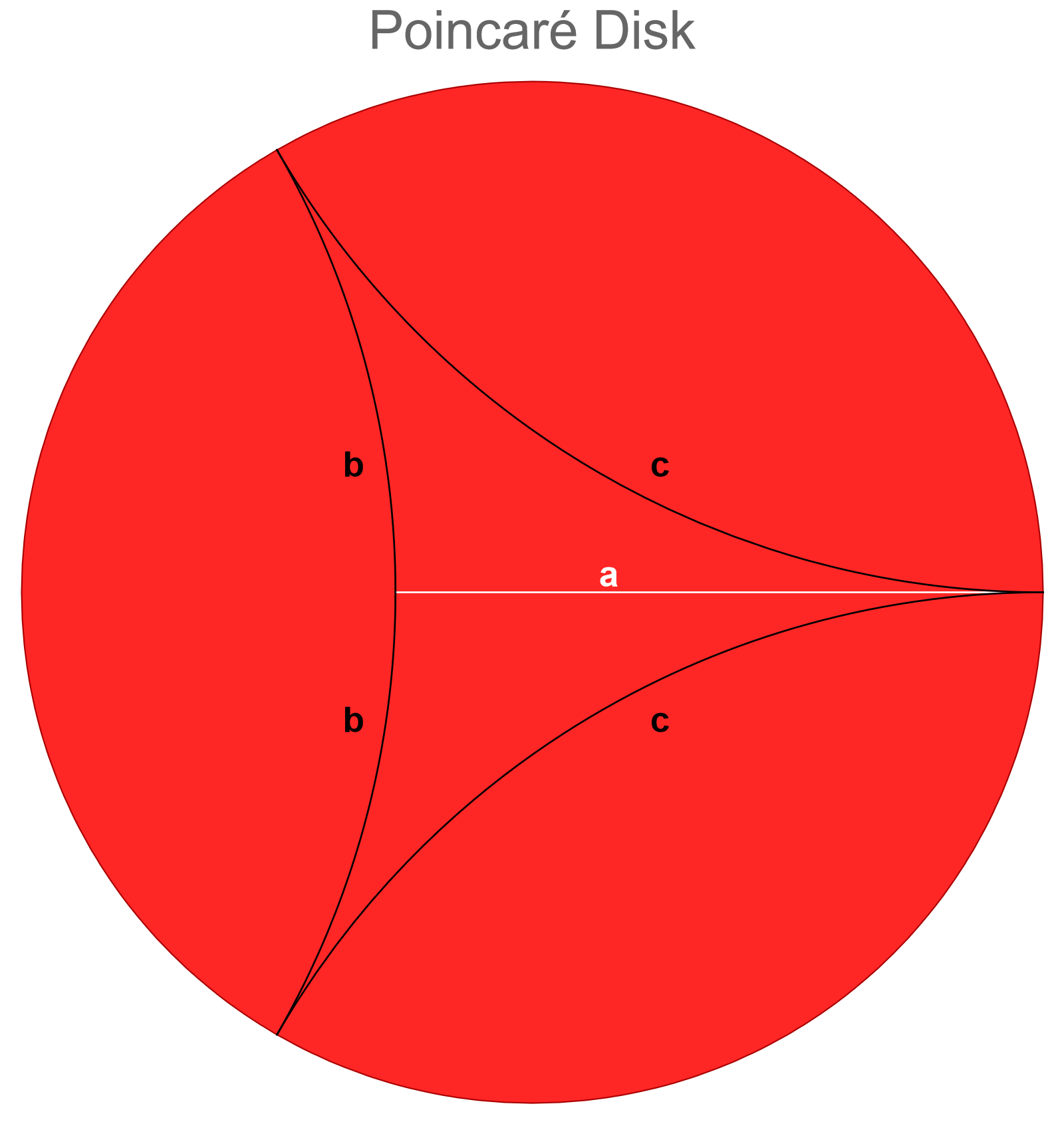}
    \includegraphics[scale=0.095]{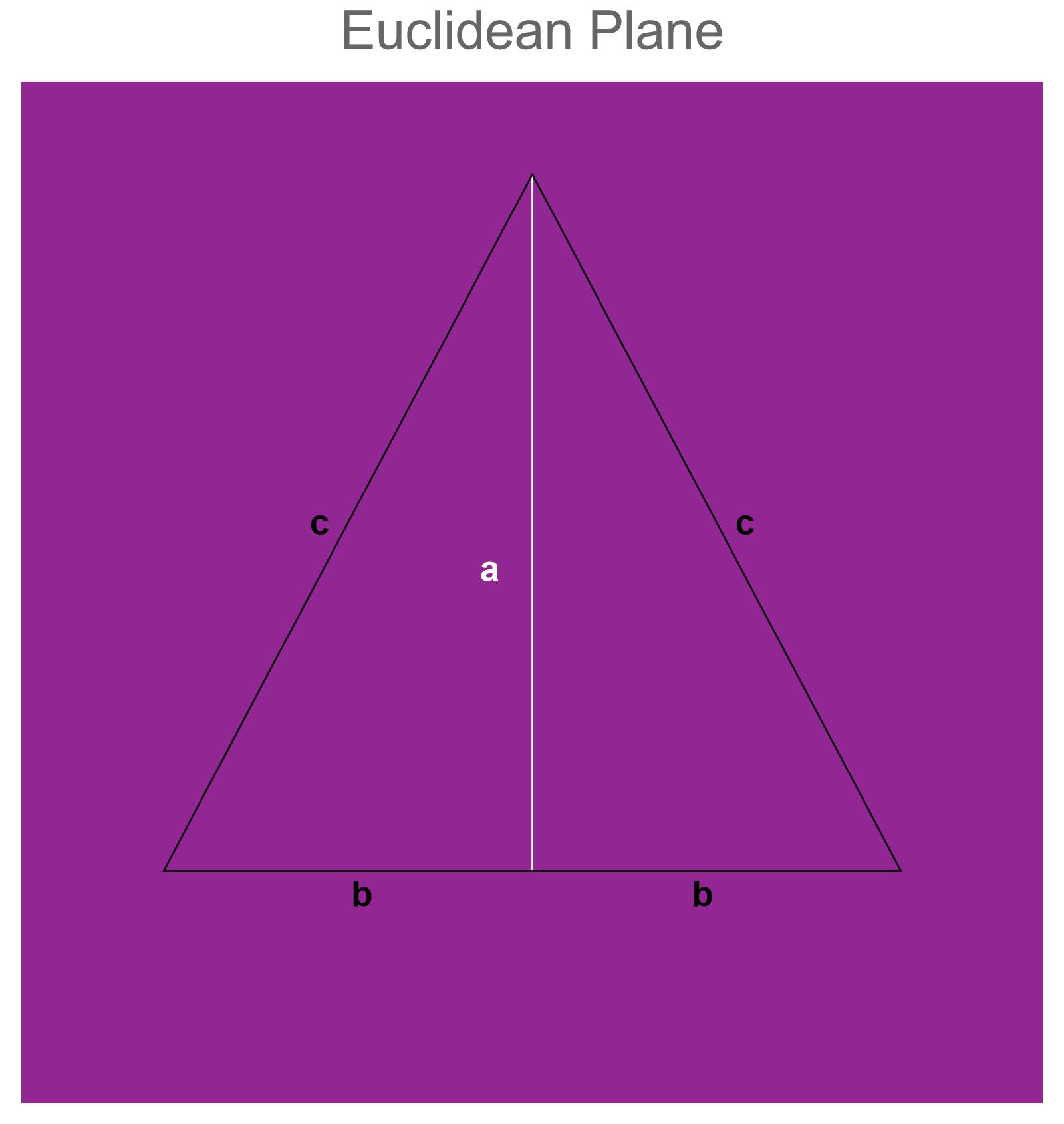}
    \caption{The triangle construction illustrated in Figure \ref{fig:tikzTriang} embedded into surfaces of constant positive, negative and 0 curvature respectively.}
    \label{fig:embedded}
\end{figure}
\begin{figure}[H]
\centering
\includegraphics[scale=0.25]{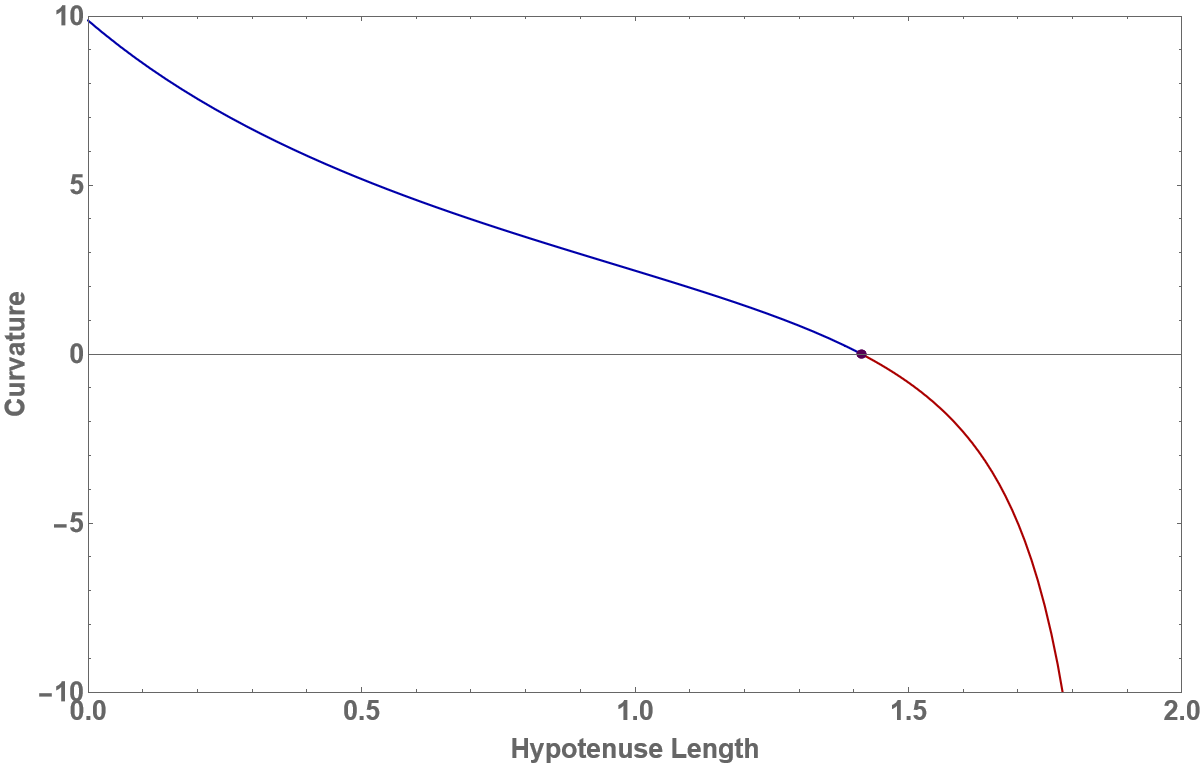}
\caption{The curvature over a right-angled triangle with two legs with length $1$ and the hypotenuse is varied on the x-axis. One can see that small errors in edge length measurements have a much larger effect in negative curvature geometries than in positive curvature geometries.}
\label{fig:curv}
\end{figure}

To get an idea of the robustness of the curvature to uncertainties in the lengths of the triangle (introduced by metric distortion, edge lengths and the possible failure to form perfect right-angled triangles) we can consider a the curvature calculated from a right-angled triangle with side lengths 1 and hypotenuse $c$. As can be seen from Figure \ref{fig:curv}, we can expect more robustness to noise for cases of positive curvature, due to small variations in length only introducing small (bounded) variation in curvature estimate.

\subsubsection{On Graphs}\label{sec:graphs}
We note here that since this discrete sectional curvature only uses distances, and is applicable to path metric spaces in general, we can use this as a definition of the curvature of a graph. It would be reasonable to expect that a graph that approximates the metric of some manifold well, and has effective edge lengths short relative to the curvature, should have sectional curvatures that correspond to that given manifold. It would be unreasonable to assume that in some chaotic graph any choice of triangle would give the expected results, and indeed in general it does not. Instead it is necessary to first choose some minimum length scale (since below this scale the possible values of the curvature has too large error margins). Next we need to average out over a sample of triangles in some region (sub-graph) of the graph. There is of course not a unique way to do this averaging, and two sensible possibilities comes to mind. The first, is to simply take some sample of curvature values and average out over them. This seems to be the most computationally viable approach, and is the one employed below. Another, as of yet unexplored, possibility is to first ask what the average triangles are (say by finding right-triangles, and then take the average hypotenuse for every right-triangle with each given side lengths), and then to compute the curvature of these averaged triangles. This is of course inefficient (since it takes many more operations on the graph, which are generally more expensive than finding roots), but should conceivably work with smaller minimum length scales. As a heuristic it seems that minimum length scales on the order of $\frac{1}{\sqrt{|K|}}$ gives reasonable results. This would correspond to the difficulty in reasonably quantifying the error in cases where we expect $0$ curvature, since there is no inherent length scale to compare our errors to. In the euclidean plane we can re-scale our metric arbitrarily to re-scale our error of curvature to whatever we want. This makes euclidean setting unreasonable for comparisons and testing. To test the robustness of discrete sectional curvature we need to construct graphs that have metrics (the combinatorial distance multiplied by some effective edge length) that are close to that of Riemannian manifolds with known curvatures. In order to have such a graph we first need to quantify what we mean by a metric approximating another metric well.

\subsection{Defining Metric-Distortion for Graphs}\label{sec:distortion}
In order to make judgements on the accuracy of the discrete sectional curvature described in Section \ref{sec:curv} we need some way to quantify the error some discrete space has in approximating a given metric space (Riemannian manifolds in our case). Taking inspiration from metric embedding theory\cite{abraham_advances_nodate} we can define a notion of metric distortion. As has been noted in \cite{chennuru_vankadara_measures_2018}, many existing measures have some undesirable characteristics, especially in terms of robustness in terms of outliers and noise. Following \cite{chennuru_vankadara_measures_2018} we can consider two metric spaces $G,X$ with metrics $d_G,d_X$ respectively. If we then have an embedding of $G$ into $X$ given by $f:G\rightarrow X$, we can define the ratios 
\begin{equation}
    \rho_f(u,v)=\frac{d_X(f(u),f(v))}{d_G(u,v)}\hspace{25pt}u,v\in G.
\end{equation}
which we will refer to as the embedding ratios. We will define $\rho_f(v,v)=1$. In contrast to what is done in \cite{chennuru_vankadara_measures_2018}, we want to recognise that the embedding ratios could as well have been defined as their inverses $[\rho_f(u,v)]^{-1}$. We would also like to recognise that re-scaling either of our metrics by some constant, leaves their structure unchanged, and we would like for it to leave our measure of distortion similarly unchanged. We can then think of the logs of the embedding ratios as forming some distribution, which is translated when a metric is re-scaled, and reflected over the x-axis when the inverse definitions $[\rho_f(u,v)]^{-1}$ of the log ratios are used. The natural quantities to define on this distribution that is invariant under those transformations is then the standard deviation, variance and mean deviation. For our purposes the mean deviation gives more consistent results, although the other options are also perfectly reasonable. We therefore define the distortion of an embedding to be the mean deviation of the logs of the embedding ratios, which for some finite space $G$ with $|G|$ points takes the form:
\begin{equation}
    \text{dist}(f)=\frac{1}{|G|^2}\sum_{u,v\in G}|\log(\rho_f(u,v))-\log(\overline{\rho_f}_\text{geom}) |,
\end{equation}
where $$\overline{\rho_f}_\text{geom}=\left[\prod_{u,v\in G}\rho_f(u,v)\right]^{|G|^{-2}}$$ is the geometric mean of the metric ratios. We can see then if we re-scale our metric $d'_G(u,v)=d_G(u,v)\times\overline{\rho_f}_\text{geom}$ to get the re-scaled embedding fractions $\rho_f'(u,v)=\frac{\rho_f(u,v)}{\overline{\rho_f}_\text{geom}}$ then the distortion reduces to 
\begin{equation}
    \text{dist}(f)=\frac{1}{|G|^2}\sum_{u,v\in G}|\log(\rho_f'(u,v)) |.
\end{equation}
If our space $G$ is a graph, with $d_G$ being the combinatorial metric on that graph (the metric that counts the least amount of edges one needs to traverse to get from one vertex to another), then we will call $l_e=\overline{\rho_f}_\text{geom}$ the effective edge length of the graph, with respect to the embedding $f$. For our purposes $f$ will be implicit, so it can be neglected without causing confusion.

\subsection{Generating Low Metric-Distortion Graphs}
If we want to study how the curvature defined in Section \ref{sec:curv} acts in discrete metric spaces, graphs are the natural candidate. Therefore it is necessary to obtain graphs of low distortion, to see how well the estimate of curvature acts, as a function of that distortion. With no know algorithm to find a graph of minimal distortion, some algorithm for generating graphs with low distortion is necessary. Naively some lattice or triangulation might seem to be the natural answer. Lattices however, will in general have poor metric distortion, and a key observation is that distortion doesn't decrease with graph size, giving no sense of convergence. This is due to the regularity of lattices, so distortions compound on longer length scales. Then next obvious option is triangulation. By using the Mathematica function   \textbf{MeshConnectivityGraph@DelaunayMesh}  to get a triangulation of some set of points in $R^2$. It turns out we can do even better than the already low distortion given by the triangulation by employing a graph sprinkling algorithm.

\begin{figure}[H]
    \centering
    \includegraphics[scale=0.25]{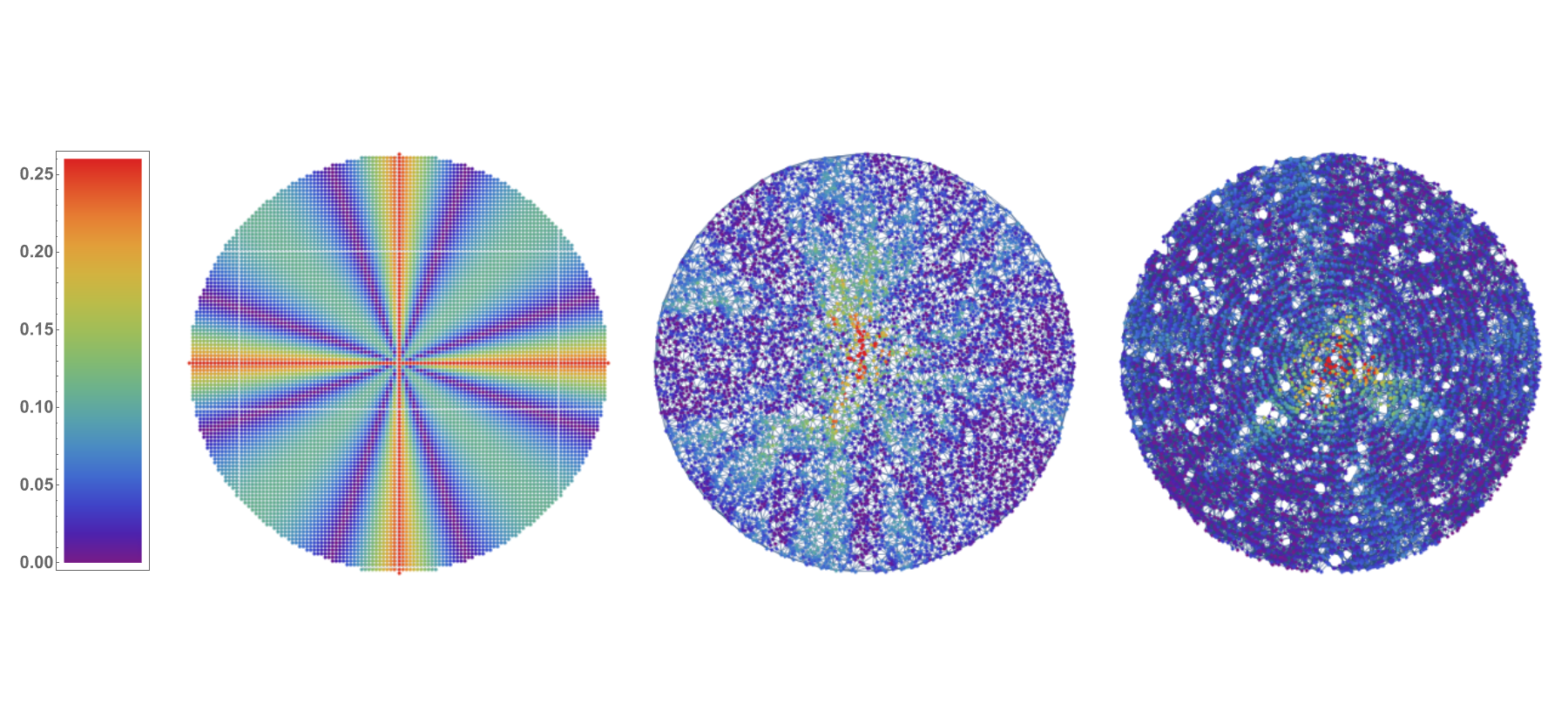}
    \caption{Visualising the distortions between the graph center of graphs representing the unit disk. From left to right we have a lattice, a triangulization and a random geometric graph of the unit disk in the euclidean plane. We can see the large distortions that is ubiquitous in lattice graphs with scaled combinatorial metrics, the smaller distortions (crucially the magnitude of the distortion seems to decrease with length scale) in a triangulization, and then the extremely small distortions seen in random geometric graphs (above some minimum length scale). The `holes' in the last graph are presumably due to the embedding algorithm, and are irrelevant for our purposes.}
    \label{fig:distVis}
\end{figure}

\begin{figure}[H]
    \centering
    \includegraphics[scale=0.2]{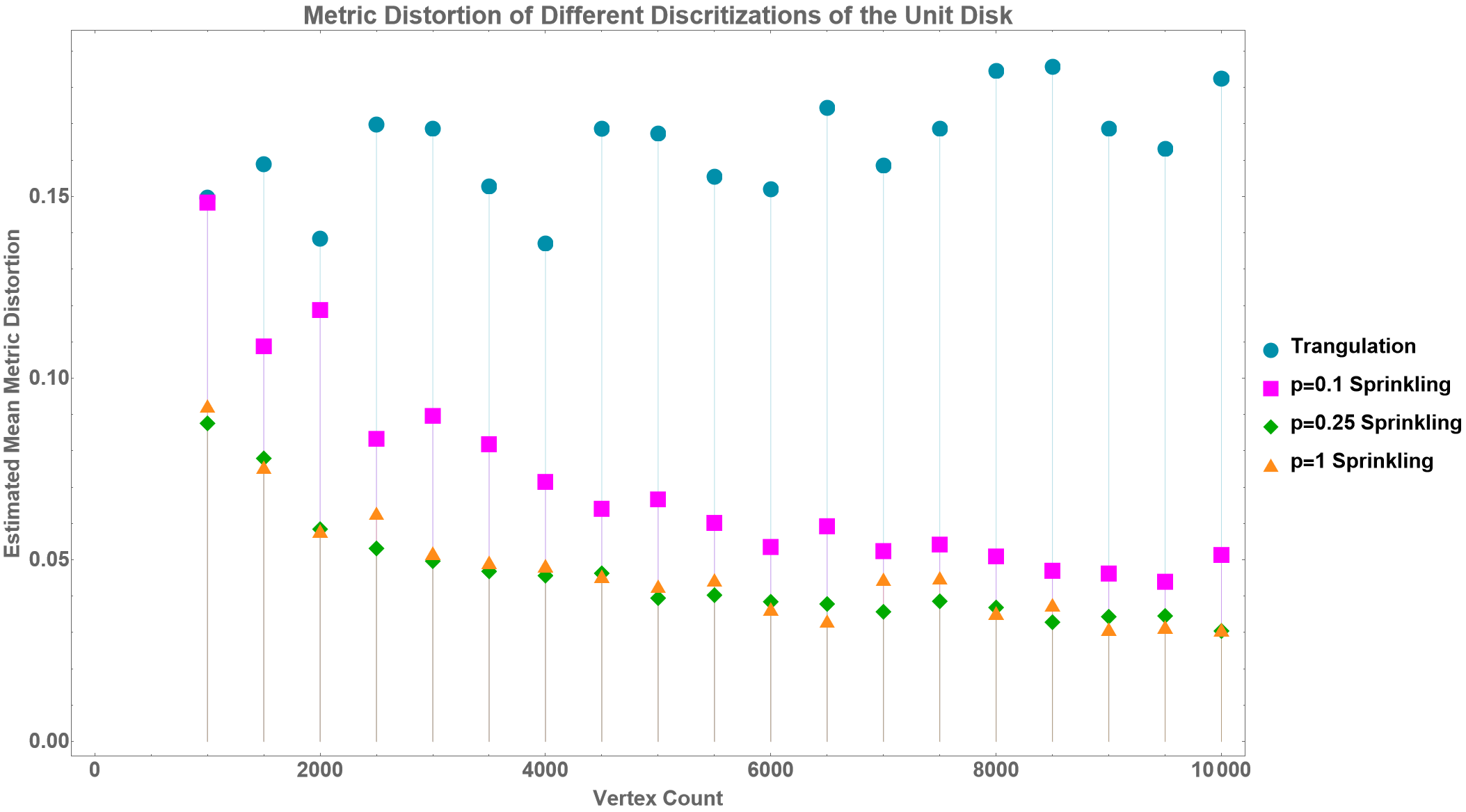}
    \caption{By generating random points in the unit disk, we can construct a graph from it in multiple ways. One is through Delaunay triangulation, another is to do graph sprinkling with $p=0.1,\,p=0.25,\,p=1$ for example.}
    \label{fig:TriangDist}
\end{figure}

\subsection{Constructing Hard Annulus Random Geometric Graphs}
In this paper, we improve upon earlier constructions of Random Geometric Graphs \cite{penrose_random_2003}, by including an additional parameter in order to lower metric-distortion. This is somewhat equivalent to a deterministic version of the 'soft annulus' connection function described in \cite{dettmann_random_2016}. If we have some region we want to approximate with a graph, we can generate random points (uniform with respect to the volume measure) which will be our vertices. If two vertices' corresponding points are distance (inside the region) $d$ apart, then we connect them if $|d-l|\leq l p$ where $l$ is our chosen \textbf{connection length} and $0<p<1$ is the \textbf{tolerance}. Theoretically $l$ and $p$ can be tuned to values that minimize the average distortion, allowing a very low distortion for a given amount of vertices. This would be very computationally expensive, but it turns out it already gives great results if a $p$ value of around $0.25$ is chosen, and the approximate minimum length $l$ that gives a connected graph (assuming the sprinkling region is connected) is found by binary search. We can see in Figure \ref{fig:TriangDist} that $p=0.25$ and $p=1$(corresponding to classic Random Geometric Graphs) have very similar distortions, but one finds that the amount of edges at $p=1$ is often very large in comparison with smaller $p$ values. Since the amount of edges influences memory and time requirements for working with the graphs, it is beneficial to have a smaller number of edges. A key feature of this method of constructing graphs is that the metric distortion can be lowered with increasing the number of vertices as can be seen in Figure \ref{fig:TriangDist} (and the choices of $l$ and $p$ values can be further tuned to get even closer to the theoretical minimum, but that is not necessary in this case). Note that the effective edge length will in general be related to, but not equal to $l$. In principle, our method above can lead to  graphs with even lower metric-distortion that those demonstrated here (via  optimization of $l$ and $p$); though, for the current purpose of validating our sectional curvature measure, the distortions we have achieved here are low enough to draw satisfactory conclusions. A different connection function might yield better results if the goal is ultra low distortion graphs.

\begin{figure}[H]
    \centering
    \includegraphics[scale=0.4]{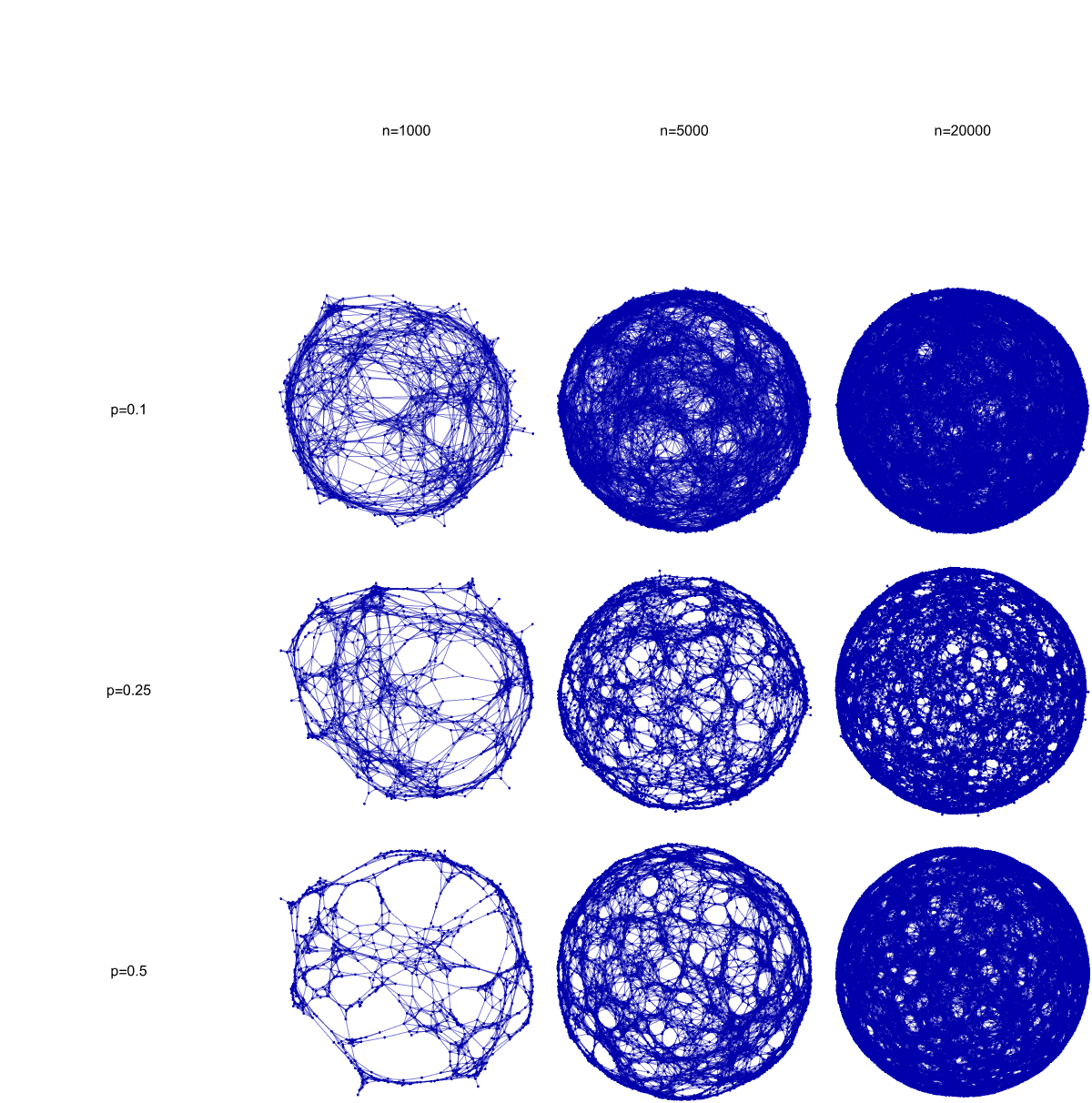}
    \caption{9 graphs sprinkled unto $S^2$, where rows correspond to $p=0.1,p=0.25,p=0.5$ from top to bottom, and columns corresponds to vertex counts of $n=1000,n=5000,n=20000$ from left to right. Note that unless specifically stated, all graph visualisations in this paper are done without vertex coordinate information from their sprinkling. We can see that for large $n$ we conveniently (in the sense of easing computations) have notably less edges around $p=0.25$.}
    \label{fig:sphereGrid}
\end{figure}

\section{Estimating Sectional Curvature of Random Geometric Graphs}\label{sec:results}
Since the discrete sectional curvature detailed above is based on the assumption of approximately constant curvature of the geodesic plane in which a given triangle lies we will first test our curvature estimator on graphs constructed using manifolds of constant curvature.

\subsection{Spheres: $S^2$ and $S^3$}

\begin{figure}[H]
    \centering
    \includegraphics[scale=0.15]{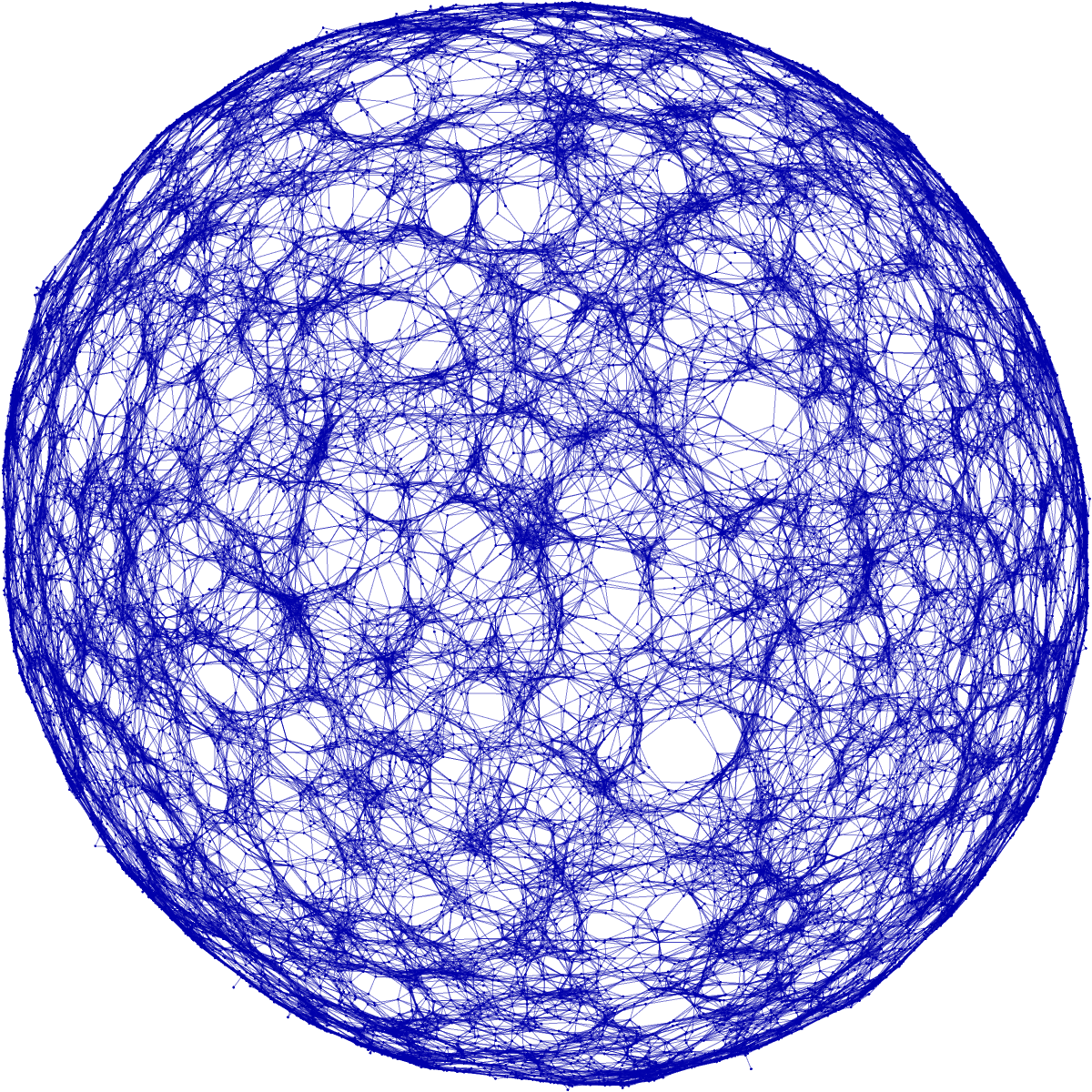}
    \includegraphics[scale=0.3]{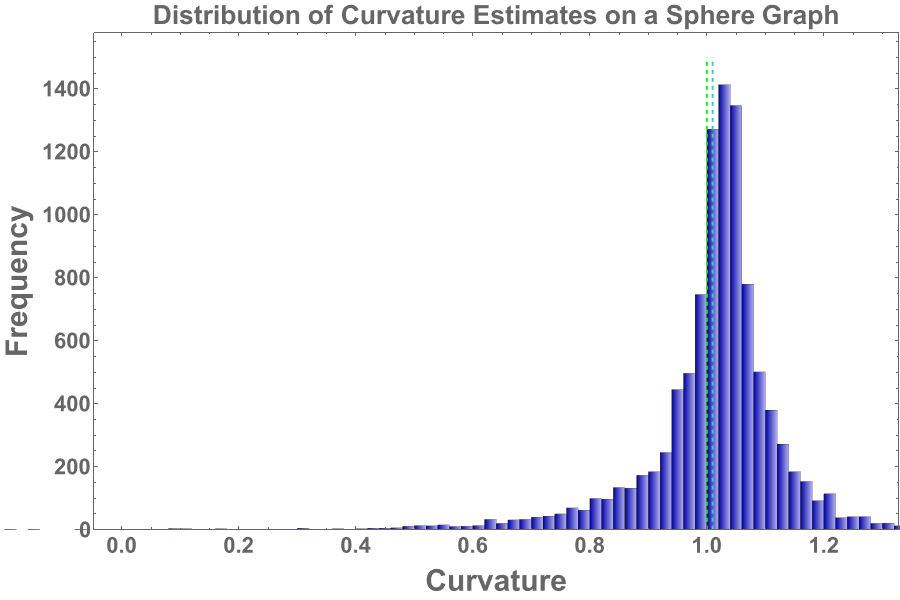}
    \caption{On the left is a random geometric graph associated with the two-dimensional unit sphere $S^2$ with $2\times 10^4$ vertices, $p = 0.25$, effective edge length $0.0436$ and an approximate metric distortion of $0.0306$. The histogram on the right shows a distribution of curvature estimates taken with $10^4$ samples. Using this, the mean curvature of the graph is estimated to be $1.0088\pm0.0013 $ (indicated by the light blue line in the distribution above), compared to the expected continuum curvature of $1$, indicated by the green dashed line. Furthermore (not shown in this figure), with much fewer samples ($10^2$), the estimated mean curvature still stands at $1.009\pm0.012$, indicating the robustness of the estimate. }
    \label{fig:bigSphere}
\end{figure}

\begin{figure}[H]
    \centering
    \includegraphics[scale=0.15]{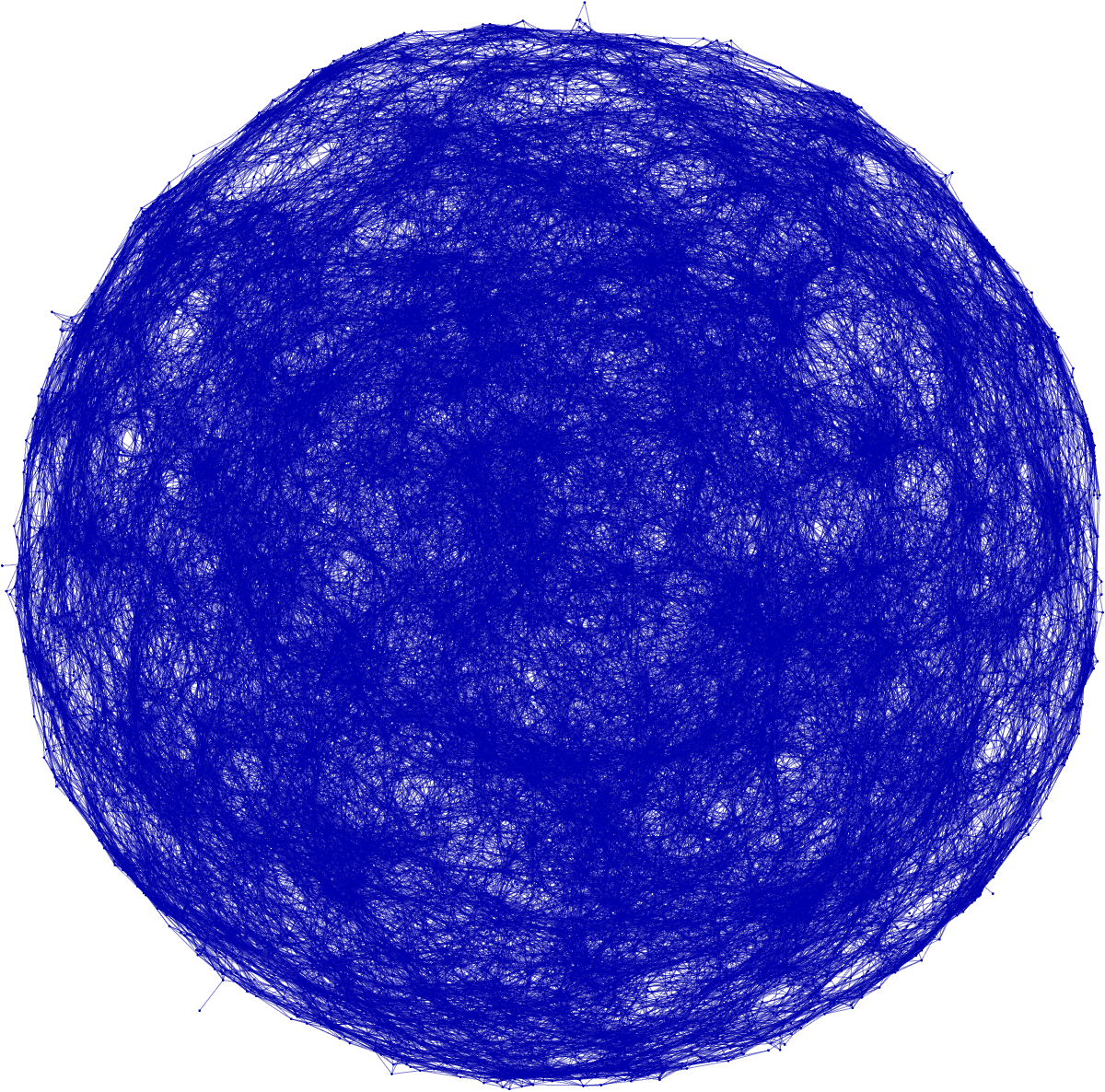}
    \includegraphics[scale=0.3]{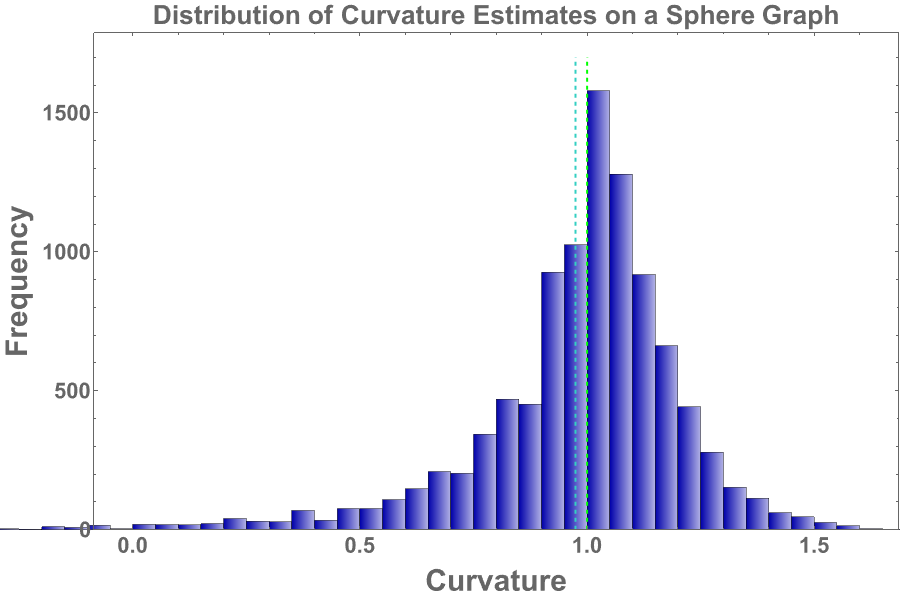}
    \caption{On the left is a graph associated with the three dimensional unit sphere $S^3$ with $3.2\times 10^4$ vertices, $p=0.25$, effective edge length $0.0858$ and an approximate metric distortion of $0.0405$. With $10^4$ samples the mean curvature is estimated at $0.974\pm0.0027 $.}
    \label{fig:big3DSphere}
\end{figure}

Here, we implement graph sprinkling onto $S^2$ and $S^3$. We find our curvature estimates approach the continuum value of $1$, and we  see that on the sphere graph in Figure \ref{fig:bigSphere} we get a very close estimate when using $2\times 10^4$ vertices. For the $S^2$ case we find that although the distribution of estimates is slightly asymmetric, it is quite narrow, and the mean corresponds well to the continuum value. For the $S^3$ case we see the distribution is slightly wider, likely due to the larger metric distortion, which comes from the need to have more vertices in $S^3$ than $S^2$ to have comparable densities of points. The mean however is still very close to the expected  continuum value, in spite of these challenges.

\subsection{Hyperbolic Plane}
\begin{figure}[H]
    \centering
    \includegraphics[scale=0.15]{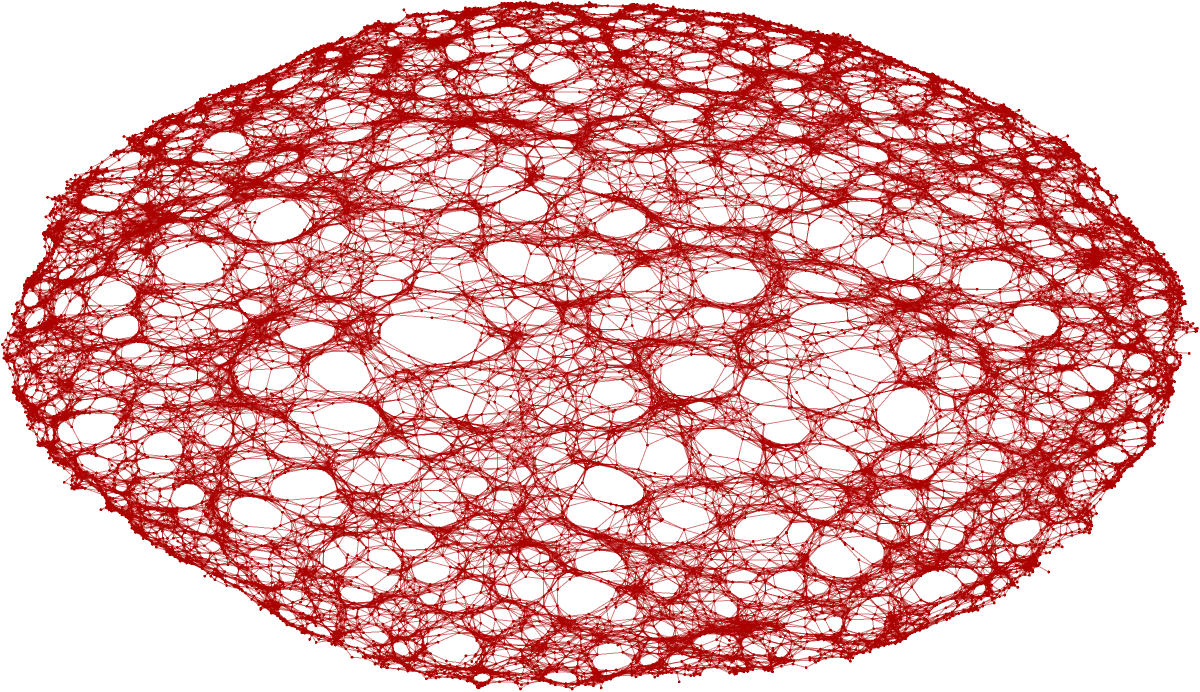}
    \includegraphics[scale=0.225]{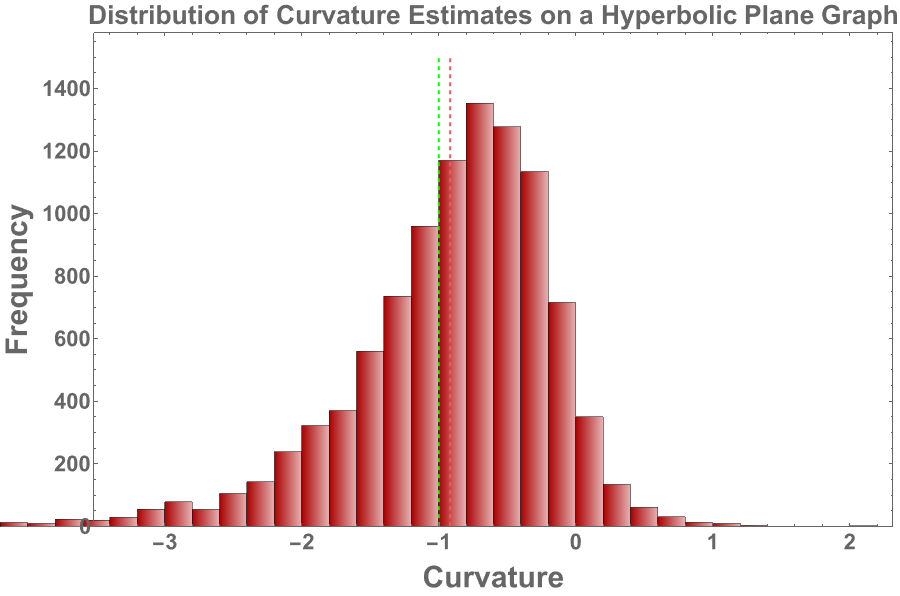}
    \caption{On the left is a graph of a circle of hyperbolic radius $\frac{\pi}{2}$ in the hyperbolic plane $H^2$ with $2\times 10^4$ vertices, $p=0.25$, effective edge length $0.0532$ and an average metric distortion of $0.0268$. The histogram on the right shows a distribution of curvature estimates taken with $10^3$ samples. This yields a mean curvature estimate of $-0.997\pm0.028 $ (compared to the expected value of $-1$), with the red dashed line indicating the $5\%$ trimmed mean of $-0.9341$. Furthermore (not shown in the figure), with much fewer samples ($10^2$),  the estimated mean curvature stands at $-0.82\pm0.08$, indicating the robustness of the estimate. }
    \label{fig:bigHP}
\end{figure}
\begin{figure}[H]
    \centering
    \includegraphics[scale=0.25]{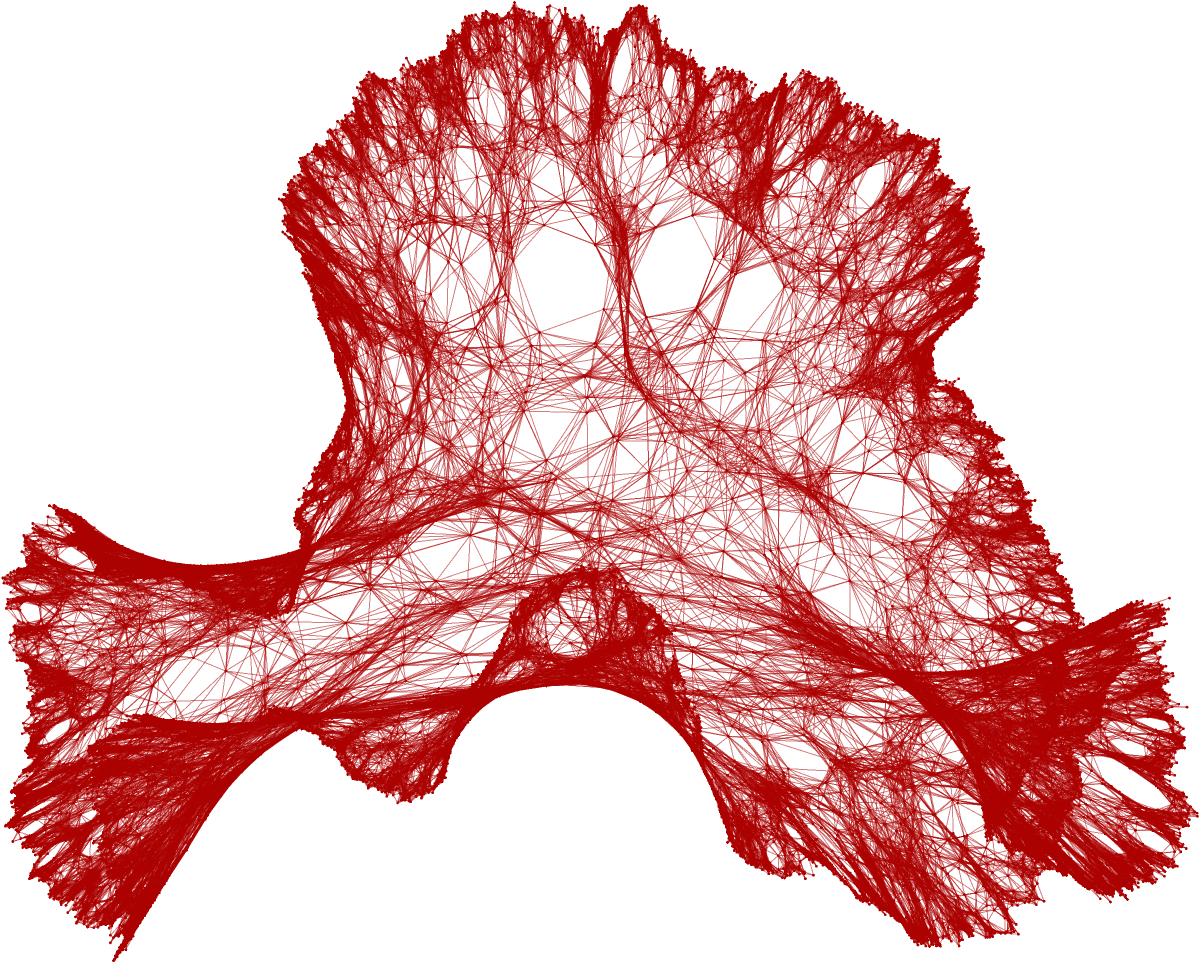}
    \caption{A graph of a disk of radius 5 in the unit hyperbolic plane.}
    \label{fig:curlHP}
\end{figure}

For the hyperbolic plane we need to choose which subset of the hyperbolic plane we sprinkle upon. We can choose the disk with hyperbolic area $4\pi$ (in analogy with the unit sphere) in $H^2$ where we would expect curvature estimates of approximately $-1$. As expected, the curvature estimates for negative curvatures are less well-behaved (in the sense of having a wider distribution) than for positive curvatures, as can be seen in the distribution in Figure \ref{fig:bigHP}. Looking at Figure \ref{fig:bigHP} one might be concerned with how flat the graph looks. This is due to the graph embedding algorithm, which shows the curvature for larger radii disks, as can be seen in Figure \ref{fig:curlHP}. One can interpret this as the curvature estimator being able to detect smaller deviations from euclidean geometry than the embedding algorithm, as a consequence of differential manifolds being locally euclidean.

\subsection{Euclidean Plane}
\begin{figure}[H]
    \centering
    \includegraphics[scale=0.15]{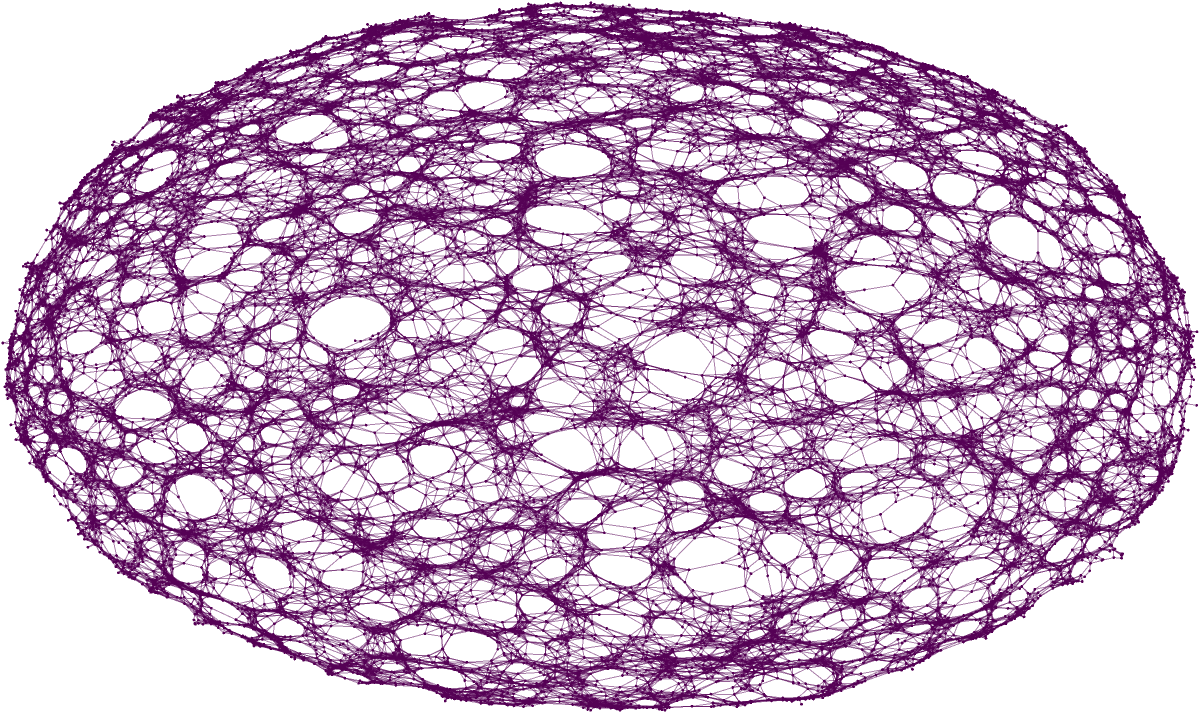}
    \includegraphics[scale=0.3]{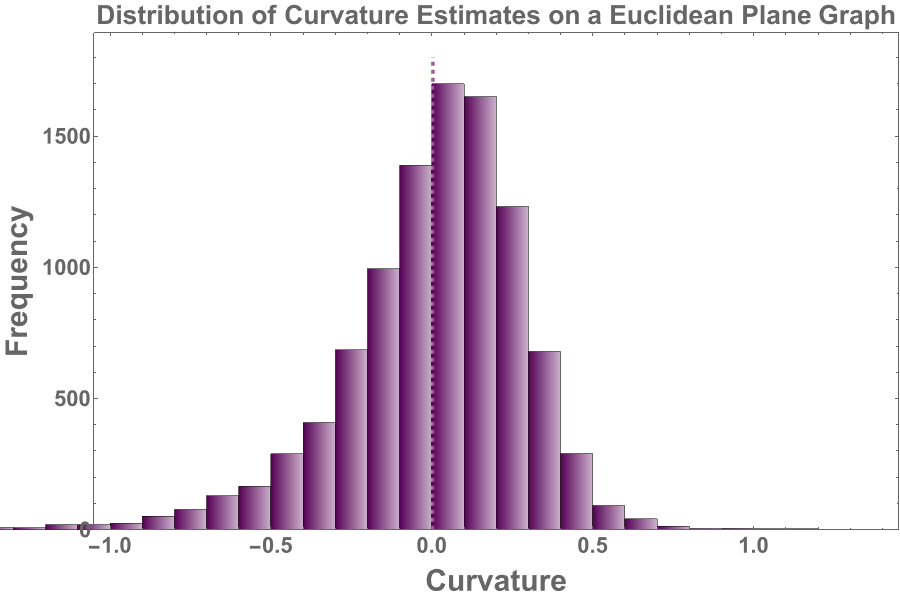}
    \caption{On the left is a graph of a radius 2 disk in $\mathbb{R}^2$ with $2\times 10^4$ vertices, $p=0.25$, effective edge length $0.0268$ and an average metric distortion of $0.0240$. With 100 samples the estimated curvature is $-0.04\pm0.16$ against an expected $0$. With $10^3$ samples the estimate is $0.09\pm0.04 $, and the distribution can be seen on the right, with the pink dashed line indicating the median of $0.2275$.}
    \label{fig:bigP}
\end{figure}
\begin{figure}[H]
    \centering
    \includegraphics[scale=0.35]{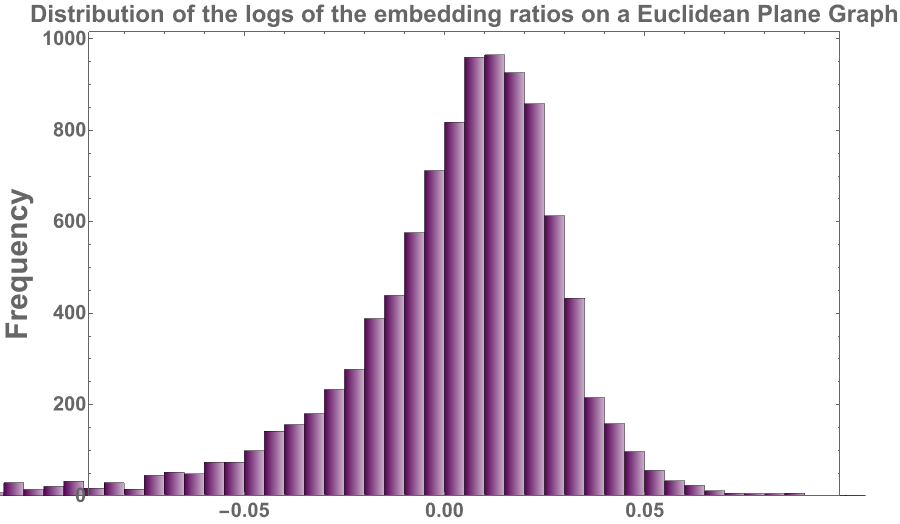}
    \caption{Distribution of the logs of embedding ratios of the graph shown in Figure \ref{fig:bigP}.}
    \label{fig:planDist}
\end{figure}

Interpreting the error in curvature of a euclidean graph is essentially impossible, since there is no intrinsic length scale present in the euclidean plane. Nevertheless, the distribution of the curvature estimates can be seen in Figure \ref{fig:bigP}, where the shape is the only relevant information, given that the x-axis can be arbitrarily re-scaled. We notice however that this distribution is qualitatively similar to the sphere and hyperbolic plane cases. Looking at the distribution formed by the logs of the normalised embedding ratios in Figure \ref{fig:planDist}, we  see that it also has a qualitatively similar distribution, demonstrating how metric distortion influences curvature estimates.

\section{Convergence Tests of Sectional Curvature Estimates}
\label{sec:convergence}

In order to test how our curvature estimates converge to the sectional curvature of the underlying manifold when the metric distortion tends to zero, there are two types of errors to consider. The most important is the absolute error of the mean curvature estimate, since the mean of a sample of curvature estimates serves as a good estimator. The mean of the absolute error of each individual sample however would also be of interest, since if it goes to zero quickly, that would then imply that less samples are needed for a robust estimate. There is no reason to expect there to be a linear relationship between the distortion and error, in fact there is good reason to think it is non-linear, especially at low distortions. We can however still fit a linear model, under the assumption that the relationship is roughly linear within the portion of the distortion we are probing. Then the $r^2$ value of this fit  gives us an idea of how close we are to a reasonable linear relationship between error and distortion in this region, which serves as a way to quantify the correlation between the metric distortion and the two types of errors.

\begin{figure}[H]
    \centering
    \includegraphics[scale=0.25]{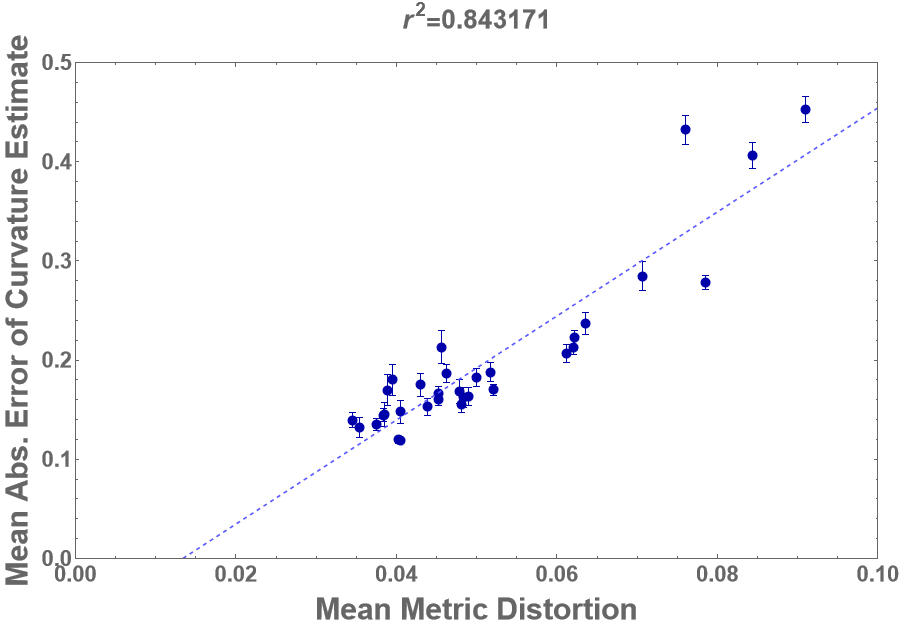}
    \includegraphics[scale=0.25]{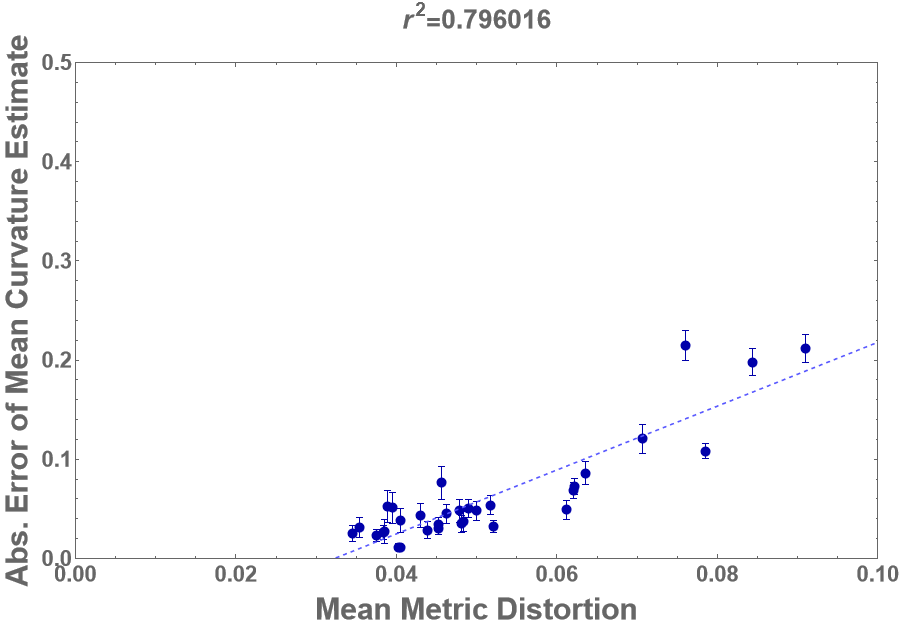}
    \caption{The Mean Absolute Error of Curvature Estimates, as well as the Absolute Error of Mean Curvature Estimates is plotted against the Mean Metric Distortion of graphs sprinkled into $S^2$. The fast convergence of the absolute error of the mean motivates using the mean estimate over a sample as a reasonable estimator, and the convergence of the mean absolute error shows that smaller sample sizes are warranted in positive curvature geometries.}
    \label{fig:sphrErDist}
\end{figure}

\begin{figure}[H]
    \centering
    \includegraphics[scale=0.25]{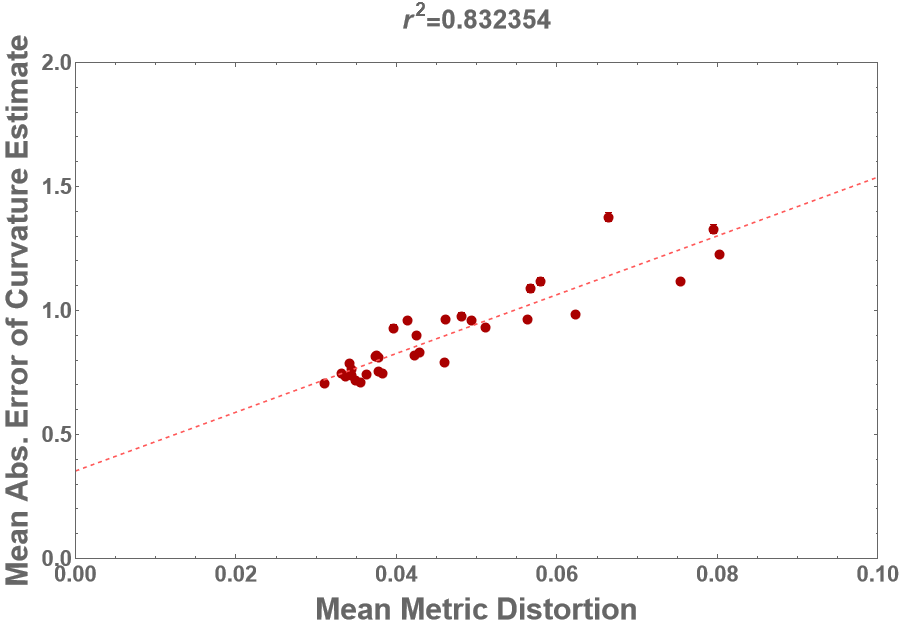}
    \includegraphics[scale=0.25]{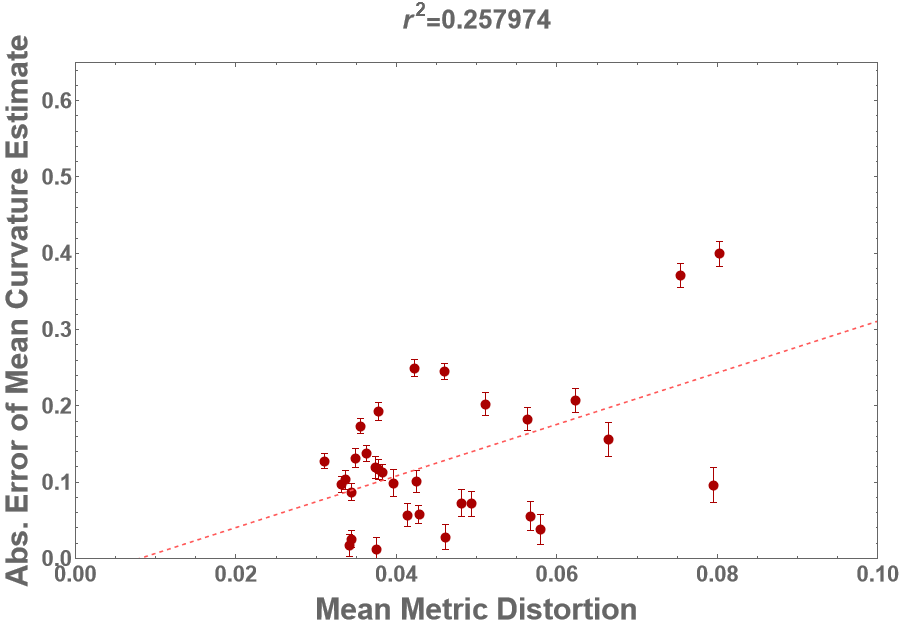}
    \caption{The Mean Absolute Error of Curvature Estimates, as well as the Absolute Error of Mean Curvature Estimates is plotted against the Mean Metric Distortion of the graphs sprinkled into $H^2$. The fast convergence of the absolute error of the mean motivates using the mean estimate over a sample as a reasonable estimator, and the slower convergence (with the question of if it truly converges to 0 as expected) of the mean absolute error shows that larger sample sizes are warranted in negative curvature geometries.}
    \label{fig:hypErDist}
\end{figure}

\begin{figure}[H]
    \centering
    \includegraphics[scale=0.25]{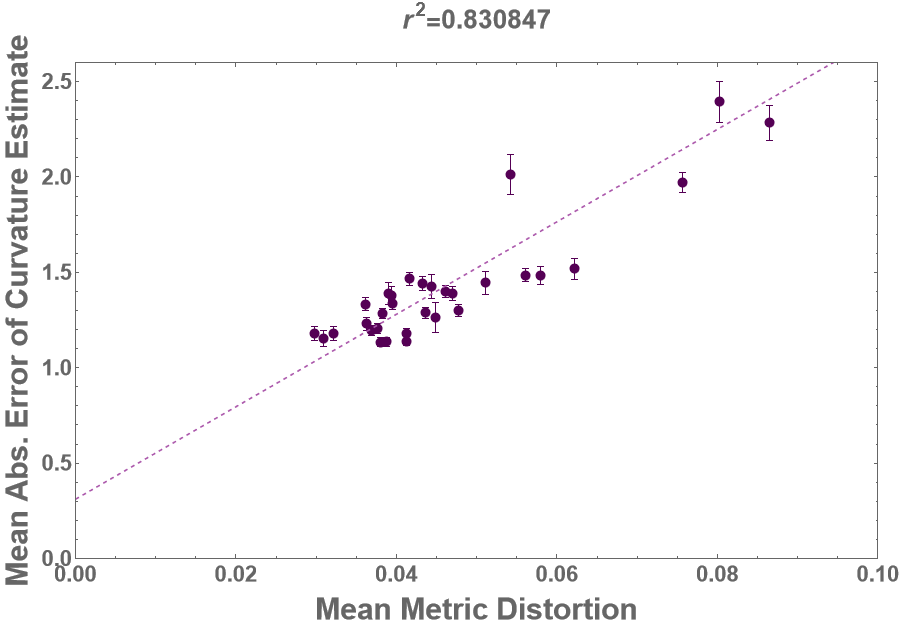}
    \includegraphics[scale=0.25]{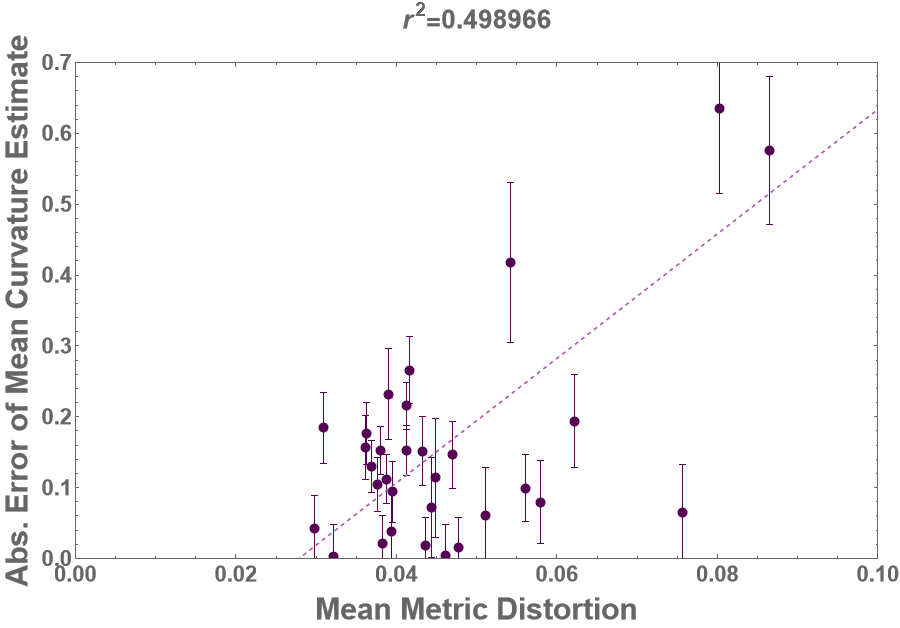}
    \caption{The Mean Absolute Error of Curvature Estimates, as well as the Absolute Error of Mean Curvature Estimates is plotted against the Mean Metric Distortion of the graphs sprinkled into $\mathbb{R}^2$. Both the mean absolute error and the absolute error of the mean shows convergence, but any statement about the speed or uniformity of the convergence is dependent on the scaling of the x and y axis, both of which are arbitrary in euclidean geometries where there are no inherent length scales.}
    \label{fig:eucErDist}
\end{figure}

\section{A Vertex-Based Curvature Estimation}  
\label{sec:vertex}

Although discrete sectional curvature is generally concerned with  the curvature over triangles, or in regions, we can still extract some notion of a vertex specific curvature by averaging out the curvatures of some triangles that contain the given vertex in the graph. We can therefore choose a sample of triangles, calculate their curvatures, and then from that we can extract some average curvature for each vertex of a graph. For some graph representing constant curvature we would want the deviation between adjacent vertices to be small relative to the deviation of the estimate overall.

\begin{figure}[H]
    \centering
    \includegraphics[scale=0.45]{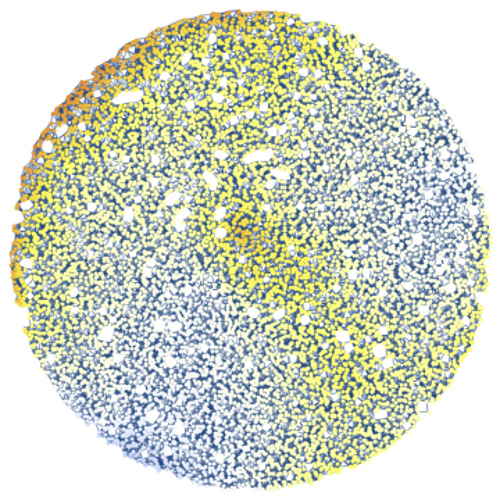}
    \includegraphics[scale=0.4]{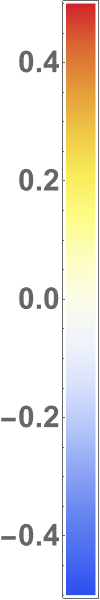}
    \caption{A vertex-based sectional curvature estimation for a graph sprinkled on a flat disc. The mean curvature estimate based on triangles associated to each vertex have been displayed (using a color scale to indicate the value of the estimate) in the graph above. Warmer colors indicate positive curvature, while the cooler ones, negative. }
    \label{fig:vertCurvs}
\end{figure}

In Figure~\ref{fig:vertCurvs} we  see that the estimated curvature per vertex is slowly varying over the graph considered here, and doesn't display the boundary effects seen in \cite{denne_convergence_2008}. We also see how the mesoscopic nature of our discrete sectional curvature is exhibited by the fact that there is little variation in curvature from vertex to vertex. There still are curvature variations over regions of the graph. These indicate that the graph has certain regions of more or less dense due to the graph sprinkling algorithm used in its construction.

\section{Comparisons to Other Discrete Curvature Definitions}
\label{sec:compare}
Many existing discrete curvature definitions such as  Ollivier-Ricci (original version), Forman-Ricci, etc, have been applied to a variety of complex networks, including random networks, power-law networks, etc. However, generic complex networks may or may not be endowed with sufficient geometric structure  (or at least are not constructed using geometry). On the other hand, random geometric graphs (and those related to those) are constructed using geometric data and hence serve as ideal testing grounds for discrete curvature definitions. This point is made clear from examples in which the Forman-Ricci curvature is computed on combinatorial graphs. There one can show that the average curvature will be negative in any combinatorial graph with average vertex degree of more than 2. 

\subsection{Wolfram-Ricci Curvature}
In the following, let us compare our sectional curvature estimate to the definition of discrete curvature used in \cite{Wolfram2020}, where the Taylor expansion for the volume of geodesic balls in constant curvature 
\begin{equation}
V_r = a r^d \left( 1 - \frac{R}{6(d+2)} r^2 \right)
\end{equation}
is used. Here $a$ is some normalization constant, related to the effective volume of each vertex, $R$ is the Ricci scalar and $d$ is the dimension. We will only be testing in graphs representing 2 dimensional surfaces here, so by noting that the Ricci scalar is then related to sectional curvature as $K=\frac{R}{2}$, so we can simplify the expression to 
\begin{equation}\label{eq:WR}
V_r = a r^2 \left( 1 - \frac{K}{48} r^2 \right)    
\end{equation}
Naive implementations of this measure turn out to be unstable for the random geometric graphs we are considering here, so instead we implement a method where vertex-count of successive geodesic balls around some vertex is taken up to a radius where $K r^2\approx 1$, where we expect the expansion to break down. Simple curve fitting is then employed, giving stable results, but with large error bars, presumably due to the equation being fitted being relatively insensitive to the curvature parameter for the small radii we are forced to consider. We can compare the errors versus metric distortions to contrast it with discrete sectional curvature, in order to ascertain the different regimes of effectiveness. The relevant results from the discrete sectional curvature will be shown in lighter colours in the background to aid the comparison.

\begin{figure}[H]
    \centering
    \includegraphics[scale=0.185]{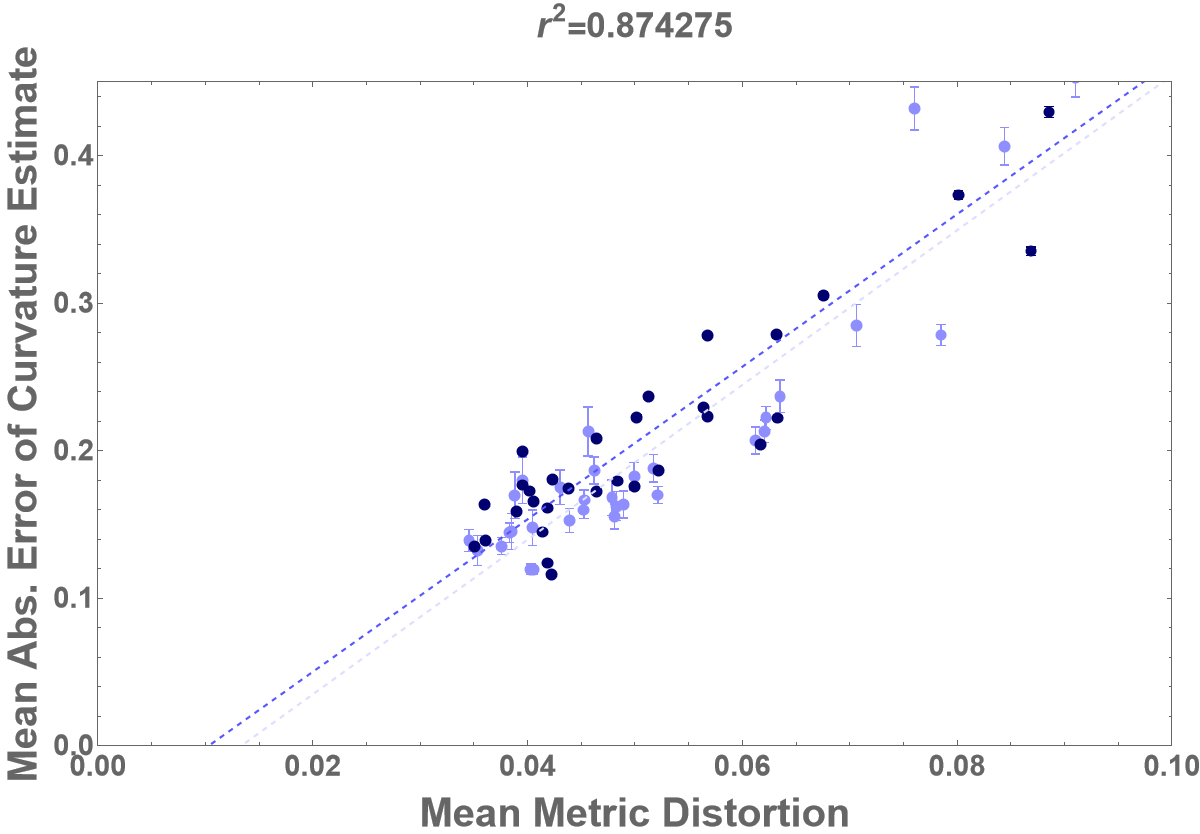}
    \includegraphics[scale=0.185]{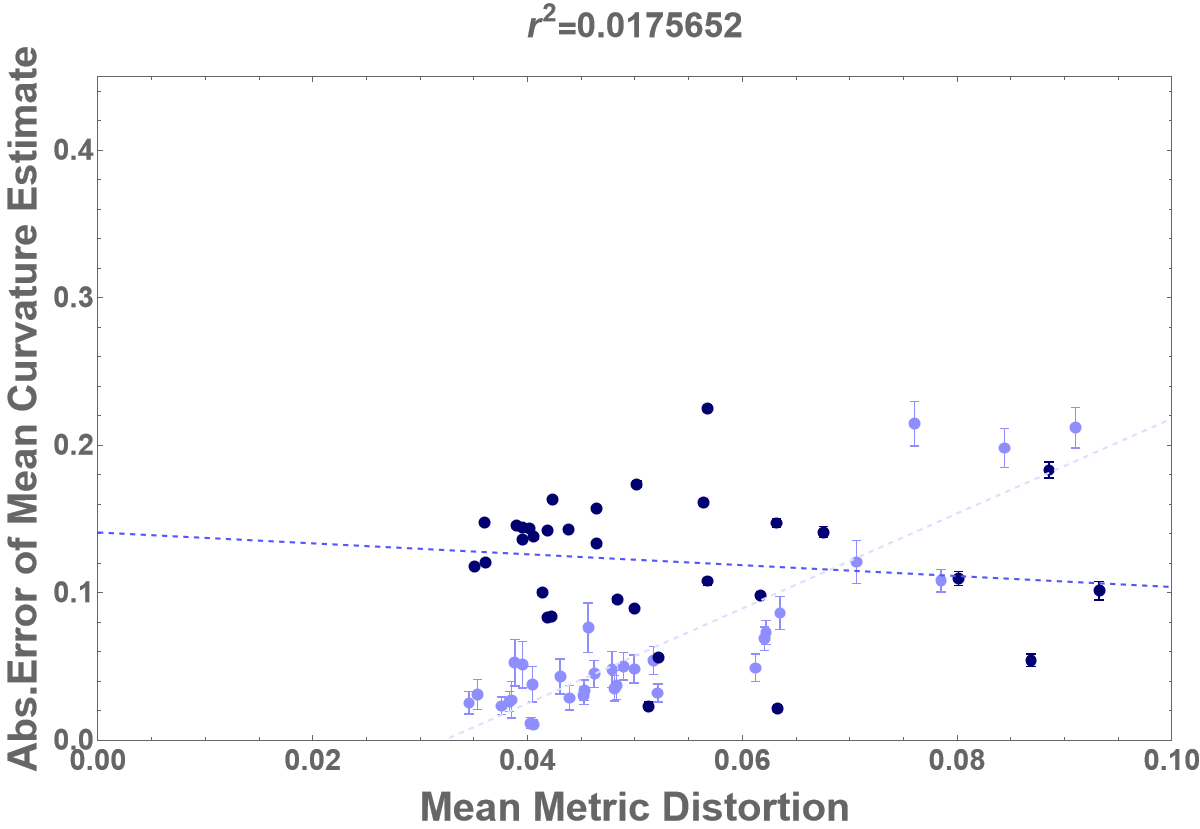}
    \caption{A comparison between the errors of Wolfram-Ricci curvature (dark blue), as well as discrete sectional curvature (very light blue), vs the mean metric distortion in positive curvature}
    \label{fig:wolfSphere}
\end{figure}

\begin{figure}[H]
    \centering
    \includegraphics[scale=0.185]{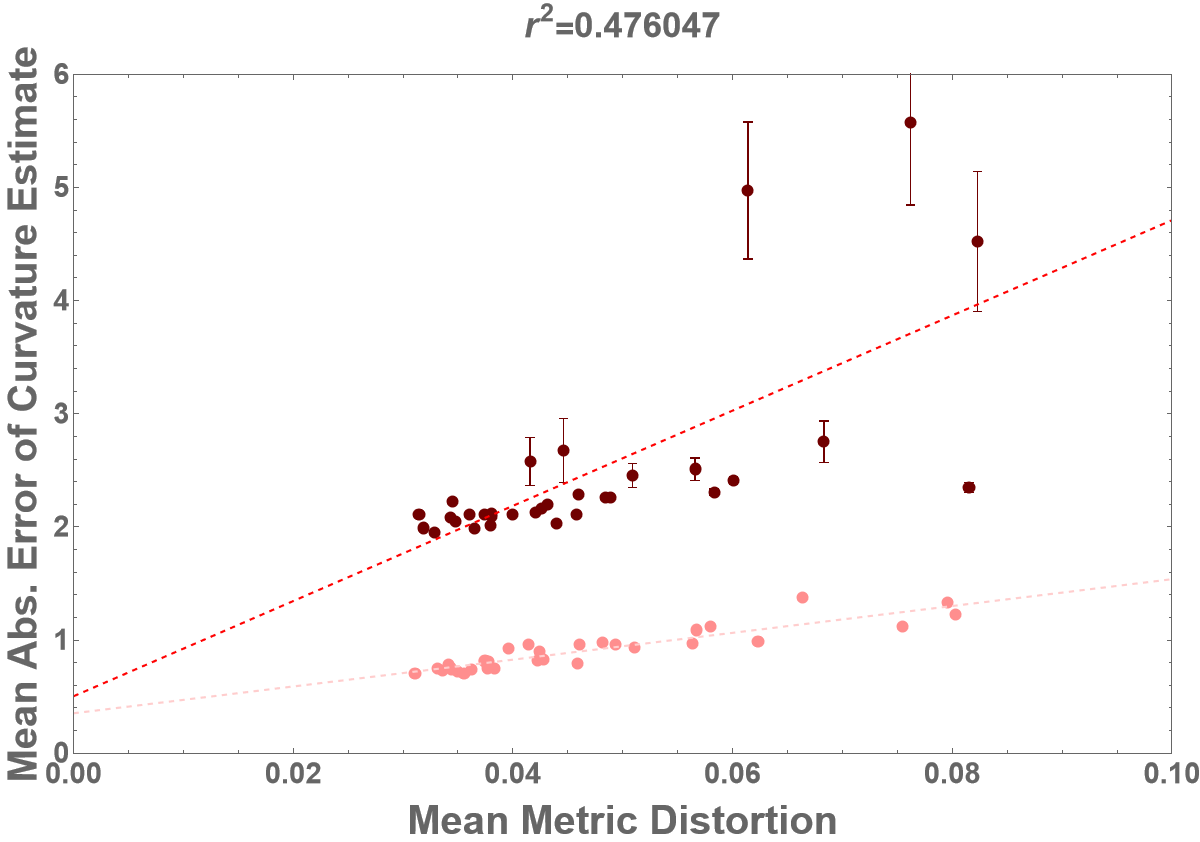}
    \includegraphics[scale=0.185]{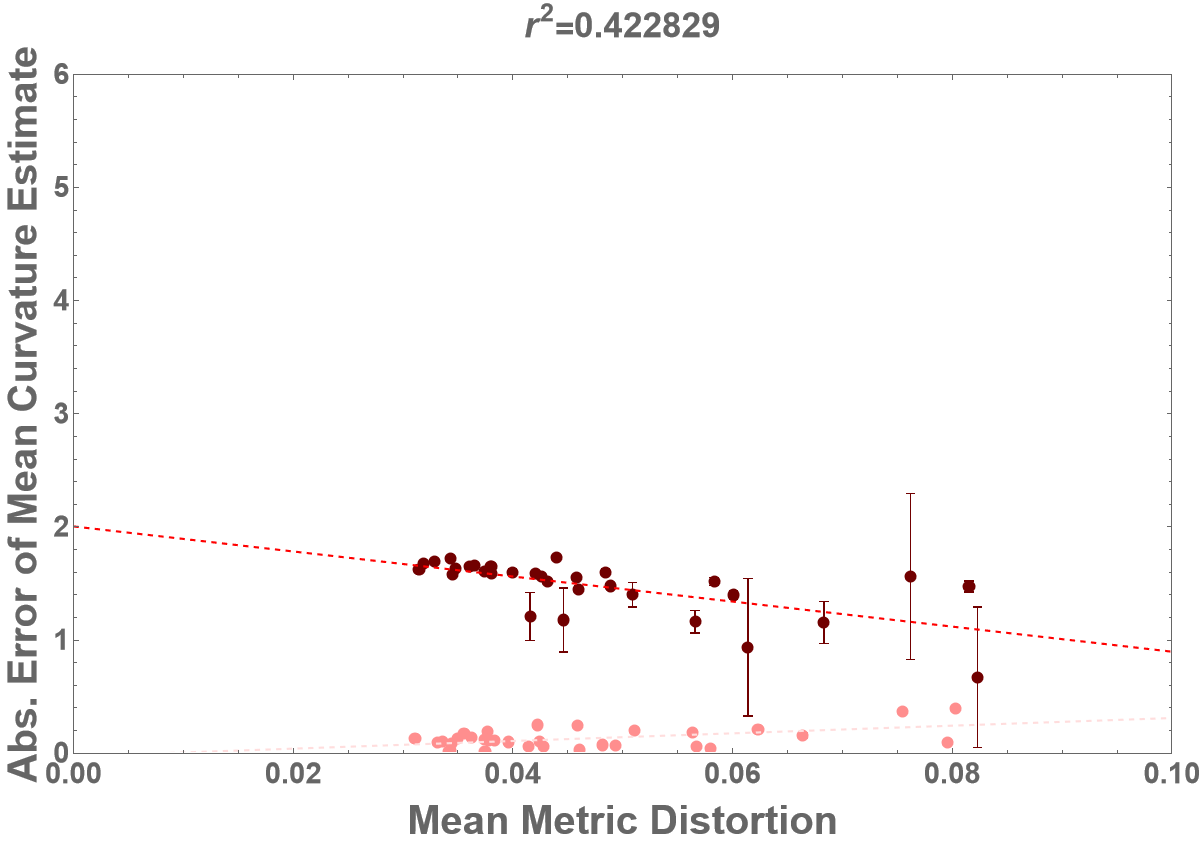}
    \caption{A comparison between the errors of Wolfram-Ricci curvature (dark red), as well as discrete sectional curvature (very light red), vs the mean metric distortion in negative curvature}
    \label{fig:wolfHyperbolicPlane}
\end{figure}

\begin{figure}[H]
    \centering
    \includegraphics[scale=0.185]{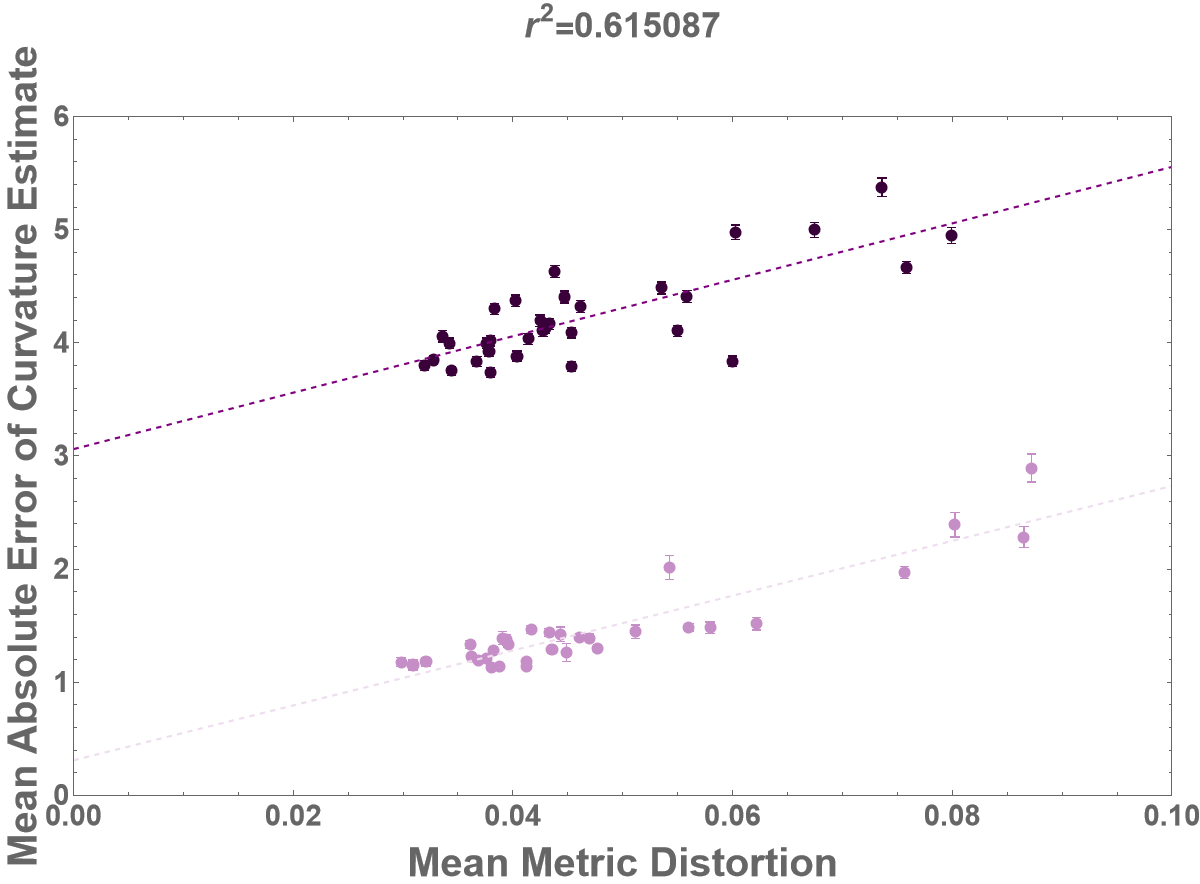}
    \includegraphics[scale=0.185]{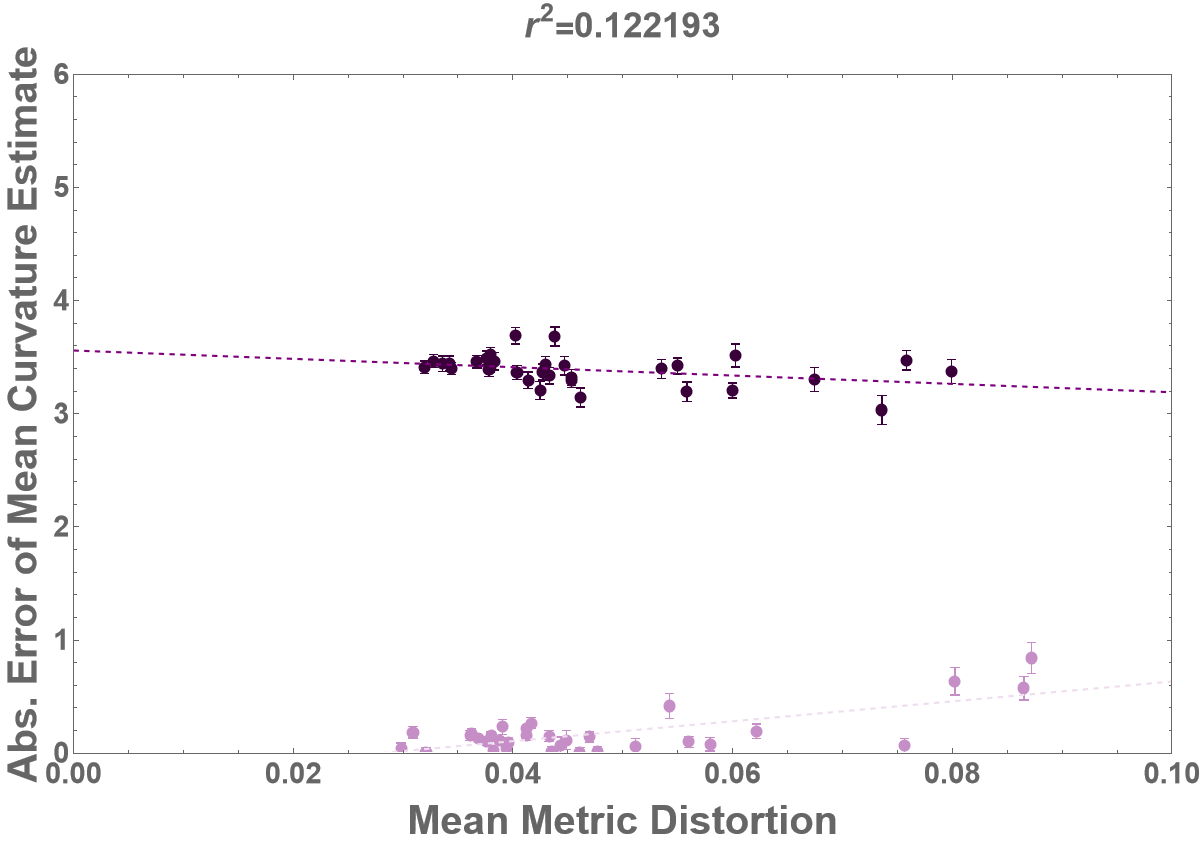}
    \caption{A comparison between the errors of Wolfram-Ricci curvature (dark purple), as well as discrete sectional curvature (very light purple), vs the mean metric distortion with $0$ curvature}
    \label{fig:wolfEuclideanPlane}
\end{figure}

As we can see the two methods are comparable in efficacy per sample for the positive curvature case, however the mean shows essentially no correlation with the metric distortion. For the negative and zero curvature cases discrete sectional curvature is notably better when considering per sample error, and we see the same lack of convergence of the mean estimate of the Wolfram-Ricci curvature.

\subsection{Mesoscopic Ollivier-Ricci Curvature}
\begin{figure}[H]
    \centering
    \includegraphics[scale=0.25]{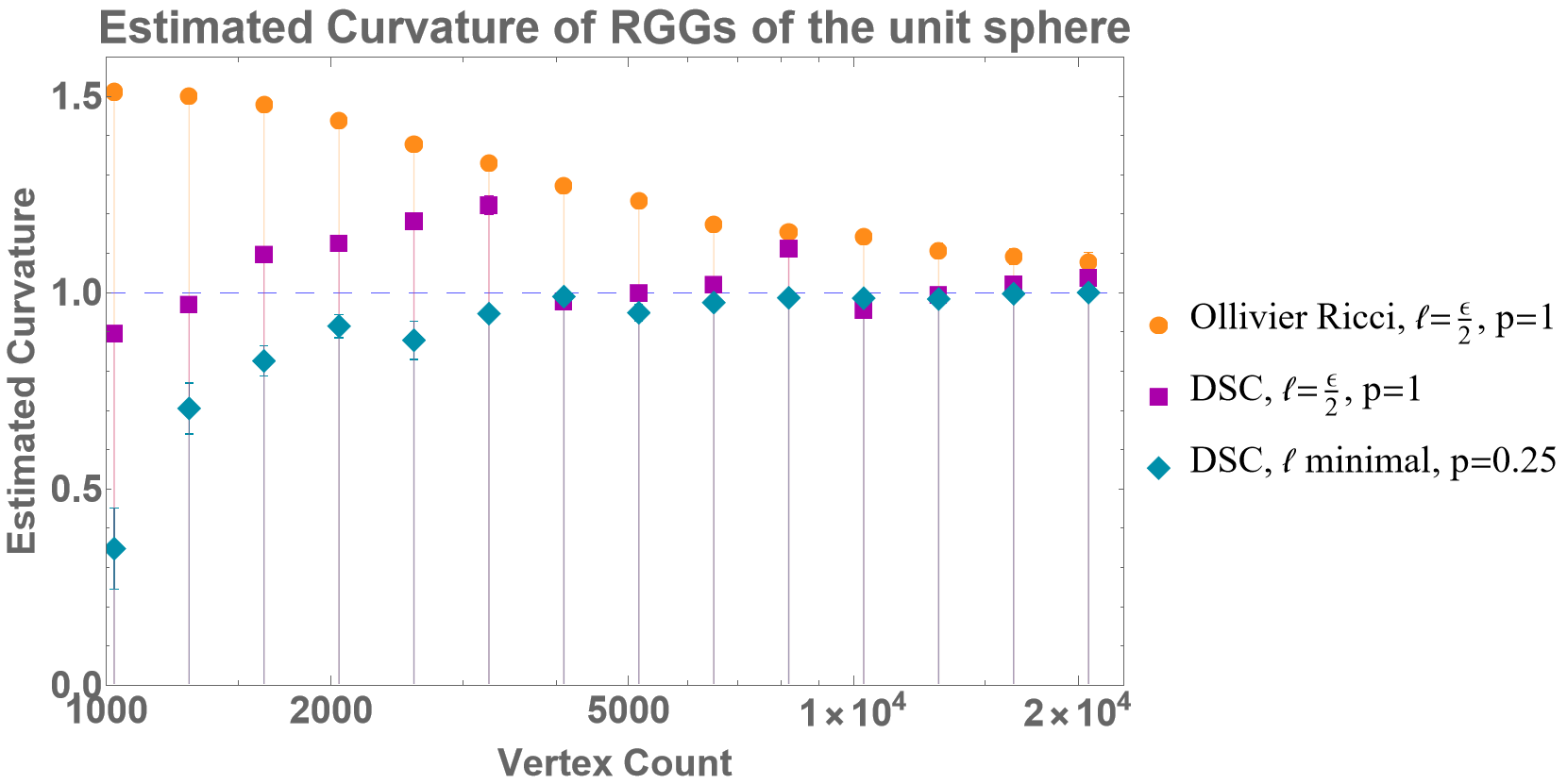}
    \caption{Convergence profiles as a function of vertex count. The discontinuities in the Discrete Sectional Curvature (DSC) with $\ell=\frac\epsilon2$ is a consequence of how $\ell$ is determined.}
    \label{fig:olComp}
\end{figure}

The difficulties of convergence of Ollivier-Ricci curvature on geometric graphs is remedied in \cite{PhysRevResearch.3.013211} by generalizing the definition to work with mesoscopic graph neighbourhoods. There the types of graphs under consideration are RGG's with connection radius $\epsilon$, which corresponds to our "hard annulus RGG's" with connection length $l=\frac{\epsilon}{2}$ and tolerance $p=1$. We note that for the mesoscopic definition used in \cite{PhysRevResearch.3.013211}, the geometric information is encoded both in graph neighbourhoods and in the metric, and therefore the graphs that were reasonably considered in the analysis were such that the average vertex degree, as well as the graph radius, diverges in the continuum limit. This is in slight contrast with the methods developed in this paper, where only the metric encodes the information, and the graphs considered generally have bounded, or in comparison small, average vertex degree, while having much larger and diverging graph radius. Nevertheless, we can look at the convergence profile showed in Figure 2(b) in \cite{PhysRevResearch.3.013211} and attempt to create a comparison, to be plotted against the data provided by the authors. We can firstly test the discrete sectional curvature by using both their vertex counts and connection radii (with $p=1$), and then we can also test convergence where only the vertex counts are used, and the connection length is taken to be the minimal length that gives a connected graph (again with $p=1$). We can see the relative profiles in Figure \ref{fig:olComp}.

We can see that the estimate of the discrete sectional curvature with $l=\frac{\epsilon}2$ is quite discontinuous as a function of vertex size, which is due to the low graph diameter, so each increase in minimum length scale is quite dramatic. It is however for larger graphs seemingly more convergent that the mesoscopic Ollivier-Ricci curvature. We can see the discontinuities disappear (or rather become small relative to other errors) when using the minimal connection length $l$. With this minum $l$ we can see very rapid convergence to the continuum value. This begs the question of how the mesoscopic Ollivier-Ricci curvature would act on the type of random geometric graphs considered in this paper, which would be a natural extension of this work.

\section{Other Applications of Discrete Curvature}
\label{sec:apps}

\subsection{Estimating the Radius of the Earth}
\begin{figure}[H]
    \centering
    \includegraphics[scale=0.23]{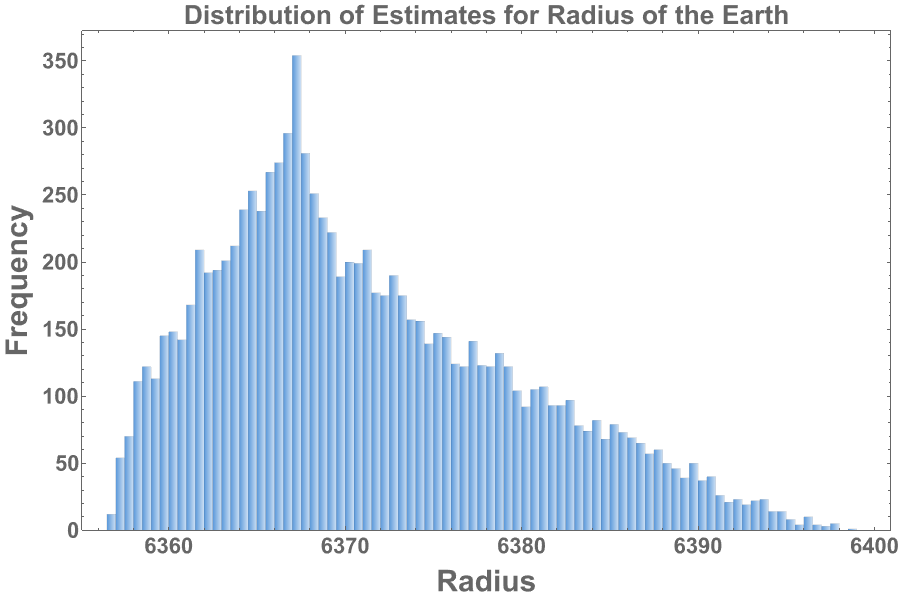}
    \includegraphics[scale=0.23]{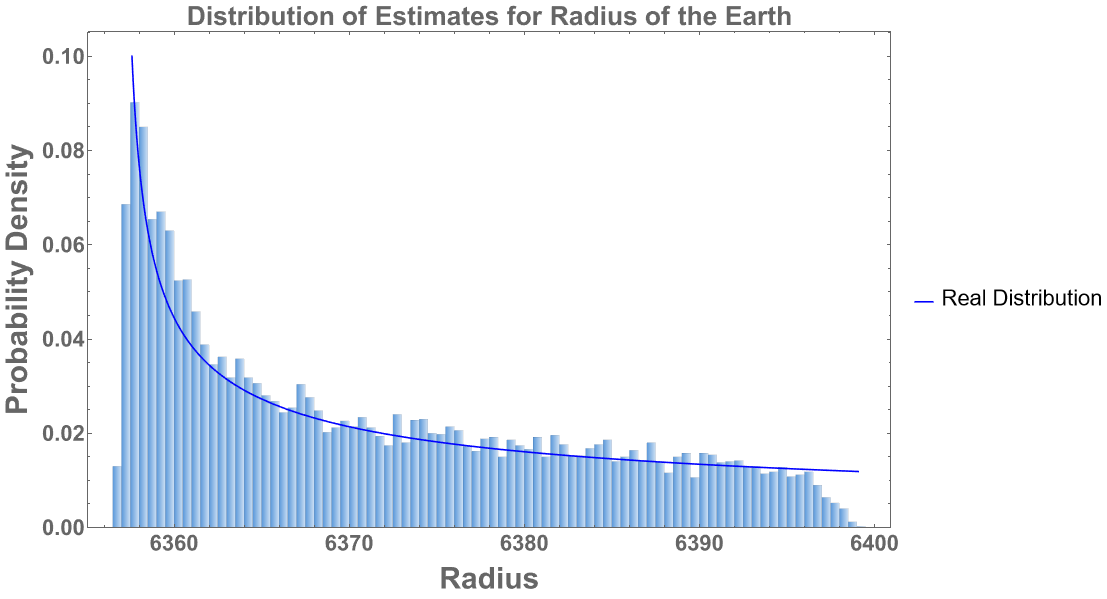}
    \caption{The left plot shows the frequency of different estimated radii over $10^4$ samples. The right plot shows the probability density of different radius estimates over $10^4$ with a maximum length scale of $\SI{6400}{\kilo\meter}$ as compared to the expected probability density for a oblate spheroid with parameters comparable to that of earth. Introducing a maximum length scale allows us to recover not only the mean curvature, but a good approximation of the distribution of sectional curvature itself.}
    \label{fig:earthRadiusDist10000}
\end{figure}

Here, we use the geographic data in-built in $Mathematica$. We model the earth as an oblate spheroid, and use geographical distances with which we can calculate curvature. If we take $R=\frac{1}{\sqrt{K}}$ as an estimate of the earths radius (which differs depending on latitude), we can investigate if this gives a good estimate. This can of course in principle be done by measuring distances in real life, provided the length scales involved are large relative to topographic fluctuations, such as mountains or hills. If we take 100 random triangles over the surface of the earth, we get an estimate for the average radius of the earth of  $\SI[separate-uncertainty=true]{6370.5\pm0.7}{\kilo\meter}$ compared to the commonly reported average of $\SI{6371}{\kilo\meter}$. If $10^4$ samples is taken (with distribution that can be seen on the left of Figure \ref{fig:earthRadiusDist10000}) the estimated average is $\SI[separate-uncertainty=true]{6371.3\pm0.4}{\kilo\meter}$. We can then introduce a maximum length scale on the order of the radius (which is small relative to the length scale introduced by the small curvature gradient of the earth), which we can see allows us to approximately reproduce the radius distribution we would expect to find on an oblate spheroid. We get this distribution by a transformation of variables of the distribution of the area of the earth as a function of parametric latitude, and transforming it with the inverse square root of the formula for the Gaussian curvature of a spheroid as a function of parametric latitude. The polar radius of the earth is taken as being approximately $\SI{6357}{\kilo\meter}$, and the equatorial radius is taken as being approximately $\SI{6378}{\kilo\meter}$. This gives us an approximate distribution of 
\begin{equation}
    p(R)=\frac{0.077088}{\sqrt{R-6357}}\qquad6357\leq R\leq 6399.07
\end{equation}
which can be seen on the right of Figure \ref{fig:earthRadiusDist10000}. This distribution has a mean of $\SI{6371.07}{\kilo\meter}$ as expected.

\begin{figure}[H]
    \centering
    \includegraphics[scale=0.4]{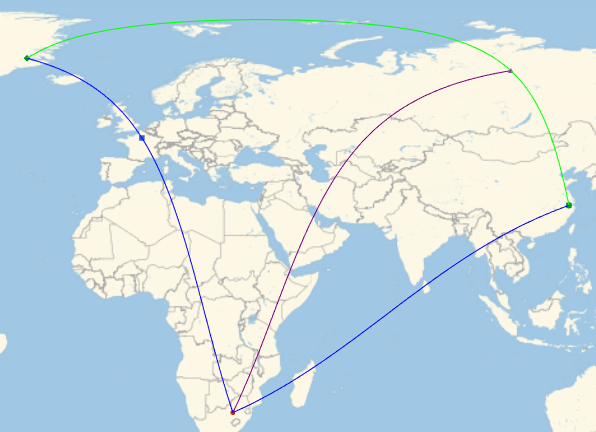}
    \caption{Constructing a large right-angled triangle on the surface of the earth. Note the distortion of angles and distances from the non-conformal projection, since Wolfram Mathematica employs the Equirectangular projection for large scale visualizations.}
    \label{fig:globalTriang}
\end{figure}

\subsection{Investigating Sectional Curvature Distribution of Fractals}
Defining and estimating geometric notions such as curvature and dimensions of fractals are challenging for many reasons \cite{winter2008curvature}. Previous work in this direction can be found in the following papers  investigating curvature measures, also known as Melnikov curvature \cite{mel1995analytic}  of fractal sets  \cite{winter2008curvature,winter2013fractal}. However, the Melnikov curvature refers to a notion of curvature pertaining to measures and is distinct from the notion of sectional curvature we are investigating here. 

\begin{figure}[H]
    \centering
    \includegraphics[scale=0.3]{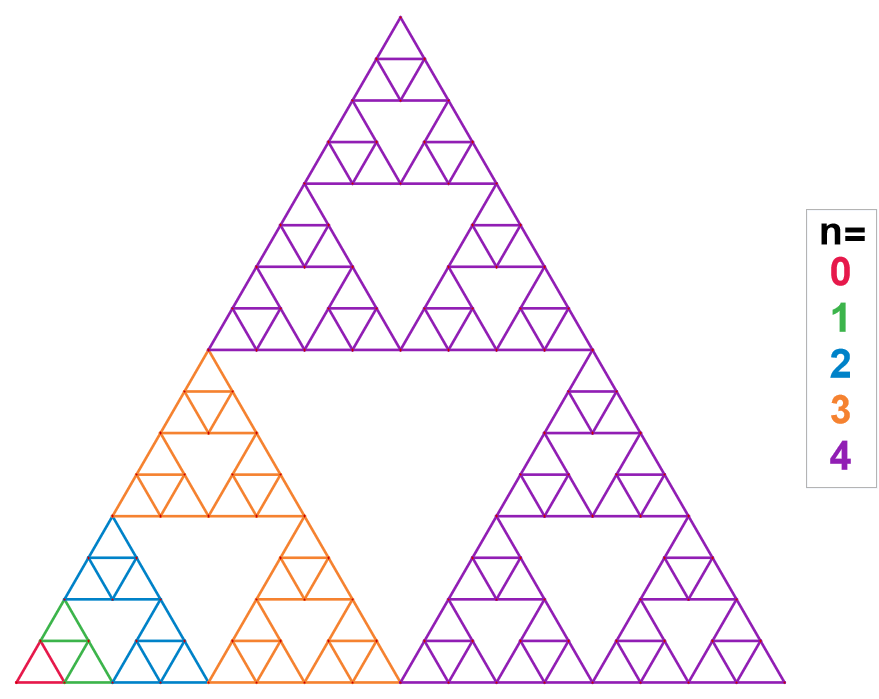}
    \caption{The Sierpinski triangle graphs $\hat{S}_3^0,\hat{S}_3^1,\hat{S}_3^2,\hat{S}_3^3,\hat{S}_3^4$, each embedded in the next iteration and color-coded according to the legend on the right.}
    \label{fig:sierTriangs}
\end{figure}

In what follows, we describe an exploratory application of our discrete sectional curvature applied to the Sierpinski triangle. First of all, notice that due to the self-similar nature of the Sierpinski triangle, barring the introduction of some length scale such as a maximum or minimum triangle size, there isn't a length scale in the system, meaning that any dimension-full quantity such as sectional curvature (with dimensions of length$^{-2}$) will formally be either vanishing or infinite. Here, we will consider the so-called Sierpinski triangle graphs $\hat{S}_3^n$ as classified in \cite{HINZ2017565}, as a sequence of better and better approximations of the Sierpinski triangle as can be seen in Figure~\ref{fig:sierTriangs} (where we can arbitrarily choose the edge lengths to be 1; we will return to this point later). We then investigate the sequence of sectional curvature distributions by considering the discrete sectional curvature of the isosceles triangles with the third side of even length (to avoid any possible ambiguities or arbitrariness) as shown in Figure~\ref{fig:distGrd} below.

\begin{figure}[H]
    \centering
    \includegraphics[scale=0.4]{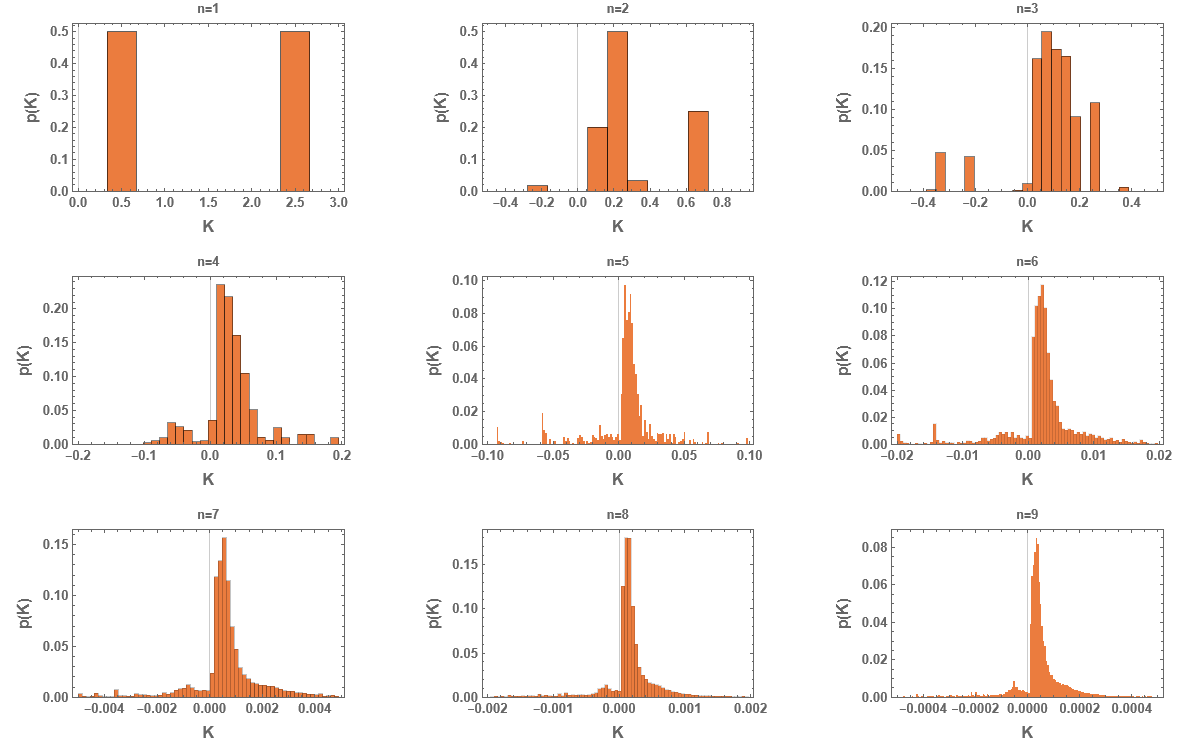}
    \caption{Sectional curvature distributions for the Sierpinski triangle graphs $\hat{S}_3^n$.}
    \label{fig:distGrd}
\end{figure}

\begin{figure}[h]
    \centering
    \includegraphics[scale=0.4]{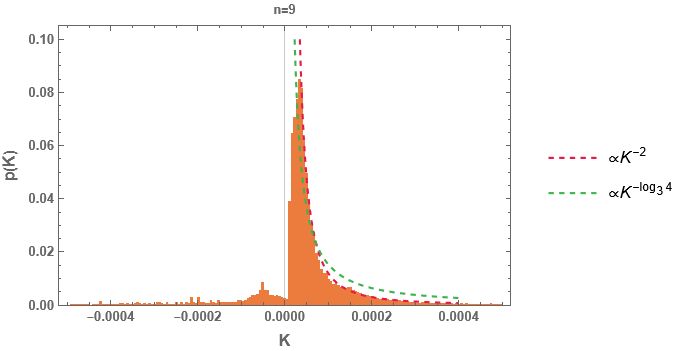}
    \caption{Naive power-law fits of the sectional curvature distribution of $\hat{S}_3^9$.}
    \label{fig:distFit}
\end{figure}

For $n=9$ in Figure~\ref{fig:distGrd} we find what appears to be a power-law behaviour of the distribution. Naively one might reason that since for each triangle with some curvature $K$ there are (in the full fractal) always $3$ triangles with half the lengths scale, and therefore curvature $4K$ each, which suggests a power law going as $K^{-\log_3 4}$. However, in Figure~\ref{fig:distFit} below we find that this naive estimate does not explain the behaviour seen. The matter of the fact is that the contribution discussed above is measure zero in the limit, since these contributions scale exponentially in $n$, while the total amount of triangles grows like $\binom{3^n}{3}$. Instead, we find that attempting to numerically fit the power-law gives an exponent of $2.08\pm0.18$, and an inverse square law can be seen to be an much better description in Figure~\ref{fig:distFit}. Hence, we see that first order attempts at analytically deriving the limiting distribution gives erroneous results. At this point, it is not obvious whether the limiting distribution for positive $K$ values should indeed be a power law resulting from analytical consideration. Since the full limiting distribution appears to be highly non-trivial to describe, we will instead focus on some of its statistics and how they scale as a function of $n$.

\begin{figure}[H]
    \centering
    \includegraphics[scale=0.2]{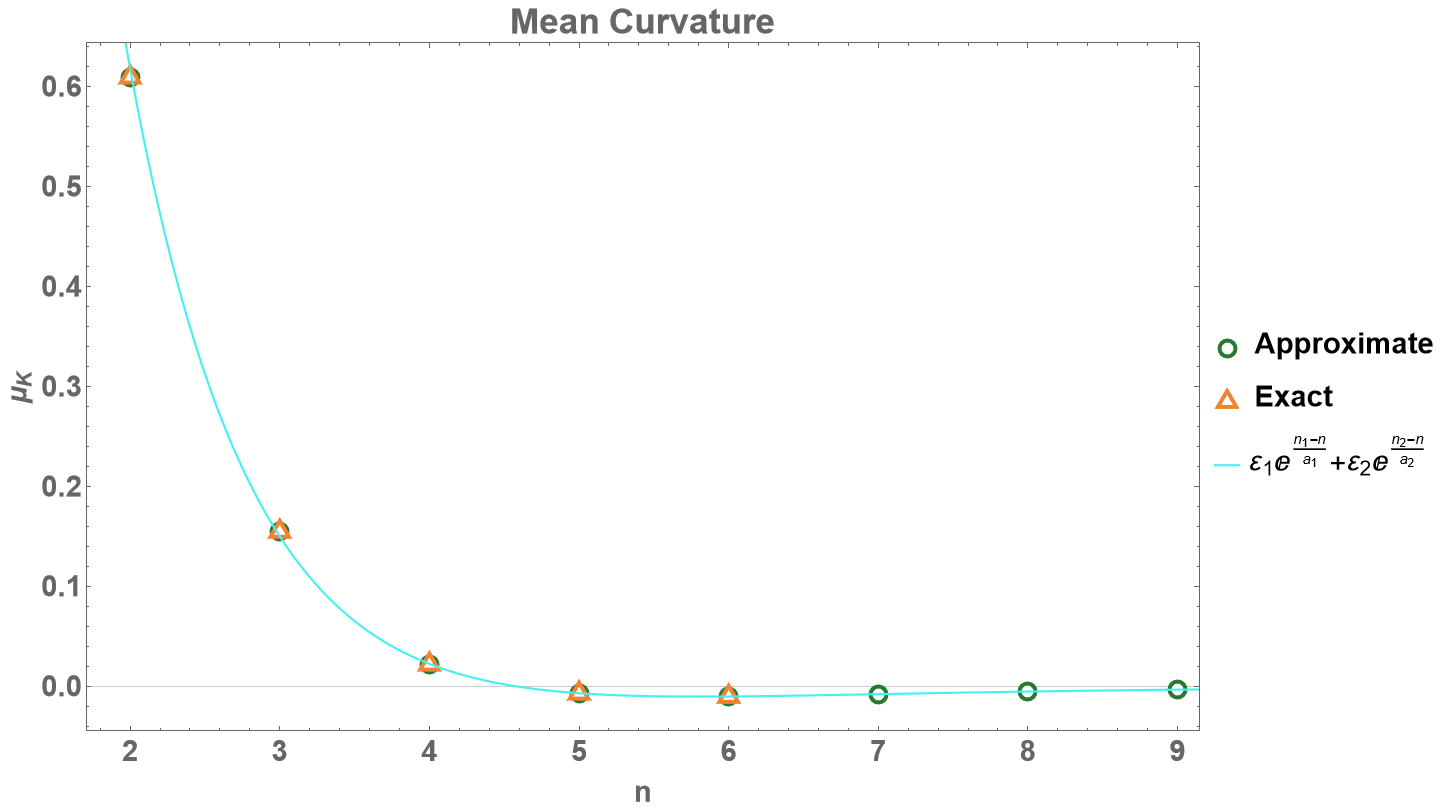}
    \includegraphics[scale=0.2]{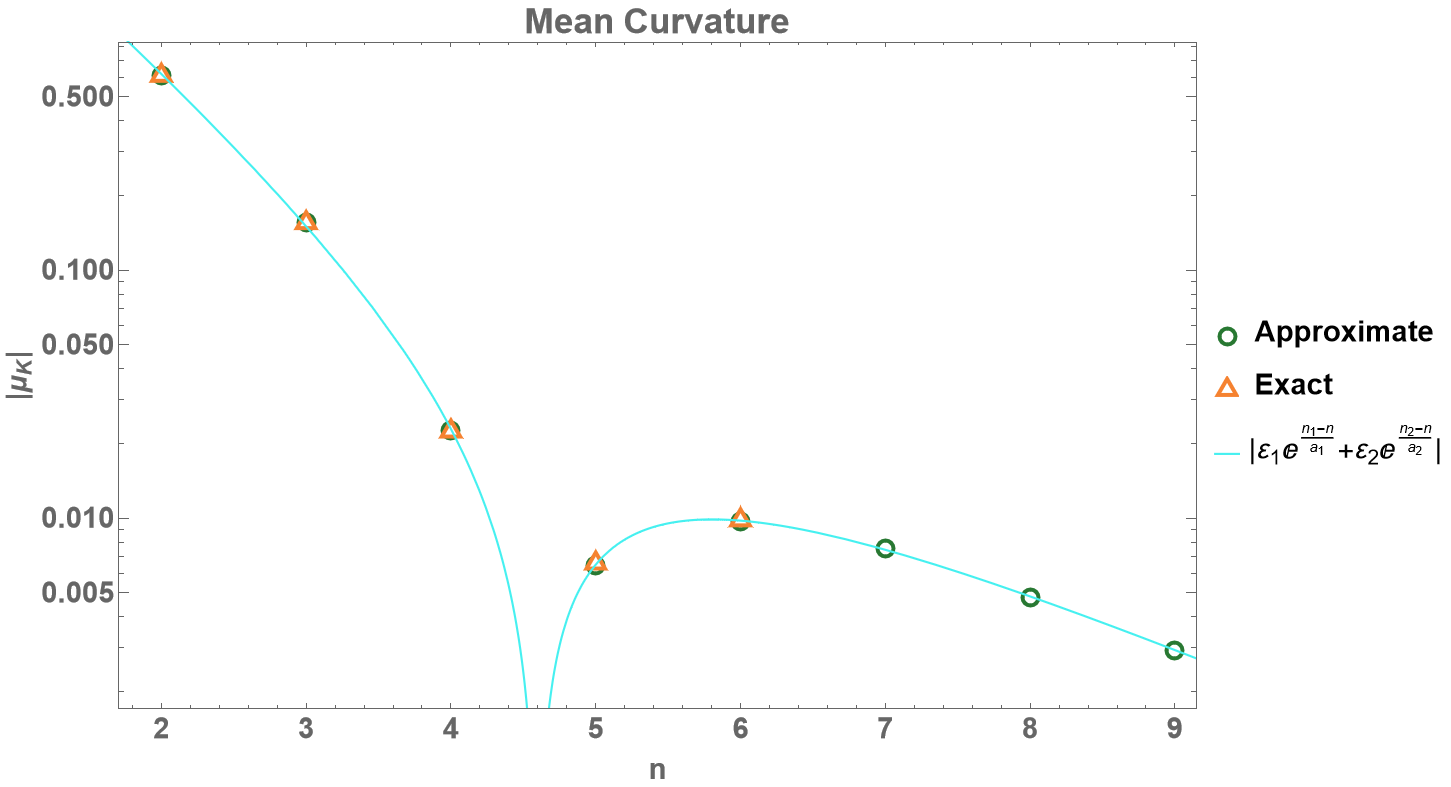}
    \includegraphics[scale=0.2]{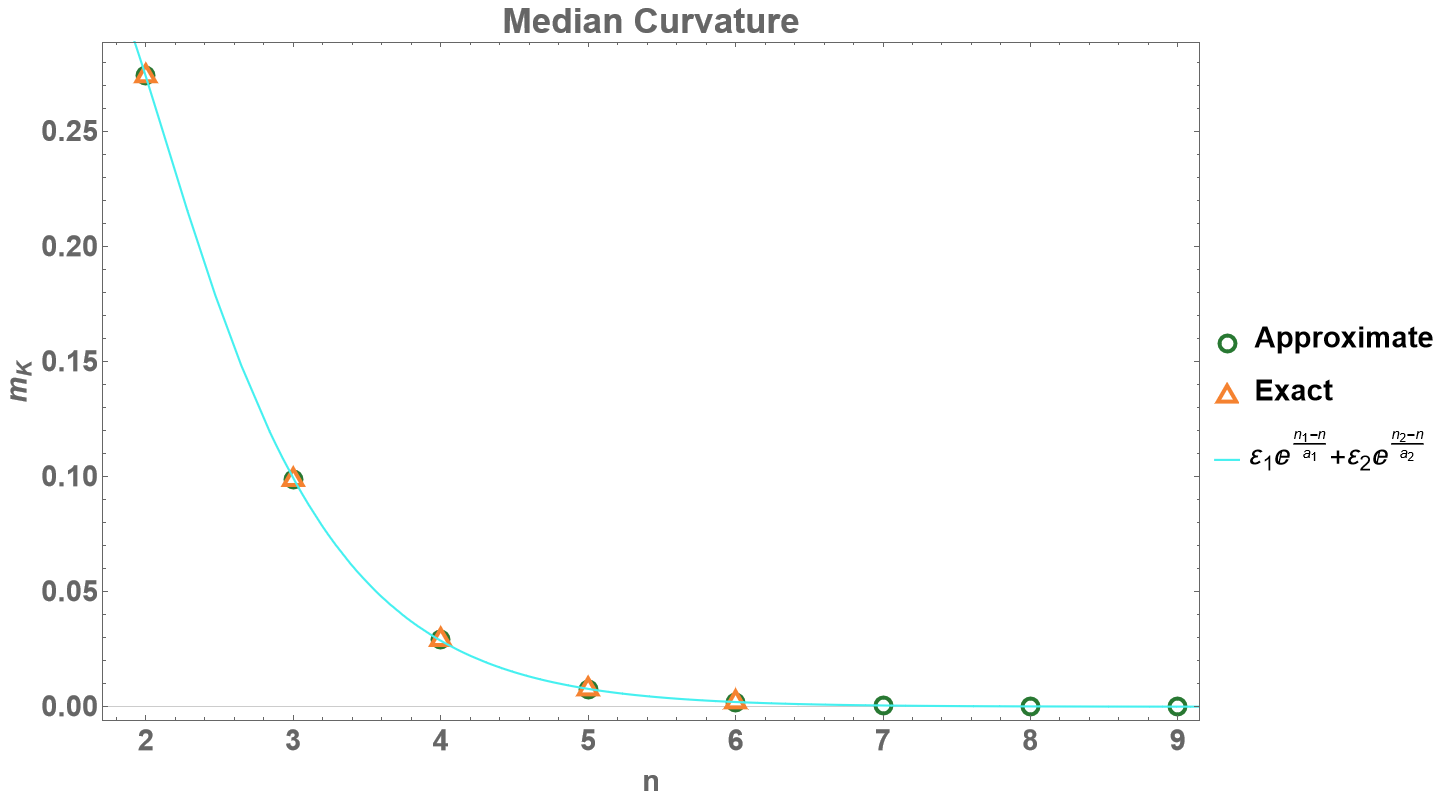}
    \includegraphics[scale=0.2]{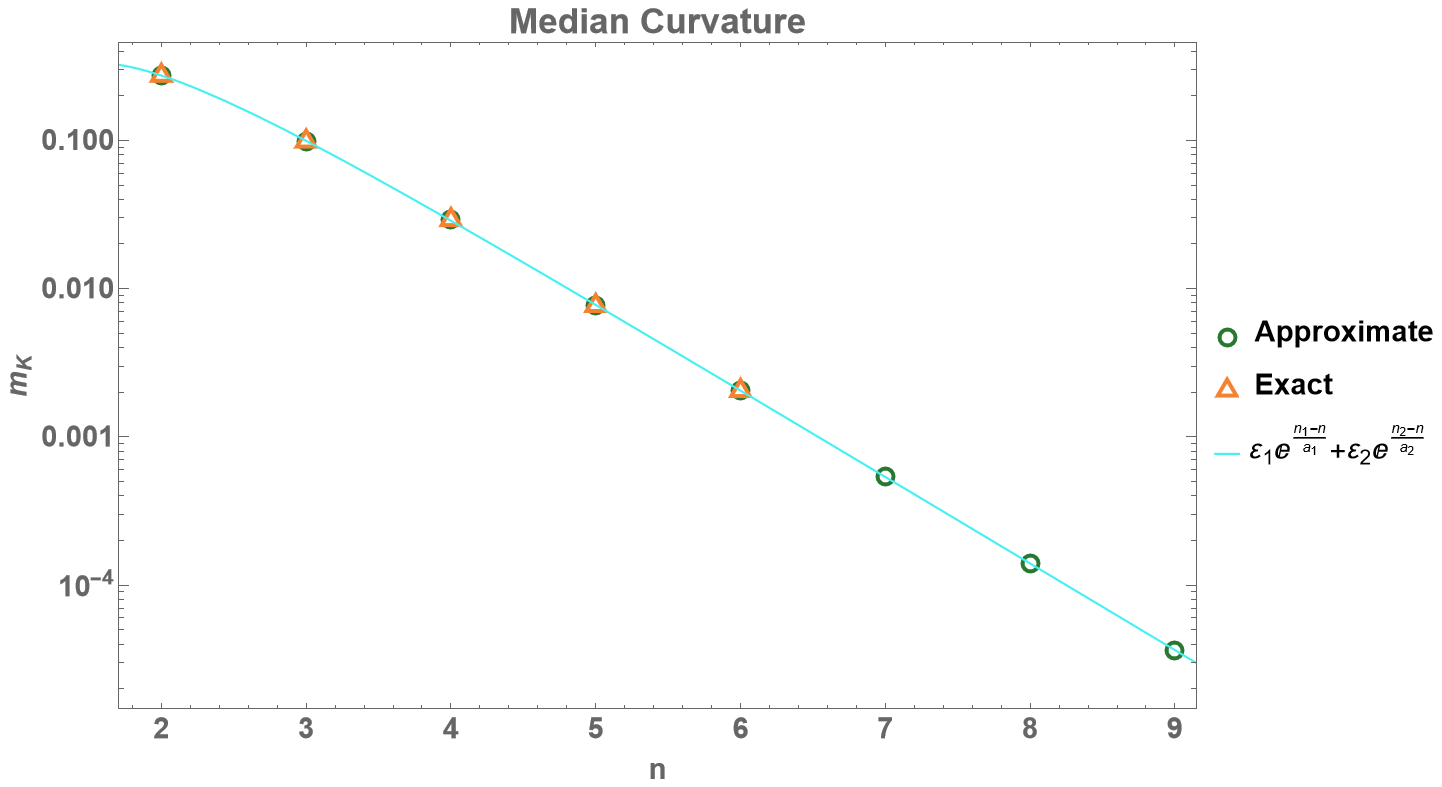}
    \includegraphics[scale=0.2]{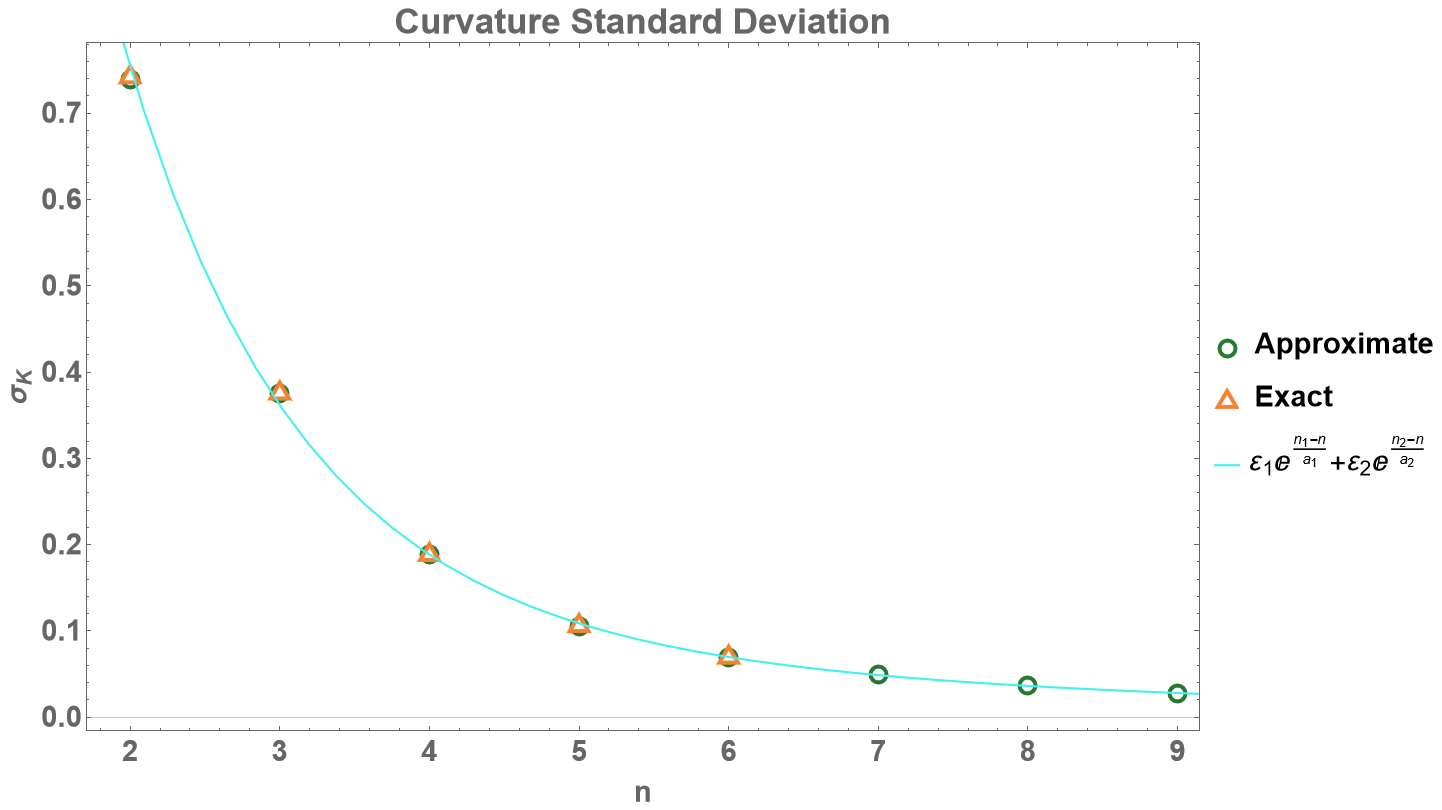}
    \includegraphics[scale=0.2]{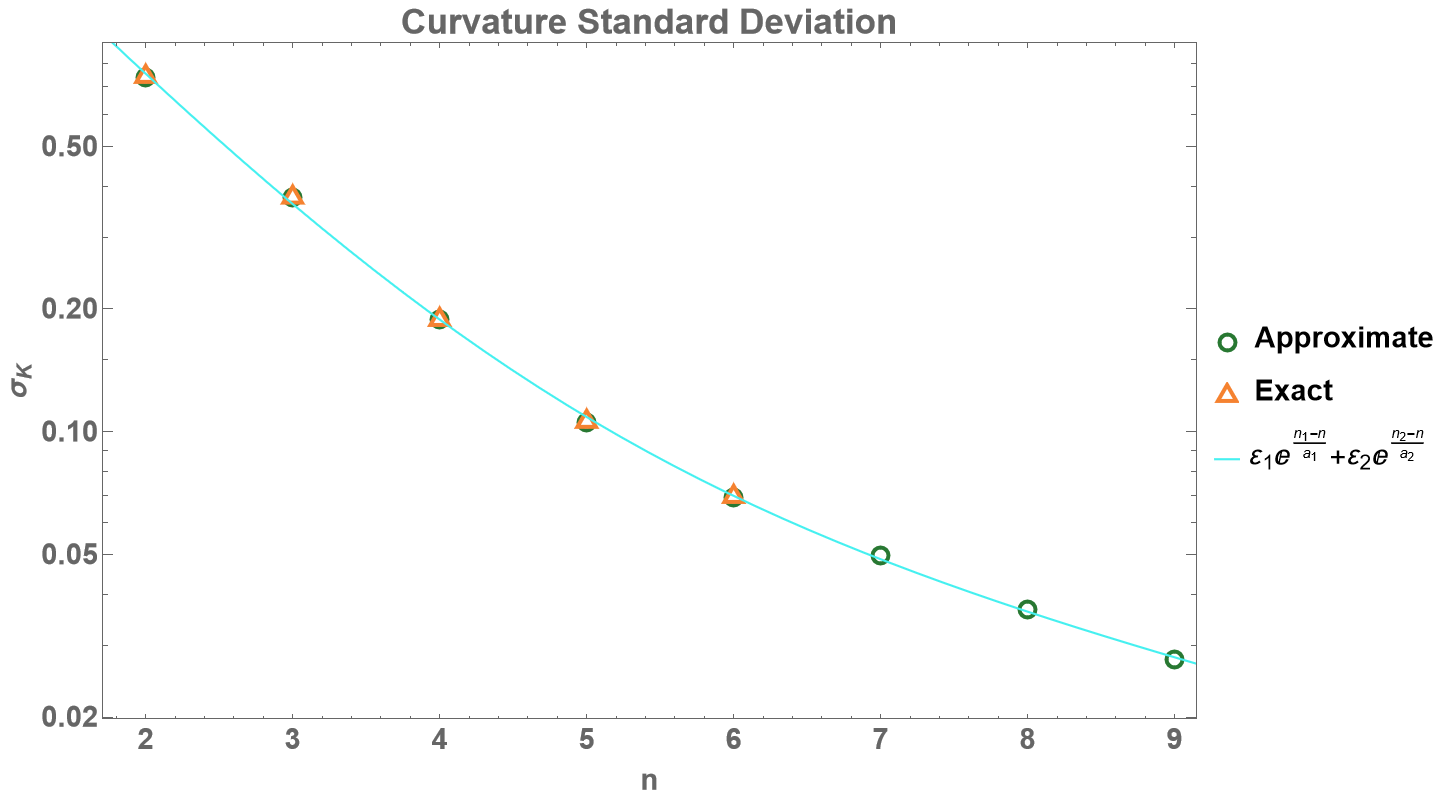}
    \caption{Fits of the scaling of the mean, median and standard deviation of the distribution of sectional curvature of the Sierpinski triangle graphs $\hat{S}_3^n$ to an empirical law. We use both exact data (enumerating all of the finite sub-triangles for a given $n$) as well as approximate data, using rejection sampling to obtain $2\times 10^7$ samples for each $n$.}
    \label{fig:srpMn}
\end{figure}

We find that as a function of $n$ (where the $n^\text{th}$ graph is $\hat{S}_3^n$) the mean, median and standard deviation of the sectional curvature distributions all seem to follow an  empirical law $\varepsilon_1 e^{\frac{n_1-n}{a_1}}+\varepsilon_2 e^{\frac{n_2-n}{a_2}}$ (which is simply taken as an ansatz chosen due to the apparent behavior of the log-plots of the statistics). Here $\varepsilon_i\in\{+1,-1\},\,n_i\in\mathbb{R}$ and $a_i\in\mathbb{R}^+$, and we note that these constants differ between the different statistics. We order the constants such that $a_i$ is decreasing, such that $\epsilon_1$ is the asymptotic sign of the corresponding statistic. The fitted constants in Table~\ref{tb:sierpConsts} are distinct for the mean, median and standard deviation respectively, and is expected to be distinct for different fractals, giving a potentially new way for classifying fractals. If other self-similar fractals give rise to similar scaling laws, one can then, in analogy to critical exponents in statistical physics, use the nature of the scaling to classify such fractals. We note that in contrast to the scale-free property of critical systems studied in statistical physics, here, we are dealing with self-similarity, the discrete form of scale invariance. This  might potentially be interesting for application of statistical physics methods to self-similar fractals. 

\begin{table}[h]
\centering
\begin{tabular}{||c|c|c|c|c|c|c||} 
 \hline
  & $\varepsilon_1$ & $n_1$ & $a_1$ & $\varepsilon_2$ & $n_2$ & $a_2$ \\ [0.5ex] 
 \hline\hline
  Mean & $-1$ & $-1.23\pm0.0.32$ & $1.78\pm0.06$ & $+1$ & $1.788\pm0.016$ & $0.856\pm0.019$  \\
 \hline
 Median & $+1$ & $1.41\pm0.06$ & $0.743\pm0.007$ & $-1$ & $1.24\pm0.05$ & $0.44\pm0.06$  \\
 \hline
  Standard Deviation & $+1$ & $-8\pm 9$ & $4.8\pm 2.4$ & $+1$ & $1.48\pm 0.13$& $1.15\pm0.23$\\ 
 \hline
\end{tabular}
\caption{Fitted constants related to the sequence of distributions of sectional curvature of Sierpinski triangle graphs}
\label{tb:sierpConsts}
\end{table}

\begin{figure}[h]
    \centering
    \includegraphics[scale=0.2]{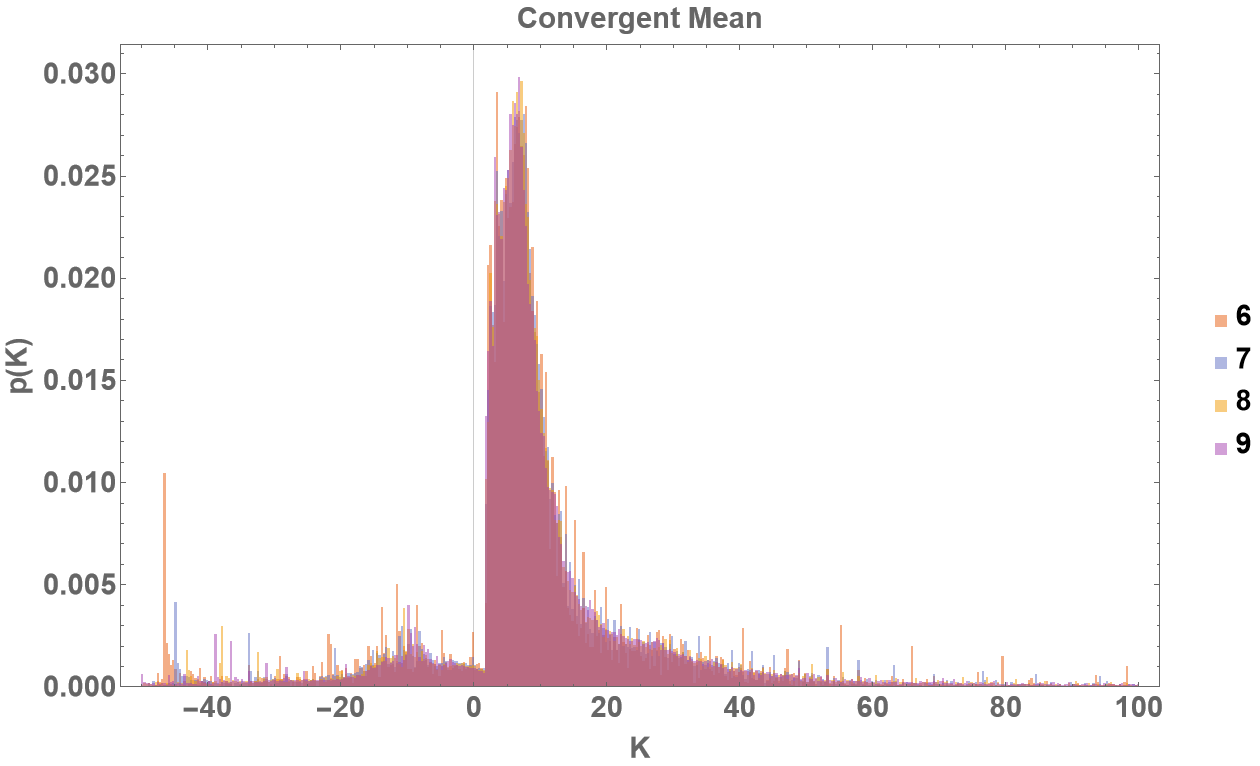}\\
    \includegraphics[scale=0.2]{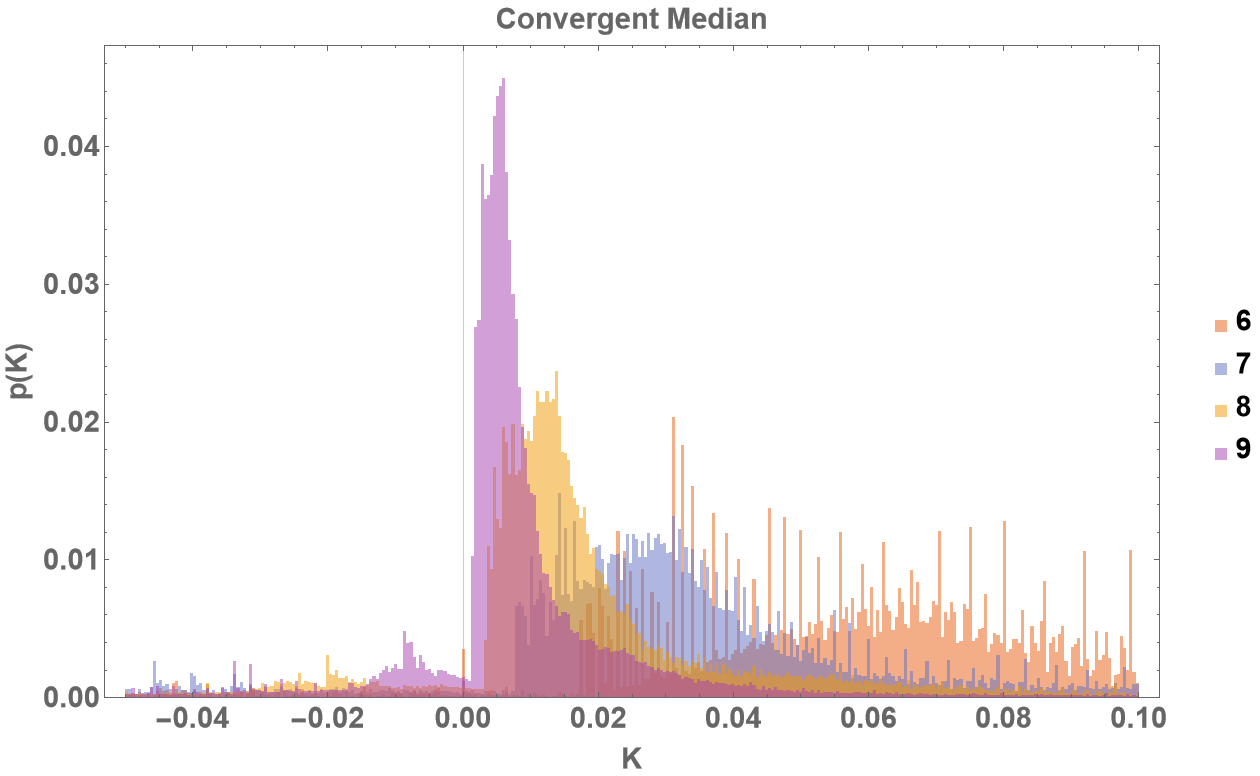}
    \includegraphics[scale=0.2]{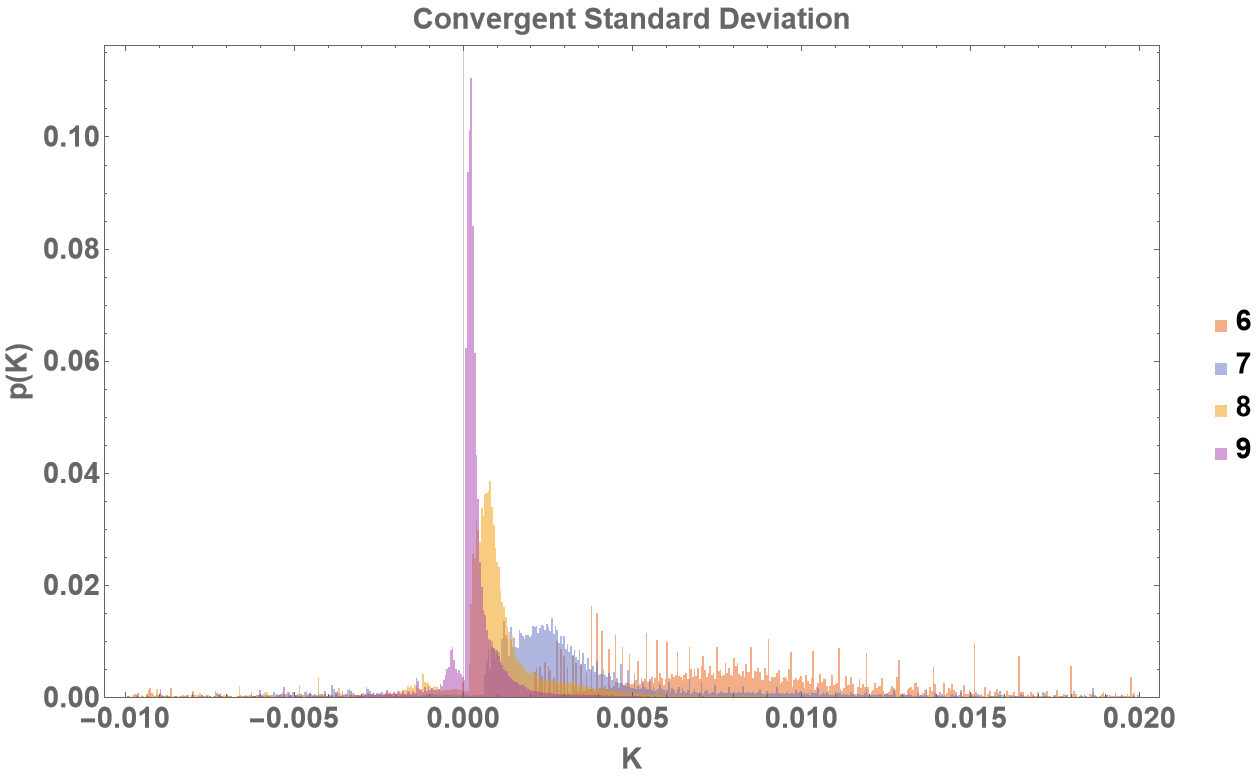}
    \caption{Sectional curvature distributions converging to finite mean, median and standard deviation respectively}
    \label{fig:shftdDistsl}
\end{figure}

Figure~\ref{fig:srpMn} shows the exponential nature of the scaling by comparing it to fits of the empirical law discussed above to the data. Note that we have chosen edge lengths of $1$ at the start, which is an arbitrary choice. If we instead had edge lengths of $\ell^n$ (multiplying the edge length by $\ell$ each iteration), the total side length of the Sierpinski triangle is $(2\ell)^n$, since the amount of edges in the side doubles each iteration. Since sectional curvature has dimensions of length$^{-2}$, the sectional curvature (distribution) would then be scaled by $\ell^{-2n}$. We can see then that for any single one of the measures investigated above we can choose $\ell=e^{-\frac{1}{2a_1}}$, such that the chosen measure would then converge to $\varepsilon_1 e^{\frac{n_1}{a_1}}$ under the assumption that the large $n$ scaling behaviour we observe continues (which remains to be rigorously shown, although there is seemingly no reason for it to break). As long as $\frac12<\ell<1$, the sequence would still geometrically limit to the scale free fractal we desire. This corresponds to the requirement that $a_1>\frac{1}{2\ln(2)}\approx0.721$, that our investigated quantities seem to adhere to. If we choose $\ell\approx0.90$ we seemingly get a limiting distribution with finite standard deviation and vanishing mean and median. If we choose $\ell\approx0.755$ we seemingly get a limiting distribution with finite (negative) mean, diverging standard deviation and vanishing median. Finally if we choose $\ell\approx0.510$ we seemingly get a distribution with finite median and diverging mean and standard deviation. These different scaled distributions can be seen in Figure~\ref{fig:shftdDistsl}. Taken together, the statistics of the discrete sectional curvature distributions of self-similar fractals potentially suggest a new avenue for the study and classification of such fractals.

\vspace{3cm}

\section{Conclusions and Discussion}
\label{sec:conclusion}

The discrete sectional curvature estimator we have proposed in this paper is applicable to any path metric space with approximately constant sectional curvature (or at least on approximately constant sectional curvature regions of such spaces). This makes it a promising candidate of discrete sectional curvature, applicable to graphs and other discrete structures endowed with a suitable path metric. We then performed extensive validation tests for this curvature estimator using hard annulus random geometric graphs corresponding to manifolds of positive, negative and zero curvature. 

For positive curvatures, we find good agreement between the estimated graph curvature and the curvature of the unit sphere; for both, the two and three dimensional spheres: $S^2$ respectively  $S^3$.  As can be seen in Figure \ref{fig:sphrErDist}, for a fixed minimum length scale both the mean absolute errors of the curvature estimates and the absolute error of the mean curvature converges convincingly close to zero as the distortion vanishes. Both plots have large $r^2$ values, indicating near perfect linear correlation between the estimate error and metric distortion in the region tested (even though there may still be small non-linear effects between curvature errors and metric distortion at extremely low metric-distortions, below the limits where our constructions are performed). Similar results have been found in the negative and zero curvature case. The relatively (to the positive curvature cases) lower $r^2$ values in Figures \ref{fig:hypErDist} and \ref{fig:eucErDist} are presumably due to reduced robustness in estimating negative or zero curvature, compared to the positive case, as noted earlier in Figure \ref{fig:curv}. However, in the zero curvature case, an additional subtlety should also be noted: due to the lack of definitive length scale in the euclidean plane, there is no way to reasonably compare errors, since the y-axis can arbitrarily be re-scaled.

To show that our method works as intended on real-world data sampled from an underlying   curved structure, we put our estimator to the test by using it to estimate the radius of the earth   (modeled as an oblate spheroid) using random geographic locations. Even more so, this example neatly demonstrates an extension of our curvature estimation method that can be applied to discrete structures with non-constant sectional curvatures (as in the case with the geometry of an oblate spheroid).  In Figure \ref{fig:earthRadiusDist10000} we see that with no maximum distance the average radius is accurately estimated, but the distribution is not the same as that of the underlying manifold. However, upon introducing a maximum length scale we find two things: firstly, the accuracy of the average curvature estimate is left intact; and additionally this yields the correct distribution of curvatures of an oblate spheroid with parameters comparable to that of earth. This shows that a specific generalisation of the constant sectional curvature estimation methods described in this paper to cases of graphs with  non-constant curvature via the addition of an appropriately chosen maximum length scale. In future work, we will investigate algorithms for deriving the maximum length scale in more generic cases of geometries with varying sectional curvatures.

A relevant future application of our discrete sectional curvature estimator would be to investigate the geometry of real-world complex networks from biology to technology and sociology. Presumably, this will reveal global geometric properties or invariants in real-world data that are otherwise difficult to compute using purely local graph-theoretic measures (for examples of such applications, see   \cite{wu2015emergent},  \cite{arsiwalla2016global, arsiwalla2016high, arsiwalla2017spectral, betella2014brainx3}, \cite{mulder2018network},  \cite{boguna2021network}).  Yet another application of our work here, may be for estimating dimensions of graphs and discrete structures. For instance,  one can invert the formula used in the  Wolfram-Ricci curvature (see Equation~\ref{eq:WR}) to instead get an estimate of dimension. There one can use our estimated average sectional curvature as input and thereby obtain an estimate of the dimension of the graph. This seems like a natural extension to existing dimension estimators   \cite{shanker_defining_2007}. If the sectional curvatures and dimensions of random geometric graphs can be determined, then the Ricci and Riemann curvatures can be computed using an averaging process as detailed in \ref{sec:Riem}. These computations will be explored in future work. 

We have shown how this view of sectional curvature opens up a new avenue in the study of fractals by showing convincing numerical evidence of highly non-trivial behaviour of the limiting distribution of sectional curvatures on Sierpinski triangle graphs. Analytical results supporting and extending the numerical results presented here will definitely be needed and will be pursued in future work. This work could potentially yield a new classification of self-similar fractals, by investigating the scaling behavior of the discrete sectional curvature distribution of such fractals.

We anticipate that the discrete sectional curvature estimator, measure of metric distortion as well as the graph sprinkling algorithm developed in this work will be useful for both applied as well as basic science research. For example, the algorithm generating low metric-distortion graphs presented here may have potential uses in computational geometry technologies that traditionally use mesh-generation techniques. Additionally, the reliability and controlled error margins of our curvature estimates on random geometric graphs, make this a particularly useful tool in data science, where one may want to infer global geometric properties from large data-sets across domains. In fundamental physics, these methods are directly relevant to discrete models of gravity. For instance, the low metric-distortion graphs generated here improve upon traditional graph sprinkling methods used in Causal Set theory  \cite{sorkin_forks_1997,surya2019causal}. Additionally, the discrete curvature definition and estimator presented here may provide alternative formulations of gravitational path integrals to those introduced in \cite{trugenberger_random_2016,trugenberger_combinatorial_2017,trugenberger2021emergent}. The estimation of the Ricci tensor and Riemann tensor on discrete spaces may be useful for a discretized formulation of the Einstein-Hilbert action  \cite{Wolfram2002a},  \cite{Wolfram2020}, \cite{gorard2020some}, and with it Einstein's field equations (once this method has also been generalised to a Lorentzian setting), similar to what is done in Regge theory \cite{regge1961general}. 

Furthermore, let us note that although the discussion in this paper has largely been about graphs, our arguments straightforwardly extend to hypergraphs (or any metric space capable of constructing right-angled triangles). This can be seen by realizing that one can always work with the 1-skeleton of a hypergraph and apply our estimator (with additional constraints on the filling volumes). This will be useful for investigating  discrete models of gravity such as causal dynamical triangulation  \cite{loll2019quantum} and constructions  of pregeometric spaces in the Wolfram model~\cite{arsiwalla2021pregeometric, arsiwalla2021homotopies, arsiwalla2020homotopic}. Besides models of quantum gravity, computations of discrete geometric quantities, as the sectional curvature discussed in this paper, may also be a useful tool for estimating global properties of quantum informational systems such as tensor networks built from ZX operators~\cite{gorard2020zx, gorard2021zx}; or, operator algebraic frameworks associated to non-classical geometries and higher structures \cite{raptis_quantization_1999,arsiwalla2022operator,zapata2022invitation}.


\clearpage
\section*{Appendix}
\appendix
\addcontentsline{toc}{section}{Appendix \ref{apn:proof}}
\addtocontents{toc}{\setcounter{tocdepth}{-1}}
\section{Uniqueness of Roots Corresponding to Sectional Curvature}
\label{apn:proof}
We want to argue that the definition of sectional curvature in this manuscript is well-defined in the sense that the considered function only has one root that can be considered the sectional curvature.
If $a,b,c$ satisfy the triangle inequalities, then we need to show that $f(x)=\cos(c\sqrt{x})-\cos(a\sqrt{x})\cos(b\sqrt{x})$ has exactly two zeroes in $x\in\left(-\infty,\frac{\pi^2}{\max{(a,b,c)}^2}\right]$ (one being the trivial $x=0$ root), unless $a^2+b^2=c^2$ in which case $x=0$ is the only (double) root.
\subsection{For Negative Arguments: $x<0$}
Let $x=-t^2$, then for $x<0$ we have $f(-t^2)=\cosh(c t)-\cosh(a t)\cosh(b t)$ for $t>0$. We can rewrite this as $$f(-t^2)=\cosh(c t)-\frac{\cosh((a+b)t)}{2}-\frac{\cosh((a-b)t)}{2}.$$
This has roots where 
$$\cosh(c t)=\frac{\cosh((a+b)t)}{2}+\frac{\cosh((a-b)t)}{2}$$

\subsubsection{When $c^2\leq a^2+b^2$}
We can note that the Maclaurin series of the left-hand-side will have factors of $c^{2k}$ with the same numeric factors as the right hand side has for $\frac12((a+b)^{2k}+(a-b)^{2k})$. Now using the power mean inequality we can let $M_{2k}=\left[\frac12((a+b)^{2k}+(a-b)^{2k})\right]^\frac{1}{2k}$ such that $M_2=\sqrt{a^2+b^2}$. We then know that $M_{2k}> M_2$, so if $c^2\leq a^2+b^2=M_2^2$ then $c^{2k}\leq M_2^{2k}< M_{2k}^{2k}=\frac12((a+b)^{2k}+(a-b)^{2k})$. Therefore each term on the right-hand-side is strictly larger than the left-hand-side, and the functions are equal for $t=0$, so there are no roots for $t>0$, meaning $f(x)$ has no roots for $x<0$ if  $c^2\leq a^2+b^2$.\\

\subsubsection{When  $c^2> a^2+b^2$}
We can now consider the case where $c^2>a^2+b^2$ and $c<a+b$ (triangle inequality). We can now assume without loss of generality (WLOG) that $c=1$ such that $a^2+b^2<1$ and $a+b>1$. In the same spirit as above, we now want to compare $c^{2k}=1$ to $\frac12((a+b)^{2k}+(a-b)^{2k})$, since this is what the terms in the Maclaurin series will be proportional to. We can note $\frac{d}{dk}((a+b)^{2k}+(a-b)^{2k})=(a+b)^{2k}\ln((a+b)^2)+(a-b)^{2k}\ln((a-b)^2)$. Here we have $a+b>1$ as well as $(a-b)^2=a^2+b^2-2a b<1$, so the first term in the derivative is strictly positive and increases with $k$. Similarly, the second term will be strictly negative, but decrease in magnitude as $k$ increases. Therefore there is some $m$ such that  for all $k>m$ $\frac12((a+b)^{2k}+(a-b)^{2k})$ is positive. So $\frac12((a+b)^{2k}+(a-b)^{2k})$ starts as 1 when $k=0$, then decreases for $k<m$, before which it starts increasing again for $k>m$. So there exists an $n>m$, such that $\frac12((a+b)^{2k}+(a-b)^{2k})>1$ for $k>n$. This means that in this case the first few terms in the Maclaurin series of the left-hand-side is larger than the corresponding terms on the right, but for higher order terms the right will be larger, so after some point $\frac{\cosh((a+b)t)}{2}+\frac{\cosh((a-b)t)}{2}$ grows faster than $\cosh(t)$, meaning there will be exactly one root.

\subsection{For Positive Arguments: $x>0$}
For positive arguments it is more difficult to conclusively show uniqueness, due to the very specific interval to be considered. We can however consider the following argument. Notice that  $\frac{d^2f(t^2)}{dt^2}|_{t=0}=a^2+b^2-c^2$, and that the first derivative of this is simply zero. So we know for some small $\epsilon>0$ that $\text{sign}(f(\epsilon))=\text{sign}(a^2+b^2-c^2)$. We can now assume without loss of generality that $a\geq b$.\\
If $c\geq a$ then $f\left(\left(\frac{\pi}{c}\right)^2\right)=-1-\cos(\frac{a}{c}\pi)\cos(\frac{b}{c}\pi)\leq 0$.
\\Otherwise if $c<a$ then $f\left(\left(\frac{\pi}{a}\right)^2\right)=\cos(\frac{c}{a}\pi)+\cos(\frac{b}{a}\pi)$\\
We have $\frac{c}{a}\pi< \pi$ and $\frac{b}{a}\pi\leq\pi$. We also know that $c+b\geq a$, so $\frac{c}{a}\pi+\frac{b}{a}\pi\geq\pi$\\
So if $a^2+b^2>c^2$ then there has to be an even amount of roots in the interval $(0,\frac{\pi}{\max(a,b,c)}]$, and conversely if $a^2+b^2<c^2$ then there has to be an odd amount of roots in the same interval. Beyond this we currently only have strong numerical evidence demonstrating a unique root in this region.

\section*{Acknowledgments}

The authors would like to thank Stephen Wolfram and Jonathan Gorard for suggesting this project and useful discussions. Additionally, we gratefully acknowledge Dmitri Krioukov,  William J. Cunningham, Carlo Trugenberger and Juergen Jost for useful suggestions. JFDP would also like to thank Adri Wessels, Ralph McDougall and Jean Weight for several productive discussions.

\bibliographystyle{eptcs}
\bibliography{refs.bib}

\begin{thebibliography}{10}
\providecommand{\bibitemdeclare}[2]{}
\providecommand{\surnamestart}{}
\providecommand{\surnameend}{}
\providecommand{\urlprefix}{Available at }
\providecommand{\url}[1]{\texttt{#1}}
\providecommand{\href}[2]{\texttt{#2}}
\providecommand{\urlalt}[2]{\href{#1}{#2}}
\providecommand{\doi}[1]{doi:\urlalt{http://dx.doi.org/#1}{#1}}
\providecommand{\bibinfo}[2]{#2}

\bibitemdeclare{article}{abraham_advances_nodate}
\bibitem{abraham_advances_nodate}
\bibinfo{author}{Ittai \surnamestart Abraham\surnameend}, \bibinfo{author}{Yair
  \surnamestart Bartal\surnameend} \& \bibinfo{author}{Ofer \surnamestart
  Neiman\surnameend}: \emph{\bibinfo{title}{Advances in Metric Embedding
  Theory}}, p.~\bibinfo{pages}{16}.

\bibitemdeclare{article}{arsiwalla2020homotopic}
\bibitem{arsiwalla2020homotopic}
\bibinfo{author}{Xerxes~D \surnamestart Arsiwalla\surnameend}
  (\bibinfo{year}{2020}): \emph{\bibinfo{title}{Homotopic Foundations of
  Wolfram Models}}.
\newblock {\sl \bibinfo{journal}{Wolfram Community. https://community. wolfram.
  com/groups/-/m}} \bibinfo{volume}{2032113}.

\bibitemdeclare{article}{arsiwalla2022operator}
\bibitem{arsiwalla2022operator}
\bibinfo{author}{Xerxes~D \surnamestart Arsiwalla\surnameend},
  \bibinfo{author}{David \surnamestart Chester\surnameend} \&
  \bibinfo{author}{Louis~H \surnamestart Kauffman\surnameend}
  (\bibinfo{year}{2022}): \emph{\bibinfo{title}{On the Operator Origins of
  Classical and Quantum Wave Functions}}.
\newblock {\sl \bibinfo{journal}{arXiv preprint arXiv:2211.01838}}.

\bibitemdeclare{article}{arsiwalla2021pregeometric}
\bibitem{arsiwalla2021pregeometric}
\bibinfo{author}{Xerxes~D \surnamestart Arsiwalla\surnameend} \&
  \bibinfo{author}{Jonathan \surnamestart Gorard\surnameend}
  (\bibinfo{year}{2021}): \emph{\bibinfo{title}{Pregeometric Spaces from
  Wolfram Model Rewriting Systems as Homotopy Types}}.
\newblock {\sl \bibinfo{journal}{arXiv preprint arXiv:2111.03460}}.

\bibitemdeclare{article}{arsiwalla2021homotopies}
\bibitem{arsiwalla2021homotopies}
\bibinfo{author}{Xerxes~D \surnamestart Arsiwalla\surnameend},
  \bibinfo{author}{Jonathan \surnamestart Gorard\surnameend} \&
  \bibinfo{author}{Hatem \surnamestart Elshatlawy\surnameend}
  (\bibinfo{year}{2021}): \emph{\bibinfo{title}{Homotopies in Multiway
  (Non-Deterministic) Rewriting Systems as $ n $-Fold Categories}}.
\newblock {\sl \bibinfo{journal}{arXiv preprint arXiv:2105.10822}}.

\bibitemdeclare{article}{arsiwalla2017spectral}
\bibitem{arsiwalla2017spectral}
\bibinfo{author}{Xerxes~D \surnamestart Arsiwalla\surnameend},
  \bibinfo{author}{Pedro~AM \surnamestart Mediano\surnameend} \&
  \bibinfo{author}{Paul~FMJ \surnamestart Verschure\surnameend}
  (\bibinfo{year}{2017}): \emph{\bibinfo{title}{Spectral modes of network
  dynamics reveal increased informational complexity near criticality}}.
\newblock {\sl \bibinfo{journal}{Procedia Computer Science}}
  \bibinfo{volume}{108}, pp. \bibinfo{pages}{119--128}.

\bibitemdeclare{article}{arsiwalla2016global}
\bibitem{arsiwalla2016global}
\bibinfo{author}{Xerxes~D \surnamestart Arsiwalla\surnameend} \&
  \bibinfo{author}{Paul~FMJ \surnamestart Verschure\surnameend}
  (\bibinfo{year}{2016}): \emph{\bibinfo{title}{The global dynamical complexity
  of the human brain network}}.
\newblock {\sl \bibinfo{journal}{Applied network science}}
  \bibinfo{volume}{1}(\bibinfo{number}{1}), pp. \bibinfo{pages}{1--13}.

\bibitemdeclare{inproceedings}{arsiwalla2016high}
\bibitem{arsiwalla2016high}
\bibinfo{author}{Xerxes~D \surnamestart Arsiwalla\surnameend} \&
  \bibinfo{author}{Paul~FMJ \surnamestart Verschure\surnameend}
  (\bibinfo{year}{2016}): \emph{\bibinfo{title}{High integrated information in
  complex networks near criticality}}.
\newblock In: {\sl \bibinfo{booktitle}{International Conference on Artificial
  Neural Networks}}, \bibinfo{organization}{Springer}, pp.
  \bibinfo{pages}{184--191}.

\bibitemdeclare{inproceedings}{betella2014brainx3}
\bibitem{betella2014brainx3}
\bibinfo{author}{Alberto \surnamestart Betella\surnameend},
  \bibinfo{author}{Ryszard \surnamestart Cetnarski\surnameend},
  \bibinfo{author}{Riccardo \surnamestart Zucca\surnameend},
  \bibinfo{author}{Xerxes~D \surnamestart Arsiwalla\surnameend},
  \bibinfo{author}{Enrique \surnamestart Martinez\surnameend},
  \bibinfo{author}{Pedro \surnamestart Omedas\surnameend},
  \bibinfo{author}{Anna \surnamestart Mura\surnameend} \&
  \bibinfo{author}{Paul~FMJ \surnamestart Verschure\surnameend}
  (\bibinfo{year}{2014}): \emph{\bibinfo{title}{BrainX3: embodied exploration
  of neural data}}.
\newblock In: {\sl \bibinfo{booktitle}{Proceedings of the 2014 virtual reality
  international conference}}, pp. \bibinfo{pages}{1--4}.

\bibitemdeclare{article}{boguna2021network}
\bibitem{boguna2021network}
\bibinfo{author}{Marian \surnamestart Boguna\surnameend}, \bibinfo{author}{Ivan
  \surnamestart Bonamassa\surnameend}, \bibinfo{author}{Manlio \surnamestart
  De~Domenico\surnameend}, \bibinfo{author}{Shlomo \surnamestart
  Havlin\surnameend}, \bibinfo{author}{Dmitri \surnamestart
  Krioukov\surnameend} \& \bibinfo{author}{M~\surnamestart Serrano\surnameend}
  (\bibinfo{year}{2021}): \emph{\bibinfo{title}{Network geometry}}.
\newblock {\sl \bibinfo{journal}{Nature Reviews Physics}}
  \bibinfo{volume}{3}(\bibinfo{number}{2}), pp. \bibinfo{pages}{114--135}.

\bibitemdeclare{inproceedings}{chennuru_vankadara_measures_2018}
\bibitem{chennuru_vankadara_measures_2018}
\bibinfo{author}{Leena \surnamestart Chennuru~Vankadara\surnameend} \&
  \bibinfo{author}{Ulrike \surnamestart von Luxburg\surnameend}:
  \emph{\bibinfo{title}{Measures of distortion for machine learning}}.
\newblock In \bibinfo{editor}{S.~\surnamestart Bengio\surnameend},
  \bibinfo{editor}{H.~\surnamestart Wallach\surnameend},
  \bibinfo{editor}{H.~\surnamestart Larochelle\surnameend},
  \bibinfo{editor}{K.~\surnamestart Grauman\surnameend},
  \bibinfo{editor}{N.~\surnamestart Cesa-Bianchi\surnameend} \&
  \bibinfo{editor}{R.~\surnamestart Garnett\surnameend}, editors: {\sl
  \bibinfo{booktitle}{Advances in Neural Information Processing Systems}},
  \bibinfo{volume}{31}, \bibinfo{publisher}{Curran Associates, Inc.}
\newblock
  \urlprefix\url{https://proceedings.neurips.cc/paper/2018/file/4c5bcfec8584af0d967f1ab10179ca4b-Paper.pdf}.

\bibitemdeclare{incollection}{denne_convergence_2008}
\bibitem{denne_convergence_2008}
\bibinfo{author}{Elizabeth \surnamestart Denne\surnameend} \&
  \bibinfo{author}{John~M. \surnamestart Sullivan\surnameend}:
  \emph{\bibinfo{title}{Convergence and Isotopy Type for Graphs of Finite Total
  Curvature}}.
\newblock In \bibinfo{editor}{Alexander~I. \surnamestart Bobenko\surnameend},
  \bibinfo{editor}{John~M. \surnamestart Sullivan\surnameend},
  \bibinfo{editor}{Peter \surnamestart Schröder\surnameend} \&
  \bibinfo{editor}{Günter~M. \surnamestart Ziegler\surnameend}, editors: {\sl
  \bibinfo{booktitle}{Discrete Differential Geometry}},
  \bibinfo{series}{Oberwolfach Seminars}, \bibinfo{publisher}{Birkhäuser}, pp.
  \bibinfo{pages}{163--174}, \doi{10.1007/978-3-7643-8621-4\_8}.
\newblock \urlprefix\url{https://doi.org/10.1007/978-3-7643-8621-4\_8}.

\bibitemdeclare{article}{dettmann_random_2016}
\bibitem{dettmann_random_2016}
\bibinfo{author}{Carl~P. \surnamestart Dettmann\surnameend} \&
  \bibinfo{author}{Orestis \surnamestart Georgiou\surnameend}:
  \emph{\bibinfo{title}{Random geometric graphs with general connection
  functions}} \bibinfo{volume}{93}(\bibinfo{number}{3}).
\newblock \doi{10.1103/physreve.93.032313}.
\newblock \urlprefix\url{https://dx.doi.org/10.1103/physreve.93.032313}.
\newblock \bibinfo{note}{Publisher: American Physical Society ({APS})}.

\bibitemdeclare{article}{devriendt_discrete_2022}
\bibitem{devriendt_discrete_2022}
\bibinfo{author}{Karel \surnamestart Devriendt\surnameend} \&
  \bibinfo{author}{Renaud \surnamestart Lambiotte\surnameend}
  (\bibinfo{year}{2022}): \emph{\bibinfo{title}{Discrete curvature on graphs
  from the effective resistance*}}.
\newblock {\sl \bibinfo{journal}{Journal of Physics: Complexity}}
  \bibinfo{volume}{3}(\bibinfo{number}{2}), p. \bibinfo{pages}{025008},
  \doi{10.1088/2632-072X/ac730d}.
\newblock \urlprefix\url{https://dx.doi.org/10.1088/2632-072X/ac730d}.

\bibitemdeclare{article}{eidi_ollivier_2020}
\bibitem{eidi_ollivier_2020}
\bibinfo{author}{Marzieh \surnamestart Eidi\surnameend} \&
  \bibinfo{author}{Jürgen \surnamestart Jost\surnameend}:
  \emph{\bibinfo{title}{Ollivier Ricci curvature of directed hypergraphs}}
  \bibinfo{volume}{10}(\bibinfo{number}{1}).
\newblock \doi{10.1038/s41598-020-68619-6}.
\newblock \urlprefix\url{https://dx.doi.org/10.1038/s41598-020-68619-6}.
\newblock \bibinfo{note}{Publisher: Springer Science and Business Media {LLC}}.

\bibitemdeclare{article}{forman_bochners_2003}
\bibitem{forman_bochners_2003}
\bibinfo{author}{\surnamestart {Forman}\surnameend}:
  \emph{\bibinfo{title}{Bochner's Method for Cell Complexes and Combinatorial
  Ricci Curvature}} \bibinfo{volume}{29}(\bibinfo{number}{3}), pp.
  \bibinfo{pages}{323--374}.
\newblock \doi{10.1007/s00454-002-0743-x}.
\newblock \urlprefix\url{https://doi.org/10.1007/s00454-002-0743-x}.

\bibitemdeclare{article}{gorard2020some}
\bibitem{gorard2020some}
\bibinfo{author}{Jonathan \surnamestart Gorard\surnameend}
  (\bibinfo{year}{2020}): \emph{\bibinfo{title}{Some Relativistic and
  Gravitational Properties of the Wolfram Model}}.
\newblock {\sl \bibinfo{journal}{Complex Systems}}
  \bibinfo{volume}{29}(\bibinfo{number}{2}).

\bibitemdeclare{article}{gorard2020zx}
\bibitem{gorard2020zx}
\bibinfo{author}{Jonathan \surnamestart Gorard\surnameend},
  \bibinfo{author}{Manojna \surnamestart Namuduri\surnameend} \&
  \bibinfo{author}{Xerxes~D \surnamestart Arsiwalla\surnameend}
  (\bibinfo{year}{2020}): \emph{\bibinfo{title}{ZX-Calculus and Extended
  Hypergraph Rewriting Systems I: A Multiway Approach to Categorical Quantum
  Information Theory}}.
\newblock {\sl \bibinfo{journal}{arXiv preprint arXiv:2010.02752}}.

\bibitemdeclare{article}{gorard2021zx}
\bibitem{gorard2021zx}
\bibinfo{author}{Jonathan \surnamestart Gorard\surnameend},
  \bibinfo{author}{Manojna \surnamestart Namuduri\surnameend} \&
  \bibinfo{author}{Xerxes~D \surnamestart Arsiwalla\surnameend}
  (\bibinfo{year}{2021}): \emph{\bibinfo{title}{Zx-calculus and extended
  wolfram model systems II: fast diagrammatic reasoning with an application to
  quantum circuit simplification}}.
\newblock {\sl \bibinfo{journal}{arXiv preprint arXiv:2103.15820}}.

\bibitemdeclare{article}{HINZ2017565}
\bibitem{HINZ2017565}
\bibinfo{author}{Andreas~M. \surnamestart Hinz\surnameend},
  \bibinfo{author}{Sandi \surnamestart Klavžar\surnameend} \&
  \bibinfo{author}{Sara~Sabrina \surnamestart Zemljič\surnameend}
  (\bibinfo{year}{2017}): \emph{\bibinfo{title}{A survey and classification of
  Sierpiński-type graphs}}.
\newblock {\sl \bibinfo{journal}{Discrete Applied Mathematics}}
  \bibinfo{volume}{217}, pp. \bibinfo{pages}{565--600},
  \doi{https://doi.org/10.1016/j.dam.2016.09.024}.
\newblock
  \urlprefix\url{https://www.sciencedirect.com/science/article/pii/S0166218X16304309}.

\bibitemdeclare{article}{PhysRevResearch.3.013211}
\bibitem{PhysRevResearch.3.013211}
\bibinfo{author}{Pim \surnamestart van~der Hoorn\surnameend},
  \bibinfo{author}{William~J. \surnamestart Cunningham\surnameend},
  \bibinfo{author}{Gabor \surnamestart Lippner\surnameend},
  \bibinfo{author}{Carlo \surnamestart Trugenberger\surnameend} \&
  \bibinfo{author}{Dmitri \surnamestart Krioukov\surnameend}
  (\bibinfo{year}{2021}): \emph{\bibinfo{title}{Ollivier-Ricci curvature
  convergence in random geometric graphs}}.
\newblock {\sl \bibinfo{journal}{Phys. Rev. Research}} \bibinfo{volume}{3}, p.
  \bibinfo{pages}{013211}, \doi{10.1103/PhysRevResearch.3.013211}.
\newblock
  \urlprefix\url{https://link.aps.org/doi/10.1103/PhysRevResearch.3.013211}.

\bibitemdeclare{article}{jost_olliviers_2014}
\bibitem{jost_olliviers_2014}
\bibinfo{author}{Jürgen \surnamestart Jost\surnameend} \&
  \bibinfo{author}{Shiping \surnamestart Liu\surnameend}:
  \emph{\bibinfo{title}{Ollivier’s Ricci Curvature, Local Clustering and
  Curvature-Dimension Inequalities on Graphs}}
  \bibinfo{volume}{51}(\bibinfo{number}{2}), pp. \bibinfo{pages}{300--322}.
\newblock \doi{10.1007/s00454-013-9558-1}.
\newblock \urlprefix\url{http://link.springer.com/10.1007/s00454-013-9558-1}.

\bibitemdeclare{article}{kamtue2018combinatorial}
\bibitem{kamtue2018combinatorial}
\bibinfo{author}{Supanat \surnamestart Kamtue\surnameend}
  (\bibinfo{year}{2018}): \emph{\bibinfo{title}{Combinatorial,
  Bakry-$\backslash$'Emery, Ollivier's Ricci curvature notions and their
  motivation from Riemannian geometry}}.
\newblock {\sl \bibinfo{journal}{arXiv preprint arXiv:1803.08898}}.

\bibitemdeclare{article}{lin_ricci_2011}
\bibitem{lin_ricci_2011}
\bibinfo{author}{Yong \surnamestart Lin\surnameend}, \bibinfo{author}{Linyuan
  \surnamestart Lu\surnameend} \& \bibinfo{author}{Shing-Tung \surnamestart
  Yau\surnameend}: \emph{\bibinfo{title}{Ricci curvature of graphs}}
  \bibinfo{volume}{63}(\bibinfo{number}{4}), pp. \bibinfo{pages}{605--627}.
\newblock \doi{10.2748/tmj/1325886283}.
\newblock \urlprefix\url{https://dx.doi.org/10.2748/tmj/1325886283}.
\newblock \bibinfo{note}{Publisher: Mathematical Institute, Tohoku University}.

\bibitemdeclare{article}{loll2019quantum}
\bibitem{loll2019quantum}
\bibinfo{author}{Renate \surnamestart Loll\surnameend} (\bibinfo{year}{2019}):
  \emph{\bibinfo{title}{Quantum gravity from causal dynamical triangulations: a
  review}}.
\newblock {\sl \bibinfo{journal}{Classical and Quantum Gravity}}
  \bibinfo{volume}{37}(\bibinfo{number}{1}), p. \bibinfo{pages}{013002}.

\bibitemdeclare{article}{mel1995analytic}
\bibitem{mel1995analytic}
\bibinfo{author}{Mark~Samuilovich \surnamestart Mel'nikov\surnameend}
  (\bibinfo{year}{1995}): \emph{\bibinfo{title}{Analytic capacity: discrete
  approach and curvature of measure}}.
\newblock {\sl \bibinfo{journal}{Sbornik: Mathematics}}
  \bibinfo{volume}{186}(\bibinfo{number}{6}), p. \bibinfo{pages}{827}.

\bibitemdeclare{article}{mulder2018network}
\bibitem{mulder2018network}
\bibinfo{author}{Daan \surnamestart Mulder\surnameend} \&
  \bibinfo{author}{Ginestra \surnamestart Bianconi\surnameend}
  (\bibinfo{year}{2018}): \emph{\bibinfo{title}{Network geometry and
  complexity}}.
\newblock {\sl \bibinfo{journal}{Journal of Statistical Physics}}
  \bibinfo{volume}{173}(\bibinfo{number}{3}), pp. \bibinfo{pages}{783--805}.

\bibitemdeclare{article}{ollivier_ricci_2007}
\bibitem{ollivier_ricci_2007}
\bibinfo{author}{Yann \surnamestart Ollivier\surnameend}:
  \emph{\bibinfo{title}{Ricci curvature of metric spaces}}
  \bibinfo{volume}{345}(\bibinfo{number}{11}), pp. \bibinfo{pages}{643--646}.
\newblock \doi{10.1016/j.crma.2007.10.041}.
\newblock \urlprefix\url{https://dx.doi.org/10.1016/j.crma.2007.10.041}.
\newblock \bibinfo{note}{Publisher: Elsevier {BV}}.

\bibitemdeclare{book}{penrose_random_2003}
\bibitem{penrose_random_2003}
\bibinfo{author}{Mathew \surnamestart Penrose\surnameend}:
  \emph{\bibinfo{title}{Random Geometric Graphs}}.
\newblock \bibinfo{series}{Oxford Studies in Probability},
  \bibinfo{publisher}{Oxford University Press},
  \doi{10.1093/acprof:oso/9780198506263.001.0001}.
\newblock
  \urlprefix\url{https://oxford.universitypressscholarship.com/10.1093/acprof:oso/9780198506263.001.0001/acprof-9780198506263}.

\bibitemdeclare{article}{raptis_quantization_1999}
\bibitem{raptis_quantization_1999}
\bibinfo{author}{I.~\surnamestart Raptis\surnameend} \&
  \bibinfo{author}{R.~\surnamestart Zapatrin\surnameend}:
  \emph{\bibinfo{title}{Quantization of Discretized Spacetimes and {the
  Correspondence} Principle}}.
\newblock \doi{10.1023/A:1003694830614}.

\bibitemdeclare{article}{regge1961general}
\bibitem{regge1961general}
\bibinfo{author}{Tullio \surnamestart Regge\surnameend} (\bibinfo{year}{1961}):
  \emph{\bibinfo{title}{General relativity without coordinates}}.
\newblock {\sl \bibinfo{journal}{Il Nuovo Cimento (1955-1965)}}
  \bibinfo{volume}{19}(\bibinfo{number}{3}), pp. \bibinfo{pages}{558--571}.

\bibitemdeclare{book}{rosen2012discrete}
\bibitem{rosen2012discrete}
\bibinfo{author}{Kenneth~H \surnamestart Rosen\surnameend} \&
  \bibinfo{author}{Kamala \surnamestart Krithivasan\surnameend}
  (\bibinfo{year}{2012}): \emph{\bibinfo{title}{Discrete mathematics and its
  applications: with combinatorics and graph theory}}.
\newblock \bibinfo{publisher}{Tata McGraw-Hill Education}.

\bibitemdeclare{article}{samal_comparative_2018}
\bibitem{samal_comparative_2018}
\bibinfo{author}{Areejit \surnamestart Samal\surnameend},
  \bibinfo{author}{R.~P. \surnamestart Sreejith\surnameend},
  \bibinfo{author}{Jiao \surnamestart Gu\surnameend}, \bibinfo{author}{Shiping
  \surnamestart Liu\surnameend}, \bibinfo{author}{Emil \surnamestart
  Saucan\surnameend} \& \bibinfo{author}{Jürgen \surnamestart
  Jost\surnameend}: \emph{\bibinfo{title}{Comparative analysis of two
  discretizations of Ricci curvature for complex networks}}
  \bibinfo{volume}{8}(\bibinfo{number}{1}).
\newblock \doi{10.1038/s41598-018-27001-3}.
\newblock \urlprefix\url{https://dx.doi.org/10.1038/s41598-018-27001-3}.
\newblock \bibinfo{note}{Publisher: Springer Science and Business Media {LLC}}.

\bibitemdeclare{article}{saucan_discrete_2019}
\bibitem{saucan_discrete_2019}
\bibinfo{author}{Emil \surnamestart Saucan\surnameend}, \bibinfo{author}{R.P.
  \surnamestart Sreejith\surnameend}, \bibinfo{author}{R.P. \surnamestart
  Vivek-Ananth\surnameend}, \bibinfo{author}{Jürgen \surnamestart
  Jost\surnameend} \& \bibinfo{author}{Areejit \surnamestart Samal\surnameend}:
  \emph{\bibinfo{title}{Discrete Ricci curvatures for directed networks}}
  \bibinfo{volume}{118}, pp. \bibinfo{pages}{347--360}.
\newblock \doi{10.1016/j.chaos.2018.11.031}.
\newblock
  \urlprefix\url{https://linkinghub.elsevier.com/retrieve/pii/S096007791831035X}.

\bibitemdeclare{article}{shanker_defining_2007}
\bibitem{shanker_defining_2007}
\bibinfo{author}{O.~\surnamestart Shanker\surnameend}:
  \emph{\bibinfo{title}{Defining Dimension of a Complex Network}}
  \bibinfo{volume}{21}, pp. \bibinfo{pages}{321--326}.
\newblock \doi{10.1142/S0217984907012773}.

\bibitemdeclare{article}{sorkin_forks_1997}
\bibitem{sorkin_forks_1997}
\bibinfo{author}{R.~\surnamestart Sorkin\surnameend}:
  \emph{\bibinfo{title}{Forks in the road, on the way to quantum gravity}}.
\newblock \doi{10.1007/BF02435709}.

\bibitemdeclare{article}{sreejith_forman_2016}
\bibitem{sreejith_forman_2016}
\bibinfo{author}{R.~P. \surnamestart Sreejith\surnameend},
  \bibinfo{author}{Karthikeyan \surnamestart Mohanraj\surnameend},
  \bibinfo{author}{Jürgen \surnamestart Jost\surnameend},
  \bibinfo{author}{Emil \surnamestart Saucan\surnameend} \&
  \bibinfo{author}{Areejit \surnamestart Samal\surnameend}:
  \emph{\bibinfo{title}{Forman curvature for complex networks}}
  \bibinfo{volume}{2016}(\bibinfo{number}{6}), p. \bibinfo{pages}{063206}.
\newblock \doi{10.1088/1742-5468/2016/06/063206}.
\newblock \urlprefix\url{https://doi.org/10.1088/1742-5468/2016/06/063206}.
\newblock \bibinfo{note}{Publisher: {IOP} Publishing}.

\bibitemdeclare{article}{surya2019causal}
\bibitem{surya2019causal}
\bibinfo{author}{Sumati \surnamestart Surya\surnameend} (\bibinfo{year}{2019}):
  \emph{\bibinfo{title}{The causal set approach to quantum gravity}}.
\newblock {\sl \bibinfo{journal}{Living Reviews in Relativity}}
  \bibinfo{volume}{22}(\bibinfo{number}{1}), pp. \bibinfo{pages}{1--75}.

\bibitemdeclare{article}{tee2021enhanced}
\bibitem{tee2021enhanced}
\bibinfo{author}{Philip \surnamestart Tee\surnameend} \&
  \bibinfo{author}{CA~\surnamestart Trugenberger\surnameend}
  (\bibinfo{year}{2021}): \emph{\bibinfo{title}{Enhanced Forman curvature and
  its relation to Ollivier curvature}}.
\newblock {\sl \bibinfo{journal}{EPL (Europhysics Letters)}}
  \bibinfo{volume}{133}(\bibinfo{number}{6}), p. \bibinfo{pages}{60006}.

\bibitemdeclare{article}{trugenberger_combinatorial_2017}
\bibitem{trugenberger_combinatorial_2017}
\bibinfo{author}{Carlo~A. \surnamestart Trugenberger\surnameend}:
  \emph{\bibinfo{title}{Combinatorial Quantum Gravity: Geometry from Random
  Bits}} \bibinfo{volume}{2017}(\bibinfo{number}{9}), p.~\bibinfo{pages}{45}.
\newblock \doi{10.1007/JHEP09(2017)045}.
\newblock \urlprefix\url{http://arxiv.org/abs/1610.05934}.

\bibitemdeclare{article}{trugenberger_random_2016}
\bibitem{trugenberger_random_2016}
\bibinfo{author}{Carlo~A. \surnamestart Trugenberger\surnameend}:
  \emph{\bibinfo{title}{Random holographic ``large worlds'' with emergent
  dimensions}} \bibinfo{volume}{94}(\bibinfo{number}{5}), p.
  \bibinfo{pages}{052305}.
\newblock \doi{10.1103/PhysRevE.94.052305}.
\newblock \urlprefix\url{https://link.aps.org/doi/10.1103/PhysRevE.94.052305}.
\newblock \bibinfo{note}{Publisher: American Physical Society}.

\bibitemdeclare{article}{trugenberger2021emergent}
\bibitem{trugenberger2021emergent}
\bibinfo{author}{Carlo~A \surnamestart Trugenberger\surnameend}
  (\bibinfo{year}{2021}): \emph{\bibinfo{title}{Emergent time, cosmological
  constant and boundary dimension at infinity in combinatorial quantum
  gravity}}.
\newblock {\sl \bibinfo{journal}{arXiv preprint arXiv:2112.03778}}.

\bibitemdeclare{book}{winter2008curvature}
\bibitem{winter2008curvature}
\bibinfo{author}{Steffen \surnamestart Winter\surnameend}
  (\bibinfo{year}{2008}): \emph{\bibinfo{title}{Curvature measures and
  fractals}}.
\newblock \bibinfo{volume}{453}, \bibinfo{publisher}{Institute of Mathematics,
  Polish Academy of Sciences}.

\bibitemdeclare{article}{winter2013fractal}
\bibitem{winter2013fractal}
\bibinfo{author}{Steffen \surnamestart Winter\surnameend} \&
  \bibinfo{author}{Martina \surnamestart Z{\"a}hle\surnameend}
  (\bibinfo{year}{2013}): \emph{\bibinfo{title}{Fractal curvature measures of
  self-similar sets}}.
\newblock {\sl \bibinfo{journal}{Advances in Geometry}}
  \bibinfo{volume}{13}(\bibinfo{number}{2}), pp. \bibinfo{pages}{229--244}.

\bibitemdeclare{book}{Wolfram2002a}
\bibitem{Wolfram2002a}
\bibinfo{author}{Stephen \surnamestart Wolfram\surnameend}
  (\bibinfo{year}{2002}): \emph{\bibinfo{title}{A new kind of science}}.
\newblock \bibinfo{publisher}{Wolfram Media}, \bibinfo{address}{USA}.

\bibitemdeclare{article}{Wolfram2020}
\bibitem{Wolfram2020}
\bibinfo{author}{Stephen \surnamestart Wolfram\surnameend}
  (\bibinfo{year}{2020}): \emph{\bibinfo{title}{A Class of Models with the
  Potential to Represent Fundamental Physics}}.
\newblock {\sl \bibinfo{journal}{Complex Systems}}
  \bibinfo{volume}{29}(\bibinfo{number}{2}),
  \doi{10.25088/complexsystems.29.2.107}.
\newblock
  \urlprefix\url{https://www.complex-systems.com/abstracts/v29_i02_a01/}.

\bibitemdeclare{article}{wu2015emergent}
\bibitem{wu2015emergent}
\bibinfo{author}{Zhihao \surnamestart Wu\surnameend}, \bibinfo{author}{Giulia
  \surnamestart Menichetti\surnameend}, \bibinfo{author}{Christoph
  \surnamestart Rahmede\surnameend} \& \bibinfo{author}{Ginestra \surnamestart
  Bianconi\surnameend} (\bibinfo{year}{2015}): \emph{\bibinfo{title}{Emergent
  complex network geometry}}.
\newblock {\sl \bibinfo{journal}{Scientific reports}}
  \bibinfo{volume}{5}(\bibinfo{number}{1}), pp. \bibinfo{pages}{1--12}.

\bibitemdeclare{article}{zapata2022invitation}
\bibitem{zapata2022invitation}
\bibinfo{author}{Carlos \surnamestart Zapata-Carratala\surnameend} \&
  \bibinfo{author}{Xerxes~D \surnamestart Arsiwalla\surnameend}
  (\bibinfo{year}{2022}): \emph{\bibinfo{title}{An Invitation to Higher Arity
  Science}}.
\newblock {\sl \bibinfo{journal}{arXiv preprint arXiv:2201.09738}}.

\end{thebibliography}

\end{document}